    \documentclass[12pt]{amsart}
    \usepackage{anysize}
    \usepackage{amsthm}
\usepackage{amssymb}
\usepackage{times}

\input xy
\xyoption{v2}
\xyoption{all}

    \newtheorem{theorem}{Theorem}[subsection]
    \newtheorem{proposition}[theorem]{Proposition}
    \newtheorem{lemma}[theorem]{Lemma}
    \newtheorem{corollary}[theorem]{Corollary}
    
    \theoremstyle{definition}
    \newtheorem{example}[theorem]{Example}
    \newtheorem{definition}[theorem]{Definition}
    \newtheorem{remark}[theorem]{Remark}
    \theoremstyle{remark}
    \newtheorem{claim}[theorem]{Claim}

    \def\set{\setcounter{equation}
             {\value{theorem}}\addtocounter{theorem}{1}}
    \numberwithin{equation}{subsection}

\def\A{{\mathbb A}}
\def\cA{{\text{\cal A}}}

\def\cB{{\text{\cal B}}}

\def\cC{{\text{\cal C}}}
\def\sC{{\text{\sf C}}}

\def\sD{{\text{\sf D}}}

\def\F{{\mathbb F}}
\def\cF{{\text{\cal F}}}

\def\cI{{\text{\cal I}}}

\def\bL{{\mbox{\special L}}}

\def\L{{\mathbb L}}
\def\cM{{\text{\cal M}}}

\def\N{{\mathbb N}}

\def\fp{{\mathfrak p}}

\def\bP{{\mathbf P}}

\def\Q{{\mathbb Q}}

\def\cS{{\text{\cal S}}}
\def\cT{{\text{\cal T}}}

\def\Z{{\mathbb Z}}

\def\cU{{\text{\cal U}}}

\newcommand{\limdir}[1]{{\displaystyle{\mathop{\mbox{\rm lim}}_{#1}}}\,}
\newcommand{\colim}[1]{{\displaystyle{\mathop{\makebox{\rm colim}}_{#1}}}\,}
\newcommand{\liminv}[1]{{\displaystyle{\mathop{\mbox{\rm lim}}_{#1}}}\,}

\newcommand{\derotimes}{{\displaystyle{\mathop{\otimes}^\bL}}}

\newenvironment{pfclaim}{\noindent {\em Proof of the claim:\/}}
{{\smallskip}}
\newenvironment{acknowledgement}{\vskip .5cm
\setlength{\baselineskip}{0mm}
\noindent \footnotesize\rmfamily
{\em Acknowledgements\ \ \ }}

\def\Tr{\mbox{\rm Tr}}
\def\tr{\mbox{\rm tr}}
\def\gr{\mbox{\rm gr}}

\def\Ann{\mbox{\rm Ann}}

\def\flip{{\text{flip}}}

\def\Cone{\mbox{\rm Cone}}

\def\Spec{\mbox{\rm Spec}}

\def\Alhom{\mbox{\rm alHom}}
\def\Hom{\mbox{\rm Hom}}
\def\Hot{\mbox{\rm Hot}}
\def\Der{\mbox{\rm Der}}

\def\AlExt{\mbox{\rm alExt}}
\def\Ext{\mbox{\rm Ext}}
\def\hExt{{\mathbb{E}\mbox{\sf xt}}}
\def\Exal{\mbox{\rm Exal}}
\def\Exmon{\mbox{\rm Exmon}}

\def\bExun{{\bf{Exun}}}
\def\bExal{{\bf{Exal}}}
\def\bExmon{{\bf{Exmon}}}
\def\Tor{\mbox{\rm Tor}}
\def\End{\mbox{\rm End}}

\def\Coker{\mbox{\rm Coker}}
\def\Ker{\mbox{\rm Ker}}
\def\Img{\mbox{\rm Im}}

\def\Alg{{{\text{-}}\bf{Alg}}}
\def\Sym{{\text{Sym}}}
\def\AlgMorph{{\text{-}\bf{Alg.Morph}}}
\def\AlgMod{{\text{-}\bf{Alg.Mod}}}
\def\Aliso{{\text{-}\bf{al.Iso}}}

\def\Desc{{{\bf Desc}}}
\def\Et{{\bf{\acute{E}t}}}
\def\uEt{{\bf{u.\acute{E}t}}}
\def\wEt{{\bf{w.\acute{E}t}}}
\def\Mod{{\text{-}{\bf Mod}}}
\def\Mon{{\text{-}{\bf Mon}}}
\def\UniMod{{\text{-}{\bf Uni.Mod}}}

\def\Set{{\bf{Set}}}

\def\nil{{\mbox{\rm nil}}}

\def\one{\text{\bf 1}}

\def\ubar#1{{#1}_*}
\def\ubarold#1{\underline{#1}}
\def\bar#1{\overline{#1}}
 
\def\hat{\widehat} \def\tilde{\widetilde} 

\def\fm{\mathfrak m}
\def\fp{\mathfrak p}
\def\eps{\varepsilon}

\font\autf=cmcsc10
\font=ptmb7t
\font\cal=rsfs10
\font\prelim=pagko

\begin{document}
\title{Almost ring theory}
\author{Ofer Gabber}

\author{Lorenzo Ramero}

\thanks{A preliminary version of this manuscript was
prepared in 1997, while the second author was supported 
by the IHES}
        

\maketitle
\centerline{\prelim third release}

\tableofcontents

\section{Introduction}

Broadly speaking, the aim of this work is to describe 
``how to do ring theory'' within monoidal categories 
that arise as localisations of categories of modules 
over certain rings. A reader looking for forerunners 
of our themes would be drawn inevitably to Gabriel's 
``Des categories abeliennes'' \cite{Ga}, and might even 
conclude that Gabriel's memoir must have been the main 
instigation for the present article. In truth, the 
initial motivation is to be found elsewhere, namely 
in the want of adequately documented foundations for
the method of almost \'etale extensions that underpins
Faltings' approach to $p$-adic Hodge theory as presented 
in \cite{Fa2}. However, as is often the case with healthy 
offsprings, our subject matter has eventually resolved 
to venture beyond its original boundaries and pursue 
an autonomous existence.

In any case, we are glad to report that our paper 
remains true to its first vocation, which is to serve 
as a comprehensive reference, paving the way to deeper 
aspects of almost \'etale theory, especially to the 
difficult purity theorem of \cite{Fa2}. 
The notions of almost unramified and almost \'etale 
morphism are defined and their main properties are 
established, including the analogues of the classical 
lifting theorems over nilpotent extensions, and invariance
under Frobenius. Also, to any 
almost finitely presented almost flat morphism we attach 
an almost trace form, and we characterize almost \'etale 
morphisms in terms of this form. Finally, we study some 
cases of non-flat descent for almost \'etale maps.

Actually, our terminology is slightly different, in
that we replace usual modules and algebras by their ``almost'' 
counterparts, which live in the category of almost
modules, a localisation of the category of modules.
So for instance, instead of almost \'etale morphisms
of algebras, we have \'etale morphisms of almost algebras.

The categories of almost modules (or almost algebras) and
of usual modules (resp. algebras) are linked in manifold ways.
First of all we have of course the localisation functor.
Then, as it had already been observed by Gabriel, there
is a right adjoint to localisation. Furthermore, we show that
there is a {\em left\/} adjoint as well, that, to our knowledge, 
has not been exploited before, in spite of its several 
useful qualities which will establish it quickly as 
one of our main tools. The ensemble of localisation 
and right, left adjoints exhibits some remarkable
exactness properties, that are typically associated to 
open imbeddings of topoi, all of which seems to suggest the
existence of some deeper geometric structure, still 
to be unearthed. We may have encountered here an instance
of a general principle, apparently evoked first by Deligne,
according to which one should try to do algebraic geometry 
on arbitrary abelian tensor categories (notice though, that 
our categories are more general than the tannakian categories 
of \cite{DeM}).

A large part of the paper is devoted to the construction and 
study of the almost version of Illusie's cotangent complex, 
on which we base our deformation theory for almost algebras.
Faltings' original method was based on Hochschild cohomology 
rather than the cotangent complex. While Faltings' approach 
has the advantage of being more explicit and elementary, it 
also has the drawback of involving a very large number of 
long and tedious manipulations with cocycles, and requires 
a painstaking tracking of the ``epsilon book-keeping''. 
The method pre\-sented here avoids (or at least removes from 
view!) these problems, and also leads to more general results 
(especially, we can drop all finiteness assumptions from the 
statements of the lifting theorems).

Though we have strived throughout for the widest generality,
in a few places one could have gone even further : for instance 
it would have been possible to globalise all definitions and 
most results to arbitrary schemes. However, the extension to 
schemes is completely straightforward, and in practice seems 
to be scarcely useful. Similarly, there is currently not much 
incentive to study a notion of ``almost smooth morphism''.

\section{Homological theory}

\subsection{Some ring-theoretic preliminaries}\label{sec_ring.prel}

Unless otherwise stated, every ring is commutative with unit.
Our basic setup consists of a fixed base ring $V$ containing 
an ideal $\fm$ such that $\fm^2=\fm$. Starting from section 
\ref{sec_hot}, we will also assume that 
$\tilde\fm=\fm\otimes_V\fm$ is a flat $V$-module. 

\begin{example}\label{ex_rings}
i) The main example is given by a non-discrete valuation
ring $V$ with valuation $\nu:V-\{0\}\to\Gamma$ of rank one (where $\Gamma$
is the totally ordered abelian group of values of $\nu$). Then we can take 
$\fm=\{0\}\cup\{x\in V-\{0\}~|~\nu(x)>\nu(1)\}$.

ii) Suppose that $\fm=V$. This is the ``classical limit''.  In this 
case almost ring theory reduces to usual ring theory. Thus, all the 
discussion that follows specialises to, and sometimes gives 
alternative proofs for, statements about rings and their modules.
\end{example}

We define a uniform structure on the set $\cI$ of ideals of $V$ 
as follows. For every finitely generated ideal $\fm_0\subset\fm$ 
the subset of $\cI\times\cI$ given by 
$\{(I,J)~|~\fm_0\cdot J\subset I~\text{and}~\fm_0\cdot I\subset J\}$
is an entourage for the uniform structure, and the subsets of this
kind form a fundamental system of entourages (cp. 
\cite{Bou} Ch.II \S 1). The uniform structure induces a topology
on $\cI$ and moreover the notion of convergent (resp. Cauchy) 
sequence of ideals is well defined. We will
also need a notion of ``Cauchy product'' : let 
$\prod^\infty_{n=0}I_n$ be a formal infinite product of ideals.
We say that the formal product {\em satisfies the Cauchy condition\/}
(or briefly : {\em is a Cauchy product\/})
if, for every neighborhood $\cU$ of $V\in\cI$ there exists $n_0\ge 0$
such that $\prod_{m=n}^{n+p}I_m\in\cU$ for all $n\ge n_0$ and all
$p\ge 0$.

\begin{remark}\label{rem_unif} Suppose that 
$J_0\subset J_1\subset...$ is an increasing infinite sequence 
of ideals of $I$ such that $\limdir{k\to\infty}J_k=V$ (convergence 
for the above uniform structure on $\cI$). Then one checks easily 
that $\bigcup_{k=0}^\infty\fm\cdot J_k=\fm$.
\end{remark}

Let $M$ be a given $V$-module. We say that $M$ is {\em almost zero\/} 
if $\fm\cdot M=0$. A map $\phi$ of $V$-modules is an {\em almost 
isomorphism\/} if both $\Ker(\phi)$ and $\Coker(\phi)$ are almost zero 
$V$-modules.

\begin{remark}\label{rem_almost.zero}
(i) It is easy to check that a $V$-module $M$ is almost zero if and
only if $\fm\otimes_VM=0$. Similarly, a map $M\to N$ of $V$-modules 
is an almost isomorphism if and only if the induced map 
$\tilde\fm\otimes_VM\to\tilde\fm\otimes_VN$ is an isomorphism. 
Notice also that, if $\fm$ is flat, then $\fm\simeq\tilde\fm$.

(ii) Let $V\to W$ be a ring homomorphism. For a $V$-module $M$
set $M_W=W\otimes_VM$. We have an exact sequence 
\set\begin{equation}\label{eq_economy}
0\to K\to\fm_W\to\fm\cdot W\to 0
\end{equation} 
where $K=\Tor^V_1(V/\fm,W)$ is an almost zero $W$-module. By (i) 
it follows that $\fm\otimes_VK\simeq(\fm\cdot W)\otimes_WK\simeq 0$. 
Then, applying $\fm_W\otimes_W-$ and $-\otimes_W(\fm\cdot W)$
to \eqref{eq_economy} we derive 
$$\fm_W\otimes_W\fm_W\simeq\fm_W\otimes_W(\fm\cdot W)\simeq
(\fm\cdot W)\otimes_W(\fm\cdot W)$$ 
{\em i.e.\/} $\tilde\fm_W\simeq(\fm\cdot W)\tilde{}$. In particular, 
if $\tilde\fm$ is a flat $V$-module, then $\tilde\fm_W$ is a flat 
$W$-module. This means that our basic assumptions on the pair 
$(V,\fm)$ are stable under arbitrary base extension. Notice that 
the flatness of $\fm$ does not imply the flatness of $\fm\cdot W$. 
This partly explains why we insist that $\tilde\fm$, rather than 
$\fm$, be flat. 
\end{remark}

Before moving on, we want to analyze in some detail how
our basic assumptions relate to certain other natural 
conditions that can be postulated on the pair $(V,\fm)$.
Indeed, let us consider the following two hypotheses :
\smallskip

\noindent ({\bf A})\ \  $\fm=\fm^2$ and $\fm$ is a filtered 
colimit of principal ideals.
\smallskip

\noindent ({\bf B})\ \  $\fm=\fm^2$ and, for all integers $k>1$, 
the $k$-th powers of elements of $\fm$ generate $\fm$.
\smallskip

Clearly ({\bf A}) implies ({\bf B}). Less obvious is
the following result.
\begin{proposition}\label{prop_less.obv}
(i) ({\bf A}) implies that $\tilde\fm$ is flat.

(ii) If $\tilde\fm$ is flat then ({\bf B}) holds. 
\end{proposition}
\begin{proof} Suppose that ({\bf A}) holds, so that 
$\fm=\colim{\alpha\in I}Vx_\alpha$, where $I$ is a directed set
parametrizing elements $x_\alpha\in\fm$ (and 
$\alpha\leq\beta\Leftrightarrow Vx_\alpha\subset Vx_\beta$).
For any $\alpha\in I$ we have natural isomorphisms
\set\begin{equation}\label{eq_pi.alpha}
Vx_\alpha\simeq V/\Ann_V(x_\alpha)\simeq
(Vx_\alpha)\otimes_V(Vx_\alpha).
\end{equation}
For $\alpha\leq\beta$, let 
$j_{\alpha\beta}:Vx_\alpha\hookrightarrow Vx_\beta$ be the imbedding;
we have a commutative diagram
$$\xymatrix{
V \ar[rr]^{\mu_{z^2}} \ar[d]_{\pi_\alpha} & & V \ar[d]^{\pi_\beta} \\
(Vx_\alpha)\otimes_V(Vx_\alpha) 
\ar[rr]^{j_{\alpha\beta}\otimes j_{\alpha\beta}} & &
(Vx_\beta)\otimes_V(Vx_\beta)
}$$
where $z\in V$ is such that $x_\alpha=z\cdot x_\beta$, $\mu_{z^2}$
is multiplication by $z^2$ and $\pi_\alpha$ is the projection
induced by \eqref{eq_pi.alpha} (and similarly for $\pi_\beta$).
Since $\fm=\fm^2$, for all $\alpha\in I$ we can find $\beta$
such that $x_\alpha$ is a multiple of $x_\beta^2$. Say 
$x_\alpha=y\cdot x^2_\beta$; then we can take $z=y\cdot x_\beta$,
so $z^2$ is a multiple of $x_\alpha$ and in the above diagram
$\Ker(\pi_\alpha)\subset\Ker(\mu_{z^2})$. Hence one can define
a map $\lambda_{\alpha\beta}:(Vx_\alpha)\otimes_V(Vx_\alpha)\to V$
such that $\pi_\beta\circ\lambda_{\alpha\beta}=
j_{\alpha\beta}\otimes j_{\alpha\beta}$ and 
$\lambda_{\alpha\beta}\circ\pi_\alpha=\mu_{z^2}$. It now
follows that for every $V$-module $N$, the induced morphism
$\Tor_1^V(N,(Vx_\alpha)\otimes_V(Vx_\alpha))\to
\Tor_1^V(N,(Vx_\beta)\otimes_V(Vx_\beta))$
is the zero map. Taking the colimit we derive that $\tilde\fm$
is flat. This shows (i).
In order to show (ii) we consider, for any prime number $p$, 
the following condition 
\medskip

\noindent($*_p$)\qquad $\fm/p\cdot\fm$ is generated (as a $V$-module)
by the $p$-th powers of its elements.
\medskip

Clearly ({\bf B}) implies ($*_p$) for all $p$. In fact we have :
\begin{claim}
({\bf B}) holds if and only if ($*_p$) holds for every prime $p$. 
\end{claim}
\begin{pfclaim} Suppose that ($*_p$) holds for every prime $p$.
The polarization identity 
$$k!\cdot x_1\cdot x_2\cdot...\cdot x_k=
\sum_{I\subset\{1,2,...,k\}}(-1)^{k-|I|}\cdot
\left(\sum_{i\in I}x_i\right)^k$$
shows that if $N=\sum_{x\in\fm}Vx^k$ then $k!\cdot\fm\subset N$.
To prove that $N=\fm$ it then suffices to show that for every
prime $p$ dividing $k!$ we have $\fm=p\cdot\fm+N$. Let 
$\phi:V/p\cdot V\to V/p\cdot V$ be the Frobenius ($x\mapsto x^p$); 
we can denote by $(V/p\cdot V)^\phi$ the ring $V/p\cdot V$ 
seen as a $V/p\cdot V$-algebra via the homomorphism $\phi$. 
Also set $\phi^*M=M\otimes_{V/p\cdot V}(V/p\cdot V)^\phi$ 
for a $V/p\cdot V$-module $M$. Then the map 
$\phi^*(\fm/p\cdot\fm)\to(\fm/p\cdot\fm)$ (defined by raising 
to $p$-th power) is surjective by ($*_p$). Hence so is
$(\phi^r)^*(\fm/p\cdot\fm)\to(\fm/p\cdot\fm)$ for every $r>0$,
which says that $\fm=p\cdot\fm+N$ when $k=p^r$, hence for every
$k$.\end{pfclaim}

Next recall (see \cite{SGA4} Exp. XVII 5.5.2) that, if $M$
is a $V$-module, the module of symmetric tensors $\text{TS}^k(M)$
is defined as $(\otimes_V^kM)^{S_k}$, the invariants under the
natural action of the symmetric group $S_k$ on $\otimes_V^kM$.
We have a natural map $\Gamma^k(M)\to\text{TS}^k(M)$ that is an
isomorphism when $M$ is flat (see {\em loc. cit.\/} 5.5.2.5;
here $\Gamma^k$ denotes the $k$-th graded piece of the divided 
power algebra).

\begin{claim} The group $S_k$ acts trivially on 
$\otimes_V^k\fm$ and the map $\tilde\fm\otimes_V\fm\to\tilde\fm$ 
($x\otimes y\otimes z\mapsto x\otimes yz$) is an isomorphism.
\end{claim}
\begin{pfclaim} The first statement is reduced to the case
of transpositions and to $k=2$. There we can compute :
$x\otimes yz=xy\otimes z=y\otimes xz=yz\otimes x$.
For the second statement note that the imbedding 
$\fm\hookrightarrow V$ is an almost isomorphism, 
and apply remark \ref{rem_almost.zero}(i).
\end{pfclaim}

Suppose now that $\tilde\fm$ is flat and pick a prime $p$. 
Then $S_p$ acts trivially on $\otimes^p_V\tilde\fm$. Hence 
\set\begin{equation}\label{eq_iso.Gamma}
\Gamma^p(\tilde\fm)\simeq\otimes^p_V\tilde\fm\simeq\tilde\fm.
\end{equation}
But $\Gamma^p(\tilde\fm)$ is spanned as a $V$-module by
the products $\gamma_{i_1}(x_1)\cdot...\cdot\gamma_{i_k}(x_k)$
(where $x_i\in\tilde\fm$ and $\sum_ji_j=p$). Under the isomorphism
\eqref{eq_iso.Gamma} these elements map to 
${p\choose i_1,...,i_k}\cdot x_1^{i_1}\cdot...\cdot x_k^{i_k}$;
but such an element vanishes in $\tilde\fm/p\cdot\tilde\fm$ 
unless $i_k=p$ for some $k$. Therefore $\tilde\fm/p\cdot\tilde\fm$
is generated by $p$-th powers, so the same is true for $\fm/p\cdot\fm$,
and by the above, ({\bf B}) holds, which shows (ii).
\end{proof}

\begin{proposition}\label{prop_countably.pres}
Suppose that $\fm$ is countably generated as a $V$-module. 
Then we have :

i) $\tilde\fm$ is countably presented as a $V$-module;

ii) if $\tilde\fm$ is a flat $V$-module, then it is of 
homological dimension $\leq 1$.
\end{proposition}
\begin{proof} Let $(\eps_i)_{i\in I}$ be a countable 
generating family of $\fm$. Then $\eps_i\otimes\eps_j$ 
generate $\tilde\fm$ and 
$\eps_i\cdot\eps_j\cdot(\eps_k\otimes\eps_l)=
\eps_k\cdot\eps_l\cdot(\eps_i\otimes\eps_j)$ for all 
$i,j,k,l\in I$. For every $i\in I$, we can write 
$\eps_i=\sum_jx_{ij}\eps_j$, for certain $x_{ij}\in\fm$. 
Let $F$ be the $V$-module defined by
generators $(e_{ij})_{i,j\in I}$, subject to the 
relations:
$$\eps_i\cdot\eps_j\cdot e_{kl}=
\eps_k\cdot\eps_l\cdot e_{ij}\qquad 
e_{ik}=\sum_jx_{ij}e_{jk}\qquad\text{for all $i,j,k,l\in I$.}$$
We get an epimorphism $\pi:F\to\tilde\fm$ by 
$e_{ij}\mapsto\eps_i\otimes\eps_j$. The relations 
imply that, if $x=\sum_{k,l}y_{kl}\in\Ker(\pi)$, then 
$\eps_i\cdot\eps_j\cdot x=0$, so $\fm\cdot\Ker(\pi)=0$.
Whence $\fm\otimes_V\Ker(\pi)=0$ and $\one_\fm\otimes_V\pi$
is an isomorphism. We consider the diagram 
$$\xymatrix{\fm\otimes_VF \ar[r]^-{\sim} \ar[d]_\phi & 
\fm\otimes_V\tilde\fm \ar[d]^\psi \\ F \ar[r]^-\pi & 
\tilde\fm}$$
where $\phi$ and $\psi$ are induced by scalar multiplication.
We already know that $\psi$ is an isomorphism, and since 
$F=\fm\cdot F$, we see that $\phi$ is an epimorphism, 
so $\pi$ is an isomorphism, which shows (i). Now (ii)
follows from (i), since it is well-known that a flat
countably presented module is of homological dimension 
$\leq 1$ (see \cite{La} (Ch.I, Th.3.2) and the discussion 
in \cite{Mit} pp.49-50).
\end{proof}

\subsection{Almost categories}\label{sec_alm.cat}
If $\cC$ is a category, and $X,Y$ two objects of $\cC$, we will
usually denote by $\Hom_\cC(X,Y)$ the set of morphisms in $\cC$
from $X$ to $Y$ and by $\one_X$ the identity morphism of $X$. 
Moreover we denote by $\cC^o$ the opposite category of $\cC$ 
and by $s.\cC$ the category of simplicial objects over $\cC$, 
that is, functors $\Delta^o\to\cC$, where $\Delta$ is the category 
whose objects are the ordered sets $[n]=\{0,...,n\}$ for each integer 
$n\ge 0$ and where a morphism $\phi:[p]\to[q]$ is a non-decreasing
map. A morphism $f:X\to Y$ in $s.\cC$ is a sequence of morphisms
$f_{[n]}:X[n]\to Y[n]$, $n\ge 0$ such that the obvious diagrams commute. 
We can imbed $\cC$ in $s.\cC$ by sending each object $X$ to the 
``constant'' object $s.X$ such that $s.X[n]=X$ for all $n\ge 0$
and $s.X[\phi]=\one_X$ for all morphisms $\phi$ in $\Delta$. 

If $\cC$ is an abelian category, $\sD(\cC)$ will denote the 
derived category of $\cC$. As usual we have also the full subcategories 
$\sD^+(\cC),\sD^-(\cC)$ of complexes of objects of $\cC$ that are 
exact for sufficiently large negative (resp. positive) degree.
If $R$ is a ring, the category of $R$-modules (resp. $R$-algebras) will 
be denoted by $R\Mod$ (resp. $R\Alg$). Most of the times we will write 
$\Hom_R(M,N)$ instead of $\Hom_{R\Mod}(M,N)$. 

We denote by $\Set$ the category of sets. 
The symbol $\N$ denotes the set of non-negative 
integers; in particular $0\in\N$.
\medskip

The full subcategory $\Sigma$ of $V\Mod$ consisting of all $V$-modules
that are almost isomorphic to $0$ is clearly a Serre subcategory and
hence we can form the quotient category $V\Mod/\Sigma$.
There is a localization functor 
$$V\Mod\to V\Mod/\Sigma \qquad M\mapsto M^a$$ 
that takes a $V$-module $M$ to the same module, seen as an object 
of $V\Mod/\Sigma$. In particular, we have the object $V^a$
associated to $V$; it seems therefore natural to use the notation 
$V^a\Mod$ for the category $V\Mod/\Sigma$, and an object of 
$V^a\Mod$ will be indifferently referred to as ``a $V^a$-module''
or ``an almost $V$-module''. In case we need to stress the
dependance on the ideal $\fm$, we can write $(V,\fm)^a\Mod$.

Since the almost isomorphisms form a multiplicative system (see 
{\em e.g.\/} \cite{We} Exerc.10.3.2), it is possible to describe 
the morphisms in $V^a\Mod$ via a calculus of fractions, as follows. 
Let $V\Aliso$ be the category that has the same objects as $V\Mod$,
but such that $\Hom_{V\Aliso}(M,N)$ consists of all almost
isomorphisms $M\to N$. If $M$ is any object of $V\Aliso$ we write 
$(V\Aliso/M)$ for the category of objects of $V\Aliso$ over $M$ 
({\em i.e.\/} morphisms $\phi:X\to M$). If $\phi_i:X_i\to M$ $(i=1,2)$ are
two objects of $(V\Aliso/M)$ then $\Hom_{(V\Aliso/M)}(\phi_1,\phi_2)$
consists of all morphisms $\psi:X_1\to X_2$ in $V\Aliso$ such that
$\phi_1=\phi_2\circ\psi$. For any two $V$-modules $M,N$ we define a 
functor $\cF_N:(V\Aliso/M)^o\to V\Mod$ by
associating to an object $\phi:P\to M$ the $V$-module $\Hom_V(P,N)$ 
and to a morphism $\alpha:P\to Q$ the map 
$\Hom_V(Q,N)\to\Hom_V(P,N)~:~\beta\mapsto\beta\circ\alpha$.
Then we have
\set\begin{equation}\label{eq_alhom}
\Hom_{V^a\Mod}(M^a,N^a)=\colim{(V\Aliso/M)^o}\cF_N.
\end{equation}
However, formula \eqref{eq_alhom} can be simplified considerably,
by remarking that, for any $V$-module $M$, the natural morphism 
$\tilde\fm\otimes_VM\to M$ is an initial object of $(V\Aliso/M)$.
Indeed, let $\phi:N\to M$ be an almost isomorphism; the diagram
$$\xymatrix{
\tilde\fm\otimes_VN \ar[r]^\sim \ar[d] & 
\tilde\fm\otimes_VM \ar[d] \\
N \ar[r]^\phi & M
}$$
(cp. remark \ref{rem_almost.zero}(i)) allows one to define a 
morphism $\psi:\tilde\fm\otimes_VM\to N$ over $M$. We need to
show that $\psi$ is unique. But if   
$\psi_1,\psi_2:\tilde\fm\otimes_VM\to N$ are two maps over
$M$, then $\Img(\psi_1-\psi_2)\subset\Ker(\phi)$ is almost 
zero, hence $\Img(\psi_1-\psi_2)=0$, since 
$\tilde\fm\otimes_VM=\fm\cdot(\tilde\fm\otimes_VM)$. 
Consequently, \eqref{eq_alhom} boils down to
\set\begin{equation}\label{eq_alm.morph}
\Hom_{V^a\Mod}(M^a,N^a)=\Hom_V(\tilde\fm\otimes_VM,N).
\end{equation}
In particular $\Hom_{V^a\Mod}(M,N)$ has a natural structure of $V$-module
for any two $V^a$-modules $M,N$, {\em i.e.\/} $\Hom_{V^a\Mod}(-,-)$
is a bifunctor that takes values in the category $V\Mod$.

One checks easily (for instance using \eqref{eq_alm.morph}) that 
the usual tensor product induces a bifunctor $-\otimes_V-$ on almost
$V$-modules, which, in the jargon of \cite{DeM} makes of $V^a\Mod$ an 
{\em abelian tensor category\/}. Then an {\em almost $V$-algebra} 
is just a commutative unitary monoid in the tensor category $V^a\Mod$. 
Let us recall what this means. Quite generally, let $(\cC,\otimes,U)$ 
be any abelian tensor category, so that $\otimes:\cC\times\cC\to\cC$ 
is a biadditive functor, $U$ is the identity object
of $\cC$ (see \cite{DeM} p.105) and for any two objects $M$ 
and $N$ in $\cC$ we have a ``commutativity constraint'' ({\em i.e.\/}
a functorial isomorphism $\eta_{M|N}:M\otimes N\to N\otimes M$ 
that ``switches the two factors'') and a functorial isomorphism
$\nu_M:U\otimes M\to M$. Then a $\cC$-monoid $A$ is an
object of $\cC$ endowed with a morphism $\mu_A:A\otimes A\to A$ 
(the ``multiplication'' of $A$) satisfying the associativity condition
$$\mu_A\circ(\one_A\otimes\mu_A)=\mu_A\circ(\mu_A\otimes\one_A).$$
We say that $A$ is {\em unitary\/} if additionally $A$ is endowed 
with a ``unit morphism'' $\ubarold 1_A:U\to A$ satisfying the 
(left and right) unit property : 
$$\mu_A\circ(\ubarold 1_A\otimes\one_A)=\nu_A\qquad
\mu_A\circ(\ubarold 1_A\otimes\one_A)\circ\eta_{A|U}=
\mu_A\circ(\one_A\otimes \ubarold 1_A).$$
Finally $A$ is {\em commutative\/} if $\mu_A=\mu_A\circ\eta_{A|A}$
(to be rigorous, in all of the above one should indicate the 
associativity constraints, which we have omitted : see \cite{DeM}).
A commutative unitary monoid will also be simply called
an {\em algebra\/}. 
With the morphisms defined in the obvious way, the $\cC$-monoids 
form a category; furthermore, given a $\cC$-monoid $A$, a 
{\em left $A$-module} is an object $M$ of $\cC$ 
endowed with a morphism $\sigma_M:A\otimes M\to M$ 
such that $\sigma_M\circ(\one_A\otimes\sigma_M)=
\sigma_M\circ(\mu_A\otimes\one_M)$.
Similarly one defines right $A$-modules and $A$-bimodules. 
In the case of bimodules we have left and right morphisms
$\sigma_{M,l}:A\otimes M\to M$, $\sigma_{M,r}:M\otimes A\to M$
and one imposes that they ``commute'', {\em i.e.\/} that
$$\sigma_{M,r}\circ(\sigma_{M,l}\otimes\one_A)=
\sigma_{M,l}\circ(\one_A\otimes\sigma_{M,r}).$$
Clearly the (left resp. right) $A$-modules (and the $A$-bimodules)
form an additive category with {\em $A$-linear morphisms\/} 
defined as one expects.
One defines the notion of a submodule as an equivalence class of
monomorphisms $N\to M$ such that the composition 
$A\otimes N\to A\otimes M\to M$ factors through $N$. Now, if 
$f:M\to N$ is a morphism of left $A$-modules, then $\Ker(f)$ 
exists in the underlying abelian category $\cC$ and one checks 
easily that it has a unique structure of left $A$-module which 
makes it a submodule of $M$. {\em If moreover $\otimes$ is right 
exact\/} when either argument is fixed, then also $\Coker(f)$ 
has a unique $A$-module structure for which $N\to\Coker(f)$ is 
$A$-linear. In this case the category of left $A$-modules is 
abelian. Similarly, if $A$ is a unitary $\cC$-monoid, then one 
defines the notion of {\em unitary\/} left $A$-module by 
requiring that $\sigma_M\circ(\ubarold 1_A\otimes\one_M)=\nu_M$ 
and these form an abelian category when $\otimes$ is right exact.

Specialising to our case we obtain the category $V^a\Alg$ of 
almost $V$-algebras and, for every almost $V$-algebra $A$, the 
category $A\Mod$ of unitary left $A$-modules. Clearly the 
localization functor restricts to 
a functor $V\Alg\to V^a\Alg$ and for any $V$-algebra $R$ 
we have a localization functor $R\Mod\to R^a\Mod$.

Next, if $A$ is an almost $V$-algebra, we can define the category
$A\Alg$ of $A$-algebras. It consists of all the morphisms 
$A\to B$ of almost $V$-algebras. 

Let again $(\cC,\otimes,U)$ be any abelian tensor category.
By \cite{DeM} p.119, the endomorphism ring $\End_\cC(U)$ of $U$
is commutative. For any object $M$ of $\cC$, denote 
$M_*=\Hom_\cC(U,M)$; then $M\mapsto M_*$ defines a functor
$\cC\to\End_\cC(U)\Mod$. Moreover, if $A$ is a $\cC$-monoid,
$A_*$ is an associative $\End_\cC(U)$-algebra, with multiplication 
given as follows. For $a,b\in A_*$ let 
$a\cdot b=\mu_A\circ(a\otimes b)\circ\nu_U^{-1}$.
Similarly, if $M$ is an $A$-module, $M_*$ is an $A_*$-module
in a natural way, and in this way we obtain a functor 
from $A$-modules and $A$-linear morphisms to $A_*$-modules
and $A_*$-linear maps. Using \cite{DeM} (Prop. 1.3), one 
can also check that $\End_\cC(U)=U_*$ as $\End_\cC(U)$-algebras,
where $U$ is viewed as a $\cC$-monoid using $\nu_U$. 

All this applies especially to our categories of almost
modules and almost algebras : in this case we call $M\mapsto M_*$
the {\em functor of almost elements\/}. So, if $M$ is an 
almost module, an almost element of $M$ is just an honest 
element of $M_*$. Using \eqref{eq_alm.morph} one can show 
easily that for every $V$-module $M$ the natural map 
$M\to(M^a)_*$ is an almost isomorphism.

Let $A$ be an almost $V$-algebra. For any two $A$-modules $M,N$, 
the set $\Hom_{A\Mod}(M,N)$ has a natural structure of 
$\ubar A$-module and we obtain an internal Hom functor by letting
$$\Alhom_A(M,N)=\Hom_{A\Mod}(M,N)^a.$$
This is the functor of {\em almost homomorphisms\/} from $M$ to $N$.

For any $A$-module $M$ we have also a functor 
of tensor product $M\otimes_A-$ on $A$-modules which, in view
of the following proposition \ref{prop_adjoint} can be shown to be a 
left adjoint to the functor $\Alhom_A(M,-)$. It can be defined as
$M\otimes_AN=(\ubar M\otimes_{\ubar A}\ubar N)^a$ but an appropriate 
almost version of the usual construction would also work.

With this tensor product, $A\Mod$ is an abelian tensor category 
as well, and $A\Alg$ could also be described as the category of 
($A\Mod$)-algebras. Under this equivalence, a morphism 
$\phi:A\to B$ of almost $V$-algebras becomes the unit morphism 
$\ubarold 1_B:A\to B$ of the corresponding monoid. We will 
sometimes drop the subscript and write simply $\ubarold 1$.

\begin{remark}\label{rem_base.comp} Let $V\to W$ be a map of 
base rings, $W$ taken with the extended ideal $\fm\cdot W$. 
Then $W^a$ is an almost $V$-algebra so we have defined the 
category $W^a\Mod$ using base ring $V$ and the category 
$(W,\fm\cdot W)^a\Mod$ using base $W$. One shows easily 
that they are equivalent: we have an obvious functor 
$(W,\fm\cdot W)^a\Mod\to W^a\Mod$ and an essential inverse 
is provided by $M\mapsto M_*$. 
Similar base comparison statements hold for the categories 
of almost algebras.
\end{remark} 
\begin{proposition}\label{prop_adjoint}
i) There is a natural isomorphism $A\simeq A_*^a$ of almost 
$V$-algebras.

ii) Let $R$ be any $V$-algebra. Then the functor $M\mapsto\ubar M$ 
from $R^a\Mod$ to $R\Mod$ (resp. from $R^a\Alg$ to $R\Alg$) 
is right adjoint to the localization functor $R\Mod\to R^a\Mod$ 
(resp. $R\Alg\to R^a\Alg$).

iii) The counit of the adjunction $\ubar M^a\to M$ is a natural 
isomorphism from the composition of the two functors to the 
identity functor $\one_{A\Mod}$ (resp. $\one_{A\Alg}$).
\end{proposition}

\begin{proof} (i) has already been remarked. We show (ii). 
In light of remark \ref{rem_base.comp} (applied with 
$W=R$) we can assume that $V=R$. Let $M$ be a $V$-module
and $N$ an almost $V$-module; we have natural bijections
$$\begin{array}{r@{\:\simeq\:}l}
\Hom_{V^a\Mod}(M^a,N) & \Hom_{V^a\Mod}(M^a,(N_*)^a) 
\simeq \Hom_V(\tilde\fm\otimes_VM,N_*) \\
& \Hom_V(M,\Hom_V(\tilde\fm,N_*)) 
\simeq \Hom_V(M,\Hom_{V^a\Mod}(V,(N_*)^a)) \\
& \Hom_V(M,N_*)
\end{array}$$
which proves (ii). Now (iii) follows by inspecting 
the proof of (ii), or by \cite{Ga} (ch.III Prop.3).
\end{proof}

\begin{remark} The existence of the right adjoint follows 
also directly from \cite{Ga} (chap.III \S 3 Cor.1 or chap.V \S2). 
\end{remark}

\begin{corollary} The functor $M\mapsto M_*$ from $R^a\Mod$ 
to $R\Mod$ sends injectives to injectives and injective envelopes
to injective envelopes.
\end{corollary}
\begin{proof} The functor $M\mapsto M_*$ is right adjoint to 
an exact functor, hence it preserves injectives. Now, let $J$ 
be an injective envelope of $M$; to show that $J_*$ is an injective
envelope of $M_*$, it suffices to show that $J_*$ is an essential
extension of $M_*$. However, if $N\subset J_*$ and $N\cap M_*=0$,
then $N^a\cap M=0$, hence $\fm\cdot N=0$, but $J_*$ does not 
contain $\fm$-torsion, thus $N=0$.
\end{proof}

\begin{corollary}\label{cor_coco}
The categories $A\Mod$ and $A\Alg$ are both complete and cocomplete.
\end{corollary}
\begin{proof} We recall that the categories $\ubar A\Mod$ and $\ubar A\Alg$
are both complete and cocomplete. Now let $I$ be any small indexing category 
and $M:I\to A\Mod$ be any functor. Denote by $\ubar M:I\to\ubar A\Mod$ 
the composed functor $i\mapsto\ubar{M(i)}$. We claim that
$\colim{I} M=(\colim{I}\ubar M)^a$.
The proof is an easy application of proposition \ref{prop_adjoint}(iii).
A similar argument also works for limits and for the category $A\Alg$.
\end{proof}

Note that the essential image of $M\mapsto M_*$ is closed under limits.
Next recall that the forgetful functor $\ubar A\Alg\to\Set$ (resp.
$\ubar A\Mod\to\Set$) has a left adjoint $\ubar A[-]:\Set\to\ubar A\Alg$ 
(resp. $A^{(-)}:\Set\to\ubar A\Mod$) that assigns to a set $S$ the free 
$\ubar A$-algebra $\ubar A[S]$ (resp. the free $\ubar A$-module 
$\ubar A^{(S)}$) 
generated by $S$. If $S$ is any set, it is natural to write $A[S]$
(resp. $A^{(S)}$) for the $A$-algebra $(\ubar A[S])^a$ (resp. for
the $A$-module $(\ubar A^{(S)})^a$. This yields a left adjoint,
called the {\em free $A$-algebra\/} functor $\Set\to A\Alg$ 
(resp. the {\em free $A$-module\/} functor $\Set\to A\Mod$) 
to the ``forgetful'' functor $A\Alg\to\Set$ (resp. $A\Mod\to\Set$) 
$B\mapsto\ubar B$.

Now let $R$ be any $V$-algebra; we want to construct a left adjoint 
to the localisation functor $R\Mod\to R^a\Mod$. For a given 
$R^a$-module $M$, let
\set\begin{equation}\label{eq_left.adj}
M_!=\tilde\fm\otimes_V(M_*).
\end{equation}
We have the natural map (unit of adjunction) $R\to R^a_*$, 
so that we can view $M_!$ as an $R$-module.

\begin{proposition}\label{prop_left.adj} i) The functor 
$R^a\Mod\to R\Mod$ defined by \eqref{eq_left.adj} is left 
adjoint to localisation.

ii) The unit of the adjunction $M\to M^a_!$ is a natural isomorphism
from the identity functor $\one_{R^a\Mod}$ to the composition of the 
two functors.
\end{proposition}
\begin{proof} (i) follows easily from \eqref{eq_alm.morph} 
and (ii) follows easily from (i).
\end{proof}
\begin{corollary}\label{cor_left.exact} Suppose that $\tilde\fm$ 
is a flat $V$-module. Then we have : 

i) the functor $M\mapsto M_!$ is exact;

ii) the localisation functor $R\Mod\to R^a\Mod$ sends injectives
to injectives.
\end{corollary}
\begin{proof} By proposition \ref{prop_left.adj} it follows
that $M\mapsto M_!$ is right exact. To show that it is also left 
exact when $\tilde\fm$ is a flat $V$-module, it suffices to remark 
that $M\mapsto M_*$ is left exact. Now, by (i), the functor
$M\mapsto M^a$ is right adjoint to an exact functor, so (ii) 
is clear.
\end{proof}

Next, let $B$ be any $A$-algebra. The multiplication on
$B_*$ is inherited by $B_!$, which is therefore a non-unital
ring in a natural way. We endow the $V$-module
$V\oplus B_!$ with the ring structure determined by the rule:
$(v,b)\cdot(v',b')=(v\cdot v',v\cdot b'+v'\cdot b+b\cdot b')$
for all $v,v'\in V$ and $b,b'\in B_!$. Then $V\oplus B_!$ is a 
(unital) ring. Let $\mu_\fm:\tilde\fm\to\fm$ be the map 
defined by $x\otimes y\mapsto xy$ for all $x,y\in\fm$;
we notice that the subset of all elements of the 
form $(\mu(s),-s\otimes\ubarold 1)$ (for arbitrary 
$s\in\tilde\fm$) forms an ideal $I$ of $V\oplus B_!$. 
Set $B_{!!}=(V\oplus B_!)/I$.
Thus we have a sequence of $V$-modules
\set\begin{equation}\label{eq_left.adj.algebras}
0\to\tilde\fm\to V\oplus B_!\to B_{!!}\to 0
\end{equation}
which in general is only right exact.
\begin{definition} We say that $B$ is an {\em exact
$A$-algebra\/} if the sequence \eqref{eq_left.adj.algebras}
is exact.
\end{definition}
\begin{remark}\label{rem_exact.alg} Notice that if 
$\tilde\fm\stackrel{\sim}{\to}\fm$ ({\em e.g.\/} when $\fm$ is flat), 
then all $A$-algebras are exact. In the general case, 
if $B$ is any $A$-algebra, then $V^a\times B$ is always exact.
Indeed, we have $(V^a\times B)_*\simeq V^a_*\times B_*$
and, by remark \ref{rem_almost.zero}(i), 
$\tilde\fm\otimes_VV^a_*\simeq\tilde\fm$.
\end{remark}
Clearly we have a natural isomorphism $B\simeq B_{!!}^a$.
\begin{proposition} The functor $B\mapsto B_{!!}$ is left
adjoint to the localisation functor $A_{!!}\Alg\to A\Alg$.
\end{proposition}
\begin{proof} Let $B$ be an $A$-algebra, $C$ an
$A_{!!}$-algebra and $\phi:B\to C^a$ a morphism of 
$A$-algebras. By proposition \ref{prop_left.adj} we obtain 
a natural $A_*$-linear morphism $B_!\to C$. Together with
the structure morphism $V\to C$ this yields a map 
$\tilde\phi:V\oplus B_!\to C$ which is easily seen to be 
a ring homomorphism. It is equally clear that the ideal 
$I$ defined above is mapped to zero by $\tilde\phi$, hence 
the latter factors through a map of $A_{!!}$-algebras $B_{!!}\to C$.
Conversely, such a map induces a morhism of $A$-algebras
$B\to C^a$ just by taking localisation. It is easy to check
that the two procedures are inverse to each other, which shows
the assertion.
\end{proof}

\begin{remark}\label{rem_adjoints}
The functor of almost elements commutes with arbitrary limits, 
because all right adjoints do. It does not in general commute with 
colimits, not even with arbitrary infinite direct sums.
Dually, the functors $M\mapsto M_!$ and $B\mapsto B_{!!}$ 
commute with all colimits. In particular, the latter commutes
with tensor products.
\end{remark}

\subsection{Almost homological algebra}\label{sec_homol} 
In this section we fix an almost $V$-algebra $A$ and we consider 
various constructions in the category of $A$-modules. 

\begin{remark}\label{rem_tricky}
i) Let $M_1,M_2$ be two $A$-modules. By proposition \ref{prop_adjoint} 
it is clear that a morphism $\phi:M_1\to M_2$ of $A$-modules is 
uniquely determined by the induced morphism $\ubar{M_1}\to\ubar{M_2}$. 

ii) It is a bit tricky to deal with preimages of almost elements under
morphisms: for instance, if $\phi:M_1\to M_2$ is an epimorphism (by which
we mean that $\Coker(\phi)\simeq 0$) and $m_2\in\ubar{M_2}$, then it is not 
true in general that we can find an almost element $m_1\in\ubar{M_1}$ 
such that $\ubar\phi(m_1)=m_2$. What remains true is that for arbitrary 
$\eps\in\fm$ we can find $m_1$ such that $\ubar\phi(m_1)=\eps\cdot m_2$. 
\end{remark}

The abelian category $A\Mod$ satisfies axiom (AB5) (see {\em e.g.\/}
\cite{We} (\S A.4)) and it has a generator, namely the object $A$ itself. 
It then follows by a general result that $A\Mod$ has enough injectives.
By corollary \ref{cor_coco} any inverse system $\{M_n~|~n\in\N\}$ of 
$A$-modules has an (inverse) limit $\liminv{n\in\N}M$. As usual,
we denote by $\liminv{}^1$ the right derived functor of the inverse limit
functor. Notice that \cite{We} (Cor. 3.5.4) holds in the almost case 
since axiom (AB4*) holds in $A\Mod$ (on the other hand, it is not 
clear whether \cite{We} (Lemma 3.5.3) holds under (AB4*), since the
proof uses elements).

\begin{lemma}\label{lem_invlim}
Let $\{M_n~;~\phi_n:M_n\to M_{n+1}~|~n\in\N\}$ 
(resp. $\{N_n~;~\psi_n:N_{n+1}\to N_n~|~n\in\N\}$)
be a direct (resp. inverse) system of $A$-modules and morphisms
and $\{\eps_n~|~n\in\N\}$ a sequence of ideals of $V$ converging 
to $V$ (for the uniform structure introduced in section 
\ref{sec_ring.prel}).

i) If $\eps_n\cdot M_n=0$ for all $n\in\N$ then 
$\colim{n\in\N}M_n\simeq 0$.

ii) If $\eps_n\cdot N_n=0$ for all $n\in\N$ then 
$\liminv{n\in\N}N_n\simeq 0\simeq\liminv{n\in\N}^1N_n$.

iii) If $\eps_n\cdot\Coker(\psi_n)=0$ for all $n\in\N$ and 
$\prod_{j=0}^\infty\eps_j$ is a Cauchy product, then 
$\liminv{n\in\N}^1N_n\simeq 0$.
\end{lemma}

\begin{proof} (i) and (ii) : we remark only that 
$\liminv{n\in\N}^1N_n\simeq\liminv{n\in\N}^1N_{n+p}$ for all
$p\in\N$ and leave the details to the reader.
We prove (iii). From \cite{We} (Cor. 3.5.4) it follows easily
that $(\liminv{n\in\N}^1N_{n*})^a\simeq\liminv{n\in\N}^1N_n$.
It then suffices to show that $\liminv{n\in\N}^1N_{n*}$
is almost zero. We have $\eps^2_n\cdot\Coker(\psi_{n*})=0$
and the product $\prod_{j=0}^\infty(\eps_j^2)$ is again a Cauchy 
product. Next let $N'_n=\bigcap_{p\ge 0}\Img(N_{n+p*}\to N_{n*})$. 
If $J_n=\bigcap_{p\ge 0}(\eps_n\cdot\eps_{n+1}\cdot...\cdot\eps_{n+p})^2$
then $J_n\cdot N_{n*}\subset N'_n$ and 
$\limdir{n\to\infty}J_n=V$. In view of (ii),
$\liminv{n\in\N}^1N_{n*}/N'_n$ is almost zero, hence we reduce 
to showing that $\liminv{n\in\N}^1N'_n$ is almost zero.
But  
$$J_{n+p+q}\cdot N'_n\subset\Img(N'_{n+p+q}\to N'_n)
\subset\Img(N'_{n+p}\to N'_n)$$
for all $n,p,q\in\N$. On the other hand, by remark 
\ref{rem_unif} we get 
$\bigcup_{q=0}^\infty\fm\cdot J_{n+p+q}=\fm$, hence
$\fm\cdot N'_n\subset\Img(N'_{n+p}\to N'_n)$ and finally
$\fm\cdot N'_n=\fm^2\cdot N'_n\subset
\Img(\fm\cdot N'_{n+p}\to \fm\cdot N'_n)$ which means that 
$\{\fm\cdot N'_n\}$ is a surjective inverse system, so 
its $\liminv{}^1$ vanishes and the result follows.
\end{proof}

\begin{example} Let $(V,\fm)$ be as in example \ref{ex_rings}.
Then every ideal in $V$ is principal, so in the situation of
the lemma we can write $\eps_j=(x_j)$ for some $x_j\in V$.
Then the hypothesis in (iii) can be stated by saying that
there exists $c\in\N$ such that $x_j\neq 0$ for all $j\geq c$
and the sequence $n\mapsto\sum_{j=c}^n\nu(x_j)$ is Cauchy in 
$\Gamma$.
\end{example}

\begin{definition} Let $M$ be an $A$-module. 

i) We say that $M$ is {\em flat\/} (resp. {\em faithfully
flat\/}) if the functor $N\mapsto M\otimes_A N$, from the 
category of $A$-modules to itself is exact (resp.
exact and faithful). $M$ is {\em almost projective\/} if 
the functor {\em $N\mapsto\Alhom_A(M,N)$} is exact.

ii) We say that $M$ is {\em finitely generated\/} if there 
exists a positive integer $n$ and an epimorphism $A^n\to M$.
We say that $M$ is {\em almost finitely generated\/} if, for 
arbitrary $\eps\in\fm$, there exists a finitely generated
submodule $M_\eps\subset M$ such that $\eps\cdot M\subset M_\eps$.

iii)  We say that $M$ is {\em almost finitely presented\/} if, 
for arbitrary $\eps,\delta\in\fm$ there exist positive integers 
$n=n(\eps)$, $m=m(\eps)$ and a three term complex 
$A^m\stackrel{\psi_\eps}{\to}A^n\stackrel{\phi_\eps}{\to}M$ 
with $\eps\cdot\Coker(\phi_\eps)=0$ and
$\delta\cdot\Ker(\phi_\eps)\subset\Img(\psi_\eps)$.
\end{definition}

\begin{proposition}\label{prop_equiv.cond} 
(i) An $A$-module $M$ is almost finitely generated if 
and only if for every finitely generated ideal $\fm_0\subset\fm$ 
there exists a finitely generated submodule $M_0\subset M$ 
such that $\fm_0\cdot M\subset M_0$.

(ii) An $A$-module is almost finitely presented if and
only if, for every finitely generated ideal $\fm_0\subset\fm$
there is a complex $A^m\stackrel{\psi}{\to}A^n\stackrel{\phi}{\to}M$
with  $\fm_0\cdot\Coker(\phi)=0$ and 
$\fm_0\cdot\Ker(\phi)\subset\Img(\psi)$.
\end{proposition}
\begin{proof} (i) is easy and we leave it to the reader.
To prove (ii), take a finitely generated ideal $\fm_1\subset\fm$
such that $\fm_0\subset\fm\cdot\fm_1$, pick a morphism 
$\phi:A^n\to M$ whose cokernel is annihilated by $\fm_1$,
and apply the following lemma \ref{lem_fin.pres}.
\end{proof}
\begin{lemma}\label{lem_fin.pres}
If $M$ is almost finitely presented and $\phi:F\to M$ is
a morphism with $F\simeq A^n$, then for every finitely 
generated ideal $\fm_1\subset\fm\cdot\Ann_V(\Coker(\phi))$
there is a finitely generated submodule of $\Ker(\phi)$ 
containing $\fm_1\cdot\Ker(\phi)$.
\end{lemma}
\begin{proof}We need the following
\begin{claim}\label{cl_fin.pres}
Let $F_1$ be a finitely generated $A$-module
and suppose that we are given $a,b\in V$ and a (not necessarily 
commutative) diagram
$$\xymatrix{F_1 \ar[r]^p \ar@<.5ex>[d]^\phi & M \\
F_2 \ar@<.5ex>[u]^\psi \ar[ur]_q
}$$
such that $q\circ\phi=a\cdot p$, $p\circ\psi=b\cdot q$. Let
$I\subset V$ be an ideal such that $\Ker(q)$ has a finitely
generated submodule containing $I\cdot\Ker(q)$. Then
$\Ker(p)$ has a finitely generated submodule containing
$a\cdot b\cdot I\cdot\Ker(p)$.
\end{claim}
\begin{pfclaim} Let $R$ be the submodule of $\Ker(q)$
given by the assumption. We have 
$\Img(\psi\circ\phi-a\cdot b\cdot\one_{F_1})\subset\Ker(p)$
and $\psi(R)\subset\Ker(p)$. We take 
$R_1=\Img(\psi\circ\phi-a\cdot b\cdot\one_{F_1})+\psi(R)$.
Clearly $\phi(\Ker(p))\subset\Ker(q)$, so 
$I\cdot\phi(\Ker(p))\subset R$, hence
$I\cdot\psi\circ\phi(\Ker(p))\subset\psi(R)$ and finally
$a\cdot b\cdot I\cdot\Ker(p)\subset R_1$.
\end{pfclaim}

Now, let $\delta\in\text{Ann}_V(\Coker(\phi))$ and 
$\eps_1,\eps_2,\eps_3,\eps_4\in\fm$. By assumption
there is a complex $A^r\stackrel{t}{\to}A^s\stackrel{q}{\to}M$
with $\eps_1\cdot\Coker(q)=0$, 
$\eps_2\cdot\Ker(q)\subset\Img(t)$. Letting $F_1=F$, $F_2=A^s$, 
$a=\eps_1\cdot\eps_3$, $b=\eps_4\cdot\delta$, one checks easily 
that $\psi$ and $\phi$ can be given such that all the assumptions
of the above claim are fulfilled. So, with $I=\eps_2\cdot V$ we get
that 
$\eps_1\cdot\eps_2\cdot\eps_3\cdot\eps_4\cdot\delta\cdot\Ker(\phi)$
lies in a finitely generated submodule of $\Ker(\phi)$.
But $\fm_1$ is contained in an ideal generated by finitely
many such products 
$\eps_1\cdot\eps_2\cdot\eps_3\cdot\eps_4\cdot\delta$.
\end{proof}

The following proposition generalises a well-known characterization
of finitely presented modules over usual rings.
\begin{proposition}\label{prop_al.small} 
Let $M$ be an $A$-module.

i) $M$ is almost finitely generated if and only if, for every 
filtered inductive system $(N_\lambda,\phi_{\lambda\mu})$ (indexed
by a directed set $\Lambda$) the natural morphism
$$\nu:\colim{\Lambda}\Alhom_A(M,N_\lambda)\to
\Alhom_A(M,\colim{\Lambda}N_\lambda)$$
is a monomorphism.

ii) $M$ is almost finitely presented if and only if for every
filtered inductive sytem as above, $\nu$ is an isomorphism.
\end{proposition}
\begin{proof} The ``only if'' part in (i) (resp. (ii)) is first 
checked when $M$ is finitely generated (resp. finitely presented)
and then extended to the general case. We leave the details
to the reader and we proceed to verify the ``if'' part.
For (i), choose a set $I$ and an epimorphism $p:A^{(I)}\to M$.
Let $\Lambda$ be the directed set of finite subsets of $I$, ordered
by inclusion. For $S\in\Lambda$, let $M_S=p(A^S)$. Then 
$\colim{\Lambda}(M/M_S)=0$, so the assumption gives 
$\colim{\Lambda}\Alhom_A(M,M/M_S)=0$, {\em i.e.\/} 
$\colim{\Lambda}\Hom_A(M,M/M_S)=0$ is almost zero, so, for
every $\eps\in\fm$, the image of $\eps\cdot\one_M$ in the 
above colimit is $0$, {\em i.e.\/} there exists $S\in\Lambda$
such that $\eps\cdot M\subset M_S$, which proves the contention. 
For (ii), we present $M$ as a filtered colimit 
$\colim{\Lambda}M_\lambda$, where each $M_\lambda$ is finitely
presented (this can be done {\em e.g.\/} by taking such a
presentation of the $A_*$-module $M_*$ and applying $N\mapsto N^a$).
The assumption of (ii) gives that 
$\colim{\Lambda}\Hom_A(M,M_\lambda)\to\Hom_A(M,M)$ is an almost
isomorphism, hence, for every $\eps\in\fm$ there is 
$\lambda\in\Lambda$ and $\phi_\eps:M\to M_\lambda$ such that 
$p_\lambda\circ\phi_\eps=\eps\cdot\one_M$, where 
$p_\lambda:M_\lambda\to M$ is the natural morphism to the
colimit. If such a $\phi_\eps$ exists for $\lambda$, then 
it exists for every $\mu\geq\lambda$. Hence, if $\fm_0\subset\fm$
is a finitely generated subideal, say $\fm_0=\sum_j^kV\eps_j$,
then there exist $\lambda\in\Lambda$ and $\phi_i:M\to M_\lambda$
such that $p_\lambda\circ\phi_i=\eps_i\cdot\one_M$ for $i=1,...,k$.
Hence $\Img(\phi_i\circ p_\lambda-\eps_i\cdot\one_{M_\lambda})$
is contained in $\Ker(p_\lambda)$ and contains 
$\eps_i\cdot\Ker(p_\lambda)$. Hence $\Ker(p_\lambda)$ has a 
finitely generated submodule $L$ containing 
$\fm_0\cdot\Ker(p_\lambda)$. Choose a presentation 
$A^m\to A^n\stackrel{\pi}{\to} M_\lambda$. Then one can 
lift $\fm_0\cdot L$ to a finitely generated submodule 
$L'$ of $A^n$. Then $\Ker(\pi)+L'$ is a finitely generated 
submodule of $\Ker(p_\lambda\circ\pi)$ containing
$\fm_0^2\cdot\Ker(p_\lambda\circ\pi)$. Since we also have 
$\fm_0\cdot\Coker(p_\lambda\circ\pi)=0$ and $\fm_0$ is
arbitrary, the conclusion follows from proposition 
\ref{prop_equiv.cond}.
\end{proof}

\begin{lemma} Let $0\to M'\to M\to M''\to 0$ be an exact
sequence of $A$-modules. Then:

i) If $M'$, $M''$ are almost finitely generated (resp. presented)
then so is $M$.

ii) If $M$ is almost finitely presented, then $M''$ is almost
finitely presented if and only if $M'$ is almost finitely
generated.
\end{lemma}
\begin{proof} These facts can be deduced from proposition
\ref{prop_al.small} and remark \ref{rem_flproj}(iii), or proved 
directly.
\end{proof}

\begin{lemma} Let $\bP$ be one of the properties : ``flat'',
``almost projective'', ``almost finitely generated'', ``almost 
finitely presented''. If $B$ is a $\bP$ $A$-algebra, and
$M$ is a $\bP$ $B$-module, then $M$ is $\bP$ as an $A$-module. 
\end{lemma}
\begin{proof} Left to the reader.
\end{proof}
Let $R$ be a $V$-algebra and $M$ a flat (resp. faithfully flat)
$R$-module (in the usual sense, see \cite{Mat} p.45). Then $M^a$ 
is a flat (resp. faithfully flat) $R^a$-module. Indeed, the 
functor $M\otimes_R-$ preserves the Serre subcategory of almost 
zero modules, so by general facts it induces an exact functor 
on the localized categories (cp. \cite{Ga} p.369). For the 
faithfullness we have to show that an $R$-module $N$ is almost 
zero whenever $M\otimes_RN$ is almost zero. However, $M\otimes_RN$
is almost zero $\Leftrightarrow$ $M\otimes_R(\fm\otimes_VN)=0$
$\Leftrightarrow$ $\fm\otimes_VN=0$ $\Leftrightarrow$
$N$ is almost zero.
It is clear that $A\Mod$ has enough almost projective (resp. flat) 
objects. Let $R$ be a $V$-algebra. The localisation functor induces a 
functor $G:\sD(R)\to\sD(R^a)$ and, in view of corollary 
\ref{cor_left.exact}, $M\mapsto M_!$ induces a functor 
$F:\sD(R^a)\to\sD(R)$. We have a natural isomorphism 
$G\circ F\simeq\one_{\sD(R^a)}$ and a natural transformation
$F\circ G\to\one_{\sD(R)}$. These satisfy the triangular identities
of \cite{Ma} (p.83) so $F$ is a left adjoint to $G$. If $\Sigma$ denotes
the multiplicative set of morphisms in $\sD(R)$ which induce almost
isomorphisms on the cohomology modules, then the localised category
$\Sigma^{-1}\sD(R)$ exists (see {\em e.g.\/} \cite{We} (Th.10.3.7))
and by the same argument we get an equivalence of categories
$\Sigma^{-1}\sD(R)\simeq\sD(R^a)$.

Given an $A$-module $M$, we can derive the functors  
$M\otimes_A-$ (resp. $\Alhom_A(M,-)$, resp. $\Alhom_A(-,M)$) by taking 
flat (resp. injective, resp. almost projective) resolutions : one remarks 
that bounded above exact complexes of flat (resp. almost projective) 
$A$-modules are acyclic for the functor $M\otimes_A-$ (resp. 
$\Alhom_A(-,M)$) (recall the standard argument: if $F_\bullet$
is a complex of flat $A$-modules, let $\Phi_\bullet$ be a flat 
resolution of $M$; then 
$\text{Tot}(\Phi_\bullet\otimes_AF_\bullet)\to M\otimes_AF_\bullet$
is a quasi-isomorphism since it is so on rows, and 
$\text{Tot}(\Phi_\bullet\otimes_AF_\bullet)$ is acyclic 
since its colums are; similarly, if $P_\bullet$ is a complex of 
almost projective objects, one considers the double complex 
$\Alhom_A(P_\bullet,J^\bullet)$ where $J^\bullet$ is an 
injective resolution of $M$; cp. \cite{We} \S 2.7); then one
uses the construction detailed in \cite{We} (Th.10.5.9). We denote 
by $\Tor_i^A(M,-)$ (resp. $\AlExt_A^i(M,-)$, resp. 
$\AlExt^i_A(-,M)$) the corresponding derived functors. If $A=R^a$
for some $V$-algebra $R$, we obtain easily natural isomorphisms
$\Tor^R_i(M,N)^a\simeq\Tor_i^A(M^a,N^a)$
for all $R$-modules $M,N$. A similar result holds for $\Ext^i_R(M,N)$. 

\begin{remark}\label{rem_flproj}
i) Clearly, an $A$-module $M$ is flat (resp. almost projective) if 
and only if $\Tor^A_i(M,N)=0$ (resp. $\AlExt^i_A(M,N)=0$) for all 
$A$-modules $N$ and all $i>0$.

ii) Let $M,N$ be two flat (resp. almost projective) $A$-modules. Then
$M\otimes_AN$ is a flat (resp. almost projective) $A$-module and for
any $A$-algebra $B$, the $B$-module $B\otimes_AM$ is flat 
(resp. almost projective).

iii) Resume the notation of proposition \ref{prop_al.small}.
If $M$ is almost finitely presented, then one has also that the 
natural morphism 
$\colim{\Lambda}\AlExt_A^1(M,N_\lambda)\to
\AlExt_A^1(M,\colim{\Lambda}N_\lambda)$ is a monomorphism.
This is deduced from proposition \ref{prop_al.small}(ii), 
using the fact that $(N_\lambda)$ can be injected into an
inductive system $(J_\lambda)$ of injective almost modules 
({\em e.g.\/} $J_\lambda=E^{\text{Hom}_A(N_\lambda,E)}$, where 
$E$ is an injective cogenerator for $A\Mod$), and by 
applying $\AlExt$ sequences.
\end{remark}

\begin{lemma}\label{lem_thetwo}
Let $M$ be an almost finitely generated $A$-module. Consider
the following properties:

i) $M$ is almost projective.

ii) For arbitrary $\eps\in\fm$ there exist $n(\eps)\in\N$ and 
$A$-linear morphisms
\set\begin{equation}\label{eq_uv}
\xymatrix{
M \ar[r]^-{u_\eps} & A^{n(\eps)} \ar[r]^-{v_\eps} & M 
}\end{equation}
such that $v_\eps\circ u_\eps=\eps\cdot\one_M$.

iii) $M$ is flat.

Then (i) $\Leftrightarrow$ (ii) $\Rightarrow$ (iii).
\end{lemma}
\begin{proof} (ii)$\Rightarrow$(i): for given $\eps\in\fm$, we consider
any $A$-module $N$ and we apply the functor $\AlExt^i_A(-,N)$ to 
\eqref{eq_uv} :
$$\xymatrix{
\AlExt^i_A(M,N)\ar[r] & \AlExt^i_A(A^{n(\eps)},N)\simeq 0\ar[r] &
\AlExt^i_A(M,N)}$$
which implies $\eps\cdot\AlExt^i(M,N)=0$ for all $i>0$.
Since $\eps$ is arbitrary, (i) follows from remark \ref{rem_flproj}(i).

(i)$\Rightarrow$(ii): by hypothesis, for arbitrary  $\eps\in\fm$
we can find $n=n(\eps)$ and a morphism $\phi_\eps:A^n\to M$ such that
$\eps\cdot\Coker(\phi_\eps)=0$. Let $M_\eps$ be the image of
$\phi_\eps$, so that $\phi_\eps$ factors as 
$A^{n(\eps)}\stackrel{\psi_\eps}{\longrightarrow}M_\eps 
\stackrel{j_\eps}{\longrightarrow}M$.
Also $\eps\cdot\one_M:M\to M$ factors as 
$M\stackrel{\gamma_\eps}{\longrightarrow}M_\eps
\stackrel{j_\eps}{\longrightarrow}M.$
Since by hypothesis $M$ is almost projective, the natural morphism 
induced by $\psi_\eps$
$$\xymatrix{
\Alhom_A(M,A^n) \ar[r]^{\psi_\eps^*} & \Alhom_A(M,M_\eps)
}$$
is an epimorphism. Then for arbitrary $\delta\in\fm$ the morphism 
$\delta\cdot\gamma_\eps$ is in the image of $\psi_\eps^*$, in other words, 
there exists an $A$-linear morphism $u_{\eps\delta}:M\to A^n$ such that 
$\psi_\eps\circ u_{\eps\delta}=\delta\cdot\gamma_\eps$. If now we take
$v_{\eps\delta}=\phi_\eps$, it is clear that 
$v_{\eps\delta}\circ u_{\eps\delta}=\eps\cdot\delta\cdot\one_M$. 
This proves (ii), since the $\eps\in\fm$ satisfying the assertion of 
(ii) form an ideal.

(ii)$\Rightarrow$(iii): for a given $A$-module $N$, apply the 
functor $\Tor^A_i(-,N)$ to the sequence \eqref{eq_uv}. This yields
$\eps\cdot\Tor^A_i(M,N)=0$. Now the claim follows from remark 
\ref{rem_flproj}(i).
\end{proof}
There is a converse to lemma \ref{lem_thetwo} in case $M$ 
is almost finitely presented. Before stating it, we need 
the following lemma.
\begin{lemma}\label{lem_Tor.Ext}
Let $R$ be any ring, $M$ any $R$-module and 
$C=\Coker(\phi:R^n\to R^m)$ any finitely presented (left) $R$-module.
Let $C'=\Coker(\phi^*:R^m\to R^n)$ be the cokernel of the transpose
of the map $\phi$. Then there is a natural isomorphism
$$\Tor^R_1(C',M)\simeq\Hom_R(C,M)/\Img(\Hom_R(C,R)\otimes_RM).$$
\end{lemma}
\begin{proof} We have a spectral sequence :
$$E^2_{ij}=\Tor^R_i(H_j(\Cone(\phi^*)),M)\Rightarrow
H_{i+j}(\Cone(\phi^*)\otimes_RM).$$ 
On the other hand we have also natural isomorphisms 
$$\Cone(\phi^*)\otimes_RM\simeq 
\Hom_R(\Cone(\phi),R)[1]\otimes_RM\simeq 
\Hom_R(\Cone(\phi),M)[1].$$
Hence :
$$\begin{array}{r@{\:\simeq\:}l}
E^2_{10}\simeq E^\infty_{10}\simeq 
H_1(\Cone(\phi^*)\otimes_RM)/E^\infty_{01} &
H^0(\Hom_R(\Cone(\phi),M))/\Img(E^2_{01}) \\
& \Hom_R(C,M)/\Img(\Hom_R(C,R)\otimes_RM)
\end{array}$$
which is the claim.
\end{proof}
\begin{proposition}\label{prop_converse}
(i) Every almost finitely generated almost projective
$A$-module is almost finitely presented.

(ii) Every almost finitely presented flat $A$-module is almost 
projective.

\end{proposition}
\begin{proof} (ii) : let $M$ be such an $A$-module. Let 
$\eps,\delta\in\fm$ and pick a three term complex
$$\xymatrix{
A^m\ar[r]^\psi & A^n\ar[r]^\phi & M
}$$ 
such that
$\eps\cdot\Coker(\phi)=\delta\cdot\Ker(\phi)/\Img(\psi)=0$. 
Set $P=\Coker(\psi_*)$; this is a finitely presented
$A_*$-module and $\phi_*$ factors through a morphism
$\bar\phi_*:P\to M_*$.  Let $\gamma\in\fm$; from lemma 
\ref{lem_Tor.Ext} we see that $\gamma\cdot\bar\phi$ is
the image of some element 
$\sum_{j=1}^n\phi_j\otimes m_j\in\Hom_{A_*}(P,A_*)\otimes_{A_*}M_*$.
If we define $L=A_*^n$ and $v:P\to L$, $w:L\to M_*$ by
$v(x)=(\phi_1(x),...,\phi_n(x))$ and 
$w(y_1,...,y_n)=\sum^n_{j=1}y_j\cdot m_j$, then clearly 
$\gamma\cdot\bar\phi=w\circ v$. Let $K=\Ker(\bar\phi_*)$. 
Then $\delta\cdot K^a=0$ and 
the map $\delta\cdot\one_{P^a}$ factors through a morphism 
$\sigma:(P/K)^a\to P^a$. Similarly the map $\eps\cdot\one_M$
factors through a morphism $\lambda:M\to(P/K)^a$.
Let $\alpha=v^a\circ\sigma\circ\lambda:M\to L^a$
and $\beta=w^a:L^a\to M$. The reader can check that 
$\beta\circ\alpha=\eps\cdot\delta\cdot\gamma\cdot\one_M$.
By lemma \ref{lem_thetwo} the claim follows.

(i) : let $P$ be such an almost finitely generated almost
projective $A$-module. For any finitely generated 
ideal $\fm_0\subset\fm$ pick a morphism $\phi:A^r\to P$
such that $\fm_0\cdot\Coker(\phi)=0$.
If $\eps_1,...,\eps_k$ is a set of generators for $\fm_0$,
a standard argument shows that, for any $i\le k$, 
$\eps_i\cdot\one_P$ lifts to a morphism 
$\psi_i:P\to A^r/\Ker(\phi)$; then, since $P$ is almost 
projective, $\eps_j\psi_i$ lifts to a morphism 
$\psi_{ij}:P\to A^r$. Now claim \ref{cl_fin.pres} applies 
with $F_1=A^r$, $F_2=M=P$, $p=\phi$, $q=\one_P$ and 
$\psi=\psi_{ij}$ and shows that $\Ker(\phi)$ has a finitely 
generated submodule $M_{ij}$ containing 
$\eps_i\cdot\eps_j\cdot\Ker(\phi)$. Then the span of all 
such $M_{ij}$ is a finitely generated submodule 
of $\Ker(\phi)$ containing $\fm_0^2\cdot\Ker(\phi)$. By 
proposition \ref{prop_equiv.cond}(ii), the claim follows.
\end{proof}
\begin{definition}
For an $A$-module $M$, the {\em dual $A$-module\/} of $M$
is the $A$-module {\em $M^*=\Alhom_A(M,A)$}. $M$ is 
{\em reflexive\/} if the natural morphism 
\set\begin{equation}\label{eq_reflex}
M\to (M^*)^*\qquad m\mapsto (f\mapsto f(m))
\end{equation}
is an isomorphism of $A$-modules. 
\end{definition}
\begin{remark}\label{rem_B.struct} Notice that if $B$ is an 
$A$-algebra and $M$ any $B$-module, then by ``restriction of 
scalars'' $M$ is also an $A$-module and the dual $A$-module 
$M^*$ has a natural structure of $B$-module. 
This is defined by the rule $(b\cdot f)(m)=f(b\cdot m)$ 
($b\in\ubar B$, $m\in\ubar M$ and $f\in M^*_*$). 
With respect to this structure \eqref{eq_reflex} becomes a $B$-linear
morphism. Incidentally, notice that the two meanings of ``$M^*_*$''
coincide, {\em i.e.\/} $(M_*)^*\simeq (M^*)_*$.
\end{remark}

\begin{proposition} Let $P$ be an almost projective 
$A$-module and denote by $I_P$ the image of the natural
``evaluation morphism'' $P\otimes_AP^*\to A$. 

i) For every morphism of algebras $A\to B$ we have
$I_{B\otimes_AP}=I_P\cdot B$.

ii) $I_P=I_P^2$.

iii) $P=0$ if and only if $I_P=0$.

iv) $P$ is faithfully flat if and only if $I_P=A$. 
\end{proposition}
\begin{proof} Pick an indexing set $I$ large enough, and an 
epimorphism $\phi:F=A^{(I)}\to P$. For every $i\in I$ we have the
standard morphisms $A\stackrel{e_i}{\to}F\stackrel{\pi_i}{\to}A$
such that $\pi_i\circ e_j=\delta_{ij}\cdot\one_A$ and 
$\sum_{i\in I}e_i\circ\pi_i=\one_F$. For every $x\in\fm$ choose 
$\psi_x\in\Hom_A(P,F)$ such that $\phi\circ\psi_x=x\cdot\one_P$.
It is easy to check that $I_P$ is generated by the almost
elements $\pi_i\circ\psi_x\circ\phi\circ e_j$ 
($i,j\in I$, $x\in\fm$). (i) follows already. For (iii), the 
``only if'' is clear; if $I_P=0$, then $\psi_x\circ\phi=0$ 
for all $x\in\fm$, hence $\psi_x=0$ and therefore 
$x\cdot\one_P=0$. Next, notice that, from (i) and
(iii) we derive $P/(I_P\cdot P)=0$, {\em i.e.\/} $P=I_P\cdot P$, 
so (ii) follows directly from the definition of $I_P$.
Since $P$ is flat, to show (iv) we have only to verify that
the functor $M\mapsto P\otimes_AM$ is faithful. To this purpose, 
it suffices to check that $P\otimes_A(A/J)\neq 0$ for every 
proper ideal $J$ of $A$. This follows easily from (i) and (iii).
\end{proof}
If $E$, $F$ and $N$ are $A$-modules, there is a natural 
morphism :
\set\begin{equation}\label{eq_three.mods}
E\otimes_A\Alhom_A(F,N)\to\Alhom_A(F,E\otimes_AN).
\end{equation}
\begin{lemma}\label{lem_three.mods}
(i) The morphism \eqref{eq_three.mods} is an isomorphism in the 
following cases :

\quad a) when $E$ is flat and $F$ is almost finitely presented;

\quad b) when either $E$ or $F$ is almost finitely generated and
almost projective;

\quad c) when $F$ is almost projective and $E$ is almost finitely
presented;

\quad d) when $E$ is almost projective and $F$ is almost finitely
generated.

(ii) The morphism \eqref{eq_three.mods} is a monomorphism in the 
following cases :

\quad a) when $E$ is flat and $F$ is almost finitely generated;

\quad b) when $E$ is almost projective.

(iii) The morphism \eqref{eq_three.mods} is an epimorphism when
$F$ is almost projective and $E$ is almost finitely generated.
\end{lemma}
\begin{proof} If $F\simeq A^{(I)}$ for some finite set $I$, then 
$\Alhom_A(F,N)\simeq N^{(I)}$ and the claims are obvious. More 
generally, if $F$ is almost projective and almost finitely generated, 
for any $\eps\in\fm$ there exists a finite set $I=I(\eps)$ and morphisms 
\set\begin{equation}\label{eq_CD1}
\xymatrix{
F \ar[r]^{u_\eps} & A^{(I)} \ar[r]^{v_\eps} & F
}\end{equation}
such that $v_\eps\circ u_\eps=\eps\cdot\one_F$. We apply the natural 
transformation 
$$E\otimes_A\Alhom_A(-,N)\to\Alhom_A(-,E\otimes_AN)$$ 
to \eqref{eq_CD1} : an easy diagram chase allows then to conclude 
that the kernel and cokernel of \eqref{eq_three.mods} are killed 
by $\eps$. As $\eps$ is arbitrary, it follows that \eqref{eq_three.mods} 
is an isomorphism in this case. An analogous argument works when 
$E$ is almost finitely generated almost projective, so we get (i.b).
If $F$ is only almost projective, then we still have morphisms
of the type \eqref{eq_CD1}, but now $I(\eps)$ is no longer necessarily
finite. However, the cokernels of the induced morphisms 
$\one_E\otimes u_\eps$ and $\Alhom_A(v_\eps,E\otimes_AN)$ are
still annihilated by $\eps$. Hence, to show (iii) (resp. (i.c)) 
it suffices to consider the case when $F$ is free and $E$ is almost 
finitely generated (resp. presented). By passing to almost elements, 
we can further reduce to the analogous question for usual rings and 
modules, and by the usual juggling we can even replace $E$ by a 
finitely generated (resp. presented) $A_*$-module and $F$ by a 
free $A_*$-module. This case is easily dealt with, and (iii) and 
(i.c) follow. Case (i.d) (resp. (ii.b)) is similar : one considers
almost elements and replaces $E_*$ by a free $A_*$-module (resp.
and $F_*$ by a finitely generated $A_*$-module). In case (ii.a)
(resp. (i.a)), for every finitely generated submodule $\fm_0$ of
$\fm$ we can find, by proposition \ref{prop_equiv.cond}, a finitely 
generated (resp. presented) $A$-module $F_0$ and a morphism 
$F_0\to F$ whose kernel and cokernel are annihilated by $\fm_0$. 
It follows easily that we can replace $F$ by $F_0$ and suppose 
that $F$ is finitely generated (resp. presented). Then the argument 
in \cite{Bou2} (Ch.I \S 2 Prop.10) can be taken over {\em verbatim\/} 
to show (ii.a) (resp. (i.a)).
\end{proof}

\begin{lemma}\label{lem_alhom} i) Let $P$ be an $A$-module
and $B$ an $A$-algebra. If either $P$ or $B$ is almost 
finitely generated almost projective as an $A$-module, 
then the natural morphism 
\set\begin{equation}\label{eq_claim.tens}
B\otimes_A\Alhom_A(P,N)\to\Alhom_B(B\otimes_AP,B\otimes_AN)
\end{equation}
is an isomorphism for all $A$-modules $N$.

ii) Every almost projective almost finitely generated $A$-module
is reflexive.
\end{lemma}
\begin{proof} (i) is an easy consequence of lemma 
\ref{lem_three.mods}(i.b). To prove (ii), we we apply the natural 
transformation \eqref{eq_reflex} to \eqref{eq_CD1} : by diagram 
chase one sees that the kernel and cokernel of the morphism 
$F\to(F^*)^*$ are killed by $\eps$.
\end{proof}

\begin{lemma}\label{lem_limfingen}
Let $\{M_n~;~\phi_n:M_n\to M_{n+1}~|~n\in\N\}$ 
be a direct system of $A$-modules and suppose there exist
sequences $\{\eps_n~|~n\in\N\}$ and $\{\delta_n~|~n\in\N\}$ of 
ideals of $V$ such that

i) $\limdir{n\to\infty}\eps_n=V$ (convergence for the uniform 
structure on ideals of $V$) and $\prod^\infty_{j=0}\delta_j$ 
is a Cauchy product;

ii) for all $n\in\N$ there exist integers $N(n)$ and morphisms 
of $A$-modules $\psi_n:A^{N(n)}\to M_n$ such that 
$\eps_n\cdot\Coker(\psi_n)=0$;

iii) $\delta_n\cdot\Coker(\phi_n)=0$ for all $n\in\N$.

Then $\colim{n\in\N}M_n$ is an almost finitely generated $A$-module.
\end{lemma}
\begin{proof} Let $M=\colim{n\in\N}M_n$. For any $n\in\N$ let
$a_n=\bigcap_{m\ge 0}
(\prod_{j=n}^{n+m}\delta_j)$.
Then $\limdir{n\to\infty}a_n=V$.
For $m>n$ set $\phi_{n,m}=
\phi_m\circ...\circ\phi_{n+1}\circ\phi_n:M_n\to M_{m+1}$ and let 
$\phi_{n,\infty}:M_n\to M$ be the natural morphism. An easy 
induction shows that 
$\prod^m_{j=n}\delta_j\cdot\Coker(\phi_{n,m})=0$
for all $m>n\in\N$. Since 
$\Coker(\phi_{n,\infty})=\colim{m\in\N}\Coker(\phi_{n,m})$ we
obtain $a_n\cdot\Coker(\phi_{n,\infty})=0$ for all $n\in\N$. 
Therefore 
$\eps_n\cdot a_n\cdot\Coker(\phi_{n,\infty}\circ\psi_n)=0$
for all $n\in\N$. Since $\limdir{n\to\infty}\eps_n\cdot a_n=V$, 
the claim follows.
\end{proof}

\begin{lemma}\label{lem_limproj} Let 
$\{M_n~;~\phi_n:M_n\to M_{n+1}~|~n\in\N\}$ be a direct system of 
$A$-modules and suppose there exist sequences $\{\eps_n~|~n\in\N\}$ 
and $\{\delta_n~|~n\in\N\}$ of ideals of $V$ such that

i) $\limdir{n\to\infty}\eps_n=V$ and $\prod^\infty_{j=0}\delta_j$ 
is a Cauchy product;

ii) $\eps_n\cdot\AlExt^i_A(M_n,N)=
\delta_n\cdot\AlExt^i_A(\Coker(\phi_n),N)=0$ for all
$A$-modules $N$, all $i>0$ and all $n\in\N$;

iii) $\delta_n\cdot\Ker(\phi_n)=0$ for all $n\in\N$.

Then $\colim{n\in\N}M_n$ is an almost projective $A$-module.
\end{lemma}
\begin{proof} Let $M=\colim{n\in\N}M_n$. By the above remark 
\ref{rem_flproj}(i)
it suffices to show that $\AlExt^i_A(M,N)$ vanishes for all $i>0$ 
and all $A$-modules $N$. The maps $\phi_n$ define
a map $\phi:\oplus_nM_n\to\oplus_nM_n$ such that we have a short 
exact sequence 
$0\to\oplus_nM_n\stackrel{\one-\phi}{\longrightarrow}\oplus_nM_n
\longrightarrow M\to 0$. Applying the long exact $\AlExt$ sequence
one obtains a short exact sequence (cp. \cite{We} (3.5.10))
$$0\to\liminv{n\in\N}^1\AlExt_A^{i-1}(M_n,N)\to
\AlExt^i_A(M,N)\to\liminv{n\in\N}\AlExt^i_A(M_n,N)\to 0.$$
Then lemma \ref{lem_invlim}(ii) implies that $\AlExt^i_A(M,N)\simeq 0$ 
for all $i>1$ and moreover  $\AlExt^1_A(M,N)$ is isomorphic to 
$\liminv{n\in\N}^1\Alhom_A(M_n,N)$. Let 
$$\phi^*_n:\Alhom_A(M_{n+1},N)\to\Alhom_A(M_n,N)\qquad 
f\mapsto f\circ\phi_n$$
be the transpose of $\phi_n$ and write $\phi_n$ as a composition 
$M_n\stackrel{p_n}{\longrightarrow}\Img(\phi_n)
\stackrel{q_n}{\hookrightarrow}M_{n+1}$,
so that $\phi^*_n=q_n^*\circ p_n^*$, the composition of the 
respective transposed morphims. We have monomorphisms
$$\begin{array}{rl}
&\Coker(p^*_n)\hookrightarrow\Alhom_A(\Ker(\phi_n),N) \\
&\Coker(q^*_n)\hookrightarrow\AlExt^1_A(\Coker(\phi_n),N)
\end{array}$$
for all $n\in\N$. Hence $\delta^2_n\cdot\Coker(\phi^*_n)=0$
for all $n\in\N$. Since $\prod^\infty_{n=0}\delta^2_n$ is a 
Cauchy product, lemma \ref{lem_invlim}(iii) shows that 
$\liminv{n\in\N}^1\Alhom_A(M_n,N)\simeq 0$ and the assertion 
follows.
\end{proof}

\begin{proposition}\label{prop_comp.supp}
Suppose that $\tilde\fm$ is a flat $V$-module. Then for 
any $V$-algebra $R$ the functor $M\mapsto M_!$ commutes with 
tensor products and takes flat $R^a$-modules to 
flat $R$-modules.
\end{proposition}
\begin{proof} Let $M$ be a flat $R^a$-module and  
$N\hookrightarrow N'$ an injective map of $R$-modules.
Denote by $K$ the kernel of the induced map 
$M_!\otimes_RN\to M_!\otimes_RN'$; we have $K^a\simeq 0$.
We obtain an exact sequence
$0\to\tilde\fm\otimes_VK\to
\tilde\fm\otimes_VM_!\otimes_RN\to
\tilde\fm\otimes_VM_!\otimes_RN'$.
But one sees easily that $\tilde\fm\otimes_VK=0$ and 
$\tilde\fm\otimes_VM_!\simeq M_!$, which shows that $M_!$ is a flat 
$R$-module. Similarly, let $M,N$ be two $R^a$-modules. 
Then the natural map $M_*\otimes_RN_*\to(M\otimes_{R^a}N)_*$
is an almost isomorphism and the assertion follows from
remark \ref{rem_almost.zero}(i).
\end{proof}

\subsection{Almost homotopical algebra}\label{sec_hot}

The formalism of abelian tensor categories provides a minimal 
framework wherein the rudiments of deformation theory can be
developed. 

Let $(\cC,\otimes,U)$ be an abelian tensor category; we assume
henceforth that $\otimes$ is a right exact functor. Let $A$ be a 
given $\cC$-monoid. A {\em two-sided ideal\/} of $A$ is an
$A$-sub-bimodule $I\to A$. The quotient $A/I$
in the underlying abelian category $\cC$ has a unique 
$\cC$-monoid structure such that $A\to A/I$ is a morphism
of monoids. $A/I$ is unitary if $A$ is. For $I,J$ subobjects
of $A$ one denotes 
$IJ=\Img(I\otimes J\to A\otimes A\stackrel{\mu_A}{\to}A)$.
If $I$ is a two-sided ideal of $A$ such that $I^2=0$, then,
using the right exactness of $\otimes$ one checks that $I$ 
has a natural structure of an $A/I$-bimodule, unitary when $A$ is.
\begin{definition} A {\em $\cC$-extension\/} of a $\cC$-monoid 
$B$ by a $B$-bimodule $I$ is a short exact sequence of objects 
of\/ $\cC$
\set\begin{equation}\label{eq_def.ext}
\xymatrix{
X:\qquad 0 \ar[r] & I \ar[r] & C \ar[r]^p & B \ar[r] & 0 
}\end{equation}
such that $C$ is a $\cC$-monoid, $p$ is a morphism of\/
$\cC$-monoids, $I$ is a square zero two-sided ideal in $C$
and the $E/I$-bimodule structure on $I$ coincides with
the given bimodule structure on $I$. The $\cC$-extensions 
form a category $\bExmon_\cC$. The morphisms are commutative 
diagrams with exact rows
$$\diagram
~X:\dto & \qquad
0 \rto & I \rto \dto^f & E \rto^p \dto^g & B \rto \dto^h & 0 \\
~X': & \qquad 0 \rto & I' \rto & E' \rto^{p'} & B' \rto & 0
\enddiagram$$
such that $g$ and $h$ are morphisms of\/ $\cC$-monoids. We 
let $\bExmon_\cC(B,I)$ be the subcategory of\/ $\bExmon_\cC$ 
consisting of all $\cC$-extensions of $B$ by $I$, where the 
morphisms are all short exact sequences as above such that 
$f=\one_I$ and $h=\one_B$.
\end{definition}

We have also the variant in which all the $\cC$-monoids 
in \eqref{eq_def.ext} are required to be unitary (resp. 
to be algebras) and $I$ is a unitary $B$-bimodule (resp. 
whose left and right $B$-module actions coincide, {\em i.e.\/} 
are switched by composition with the ``commutativity 
constraints'' $\eta_{B|I}$ and $\eta_{I|B}$, see 
\ref{sec_alm.cat}); we will call $\bExun_\cC$ (resp. 
$\bExal_\cC$) the corresponding category. 
For a morphism $\phi:C\to B$ of $\cC$-monoids, and a $\cC$-extension 
$X$ in $\bExmon_\cC(B,I)$, we can pullback $X$ via $\phi$ to obtain
an exact sequence $X*\phi$ with a morphism $\phi^*:X*\phi\to X$;
one checks easily that there exists a unique structure of
$\cC$-extension on $X*\phi$ such that $\phi^*$ is a morphism
of $\cC$-extension; then $X*\phi$ is an object in $\bExmon_\cC(C,I)$. 
Similarly, given a $B$-linear morphism $\psi:I\to J$, we 
can push out $X$ and obtain a well defined object $\psi*X$ 
in $\bExmon_\cC(B,J)$ with a morphism $X\to\psi*X$ of $\bExmon_\cC$.
In particular, if $I_1$ and $I_2$ are two $B$-bimodules,
the functors $p_i*$ ($i=1,2$) associated to the natural 
projections $p_i:I_1\oplus I_2\to I_i$ establish an equivalence
of categories
\set\begin{equation}\label{eq_exal.p_i}
\bExmon_\cC(B,I_1\oplus I_2)\stackrel{\sim}{\to}
\bExmon_\cC(B,I_1)\times\bExmon_\cC(B,I_2)
\end{equation}
whose essential inverse is given by 
$(E_1,E_2)\mapsto (E_1\oplus E_2)*\delta$, where 
$\delta:B\to B\oplus B$ is the diagonal morphism.
A similar statement holds for $\bExal$ and $\bExun$.
These operations can be used to induce an abelian group 
structure on the set $\Exmon_\cC(B,I)$ of isomorphism 
classes of objects of $\bExmon_\cC(B,I)$ as follows. For any 
two objects $X,Y$ of $\bExmon_\cC(B,I)$ we can form $X\oplus Y$ 
which is an object of $\bExmon_\cC(B\oplus B,I\oplus I)$. Let  
$\alpha:I\oplus I\to I$ be the addition morphism of $I$. Then we 
set $X+Y=\alpha*(X\oplus Y)*\delta$. One can check that $X+Y\simeq Y+X$ 
for any $X,Y$ and that the trivial split $\cC$-extension $B\oplus I$ 
is a neutral element for $+$. Moreover every isomorphism class 
has an inverse $-X$. The functors $X\mapsto X*\phi$ and 
$X\mapsto\psi*X$ commute with the operation thus defined, and 
induce group homomorphisms
$$\begin{array}{l}
*\phi:\Exmon_\cC(B,I)\to\Exmon_\cC(C,I) \\
\psi*:\Exmon_\cC(B,I)\to\Exmon_\cC(B,J).
\end{array}$$
We will need the variant $\Exal_\cC(B,I)$ defined in the same 
way, starting from $\bExal_\cC(B,I)$. For instance, if $A$ is 
an almost algebra (resp. a commutative ring), we can consider 
the abelian tensor category $\cC=A\Mod$. 
In this case the $\cC$-extensions will be called 
simply $A$-extensions, and we will write $\bExal_A$ rather 
than $\bExal_\cC$. In fact the commutative unitary case will 
soon become prominent in our work, and the more general setup 
is only required for technical reasons, in the proof of 
proposition \ref{prop_lift.idemp} below, which is the abstract 
version of a well-known result on the lifting of idempotents
over nilpotent ring extensions. 

Let $A$ be a $\cC$-monoid. We form the biproduct 
$A^\dagger=U\oplus A$ in $\cC$. We denote by $p_1$, $p_2$ 
the associated projections from $A^\dagger$ onto $U$ and
respectively $A$. Also, let $i_1$, $i_2$ be the natural
monomorphisms from $U$, resp. $A$ to $A^\dagger$. $A^\dagger$
is equipped with a unitary monoid structure 
$$\mu^\dagger=i_2\circ\mu\circ(p_2\otimes p_2)+
i_2\circ\ell^{-1}_A\circ(p_1\otimes p_2)+
i_2\circ r^{-1}_A\circ(p_2\otimes p_1)+
i_1\circ u^{-1}\circ(p_1\otimes p_1)$$
where $\ell_A$, $r_A$ are the natural isomorphisms 
provided by \cite{DeM} (Prop. 1.3) and $u:U\to U\otimes U$ 
is as in {\em loc. cit.\/} \S 1. In terms of the ring
$A^\dagger_*\simeq U_*\oplus A_*$ this is the multiplication
$(u_1,b_1)\cdot(u_2,b_2)=
(u_1\cdot u_2,b_1\cdot b_2+b_1\cdot u_2+u_1\cdot b_2)$.
Then $i_2$ is a morphism of monoids and one verifies that
the ``restriction of scalars'' functor $i_2^*$ defines an 
equivalence from the category $A^\dagger\UniMod$ of unitary 
$A^\dagger$-modules to the category $A\Mod$ of all 
$A$-modules; let $j$ denote the inverse functor. A similar 
discussion applies to bimodules. Similarly, we derive 
equivalences of categories
$$\xymatrix{
\bExun_\cC(A^\dagger,j(M))\ar@<.5ex>[r]^-{*i_2} &
\bExmon_\cC(A,M)\ar@<.5ex>[l]^-{(-)^\dagger}
}$$
for all $A$-bimodules $M$.

Next we specialise to $A=U$ : for a given $U$-module $M$ 
let $e_M=\sigma_M\circ\ell_M:M\to M$; working out the 
definitions one finds that the condition that this is a 
module structure is equivalent to $e_M^2=e_M$. Let 
$U\times U$ be the product of $U$ by itself in the category 
of $\cC$-monoids. There is an isomorphism of unitary 
$\cC$-monoids $\zeta:U^\dagger\to U\times U$
given by $\zeta=i_1\circ p_1+i_2\circ p_1+i_2\circ p_2$.
Another isomorphism is $\tau\circ\zeta$, where $\tau$
is the flip $i_1\circ p_2+i_2\circ p_1$. Hence we get
equivalences of categories
$$\xymatrix{
U\Mod \ar@<.5ex>[r]^-{j} & U^\dagger\UniMod 
\ar@<.5ex>[rr]^-{(\zeta^{-1})^*} 
\ar@<.5ex>[l]^-{i_2^*}&&
(U\times U)\UniMod
\ar@<.5ex>[ll]^-{(\tau\circ\zeta)^*} 
}$$
The composition 
$i_2^*\circ(\zeta^{-1}\circ\tau\circ\zeta)^*\circ j$
defines a self-equivalence of $U\Mod$ which associates
to a given $U$-module $M$ the new $U$-module 
$M^\flip$ whose underlying object in $\cC$ is $M$ 
and such that $e_{M^\flip}=\one_M-e_M$.
The same construction applies to $U$-bimodules and finally
we get equivalences 
\set\begin{equation}\label{eq_flip}
\xymatrix{
\bExmon_\cC(U,M)\ar[r]^-\sim &
\bExmon_\cC(U,M^\flip)
}\qquad X\mapsto X^\flip
\end{equation}
for all $U$-bimodules $M$. If 
$X=(0\to M\to E\stackrel{\pi}{\to}U\to 0)$
is an extension and $X^\flip=(0\to M^\flip\to E^\flip\to U\to 0)$,
then one verifies that there is a natural isomorphism $X^\flip\to X$ 
of complexes in $\cC$ inducing $-\one_M$ on $M$, the identity on 
$U$ and carrying the multiplication morphism on $E^\flip$ to
$$-\mu_E+\ell_E^{-1}\circ(\pi\otimes\one_E)+
r_E^{-1}\circ(\one_E\otimes\pi):E\otimes E\to E.$$
In terms of the associated rings, this corresponds to 
replacing the given multiplication $(x,y)\mapsto x\cdot y$ 
of $E_*$ by the new operation
$(x,y)\mapsto\pi_*(x)\cdot y+\pi_*(y)\cdot x-x\cdot y$.
\begin{lemma}\label{lem_split.ext}
If $M$ is a $U$-bimodule whose left and right 
actions coincide, then every extension of $U$ by $M$ splits
uniquely.
\end{lemma}
\begin{proof} Using the idempotent $e_M$ we get a $U$-linear
decomposition $M\simeq M_1\oplus M_2$ where the bimodule 
structure on $M_1$ is given by the zero morphisms and the
bimodule structure on $M_2$ is given by $\ell_M^{-1}$
and $r_M^{-1}$. We have to prove that $\bExmon_\cC(U,M)$
is equivalent to a one-point category. By \eqref{eq_exal.p_i} 
we can assume that $M=M_1$ or $M=M_2$. By \eqref{eq_flip} we
have $\bExmon_\cC(U,M_2)\simeq\bExmon_\cC(U,M_2^\flip)$
and on $M_2^\flip$ the bimodule actions are the zero morphisms.
So it is enough to consider $M=M_1$. In this case,
if $X=(0\to M\to E\to U\to 0)$ is any extension, 
$\mu_E:E\otimes E\to E$ factors through a morphism
$U\otimes U\to E$ and composing with $u:U\to U\otimes U$
we get a right inverse of $E\to U$, which shows that
$X$ is the split extension. Then it is easy to see that
$X$ does not have any non-trivial automorphisms, which
proves the assertion.
\end{proof}
\begin{proposition}\label{prop_lift.idemp}
i) Let $X=(0\to I\to A\stackrel{p}{\to} A'\to 0)$ be a 
$\cC$-extension and suppose that $e'\in A'_*$ is an 
idempotent element whose left action on the 
$A'$-bimodule $I$ coincides with its right action. 
Then there exists a unique idempotent $e\in A_*$ such 
that $p_*(e)=e'$.

ii) Especially, if $A'$ is unitary and $I$ is a unitary
$A'$-bimodule, then every extension of $A'$ by $I$ is unitary.
\end{proposition}
\begin{proof} (i) : the hypothesis ${e'}^2=e'$ implies that 
$e':U\to A'$ is a morphism of (non-unitary) $\cC$-monoids. 
We can then replace $X$ by $X*e'$ and thereby assume that
$A'=U$, $p:A\to U$ and $I$ is a (non-unitary) $U$-bimodule 
and the right and left actions on $I$ coincide. 
The assertion to prove is that $\ubarold 1_U$ lifts to
a unique idempotent $e\in A_*$. However, this follows easily 
from lemma \ref{lem_split.ext}. To show (ii), we observe that,
by (i), the unit $\ubarold 1_{A'}$ of $A'_*$ lifts uniquely
to an idempotent $e\in A_*$. We have to show that $e$ is a unit
for $A_*$. Let us show the left unit property. Via $e:U\to A$
we can view the extension $X$ as an exact sequence of left 
$U$-modules. We can then split $X$ as the direct sum 
$X_1\oplus X_2$ where $X_1$ is a sequence of unitary $U$-modules
and $X_2$ is a sequence of $U$-modules with trivial actions.
But by hypothesis, on $I$ and on $A$ the $U$-module structure
is unitary, so $X=X_1$ and this is the left unit property.
\end{proof}

So much for the general nonsense; we now return to almost algebras.
As already announced, {\em from here on, we assume throughout that
$\tilde\fm$ is a flat $V$-module}.
As an immediate consequence of proposition \ref{prop_lift.idemp} 
we get natural equivalences of categories 
\set\begin{equation}\label{eq_equi.prod}
\xymatrix{
\bExal_{A_1}(B_1,M_1)\times\bExal_{A_2}(B_2,M_2) \ar[r]^-\sim &
\bExal_{A_1\times A_2}(B_1\times B_2,M_1\oplus M_2)
}\end{equation}
whenever $A_1$, $A_2$ are $V^a$-algebras, $B_i$ is 
a $A_i$-algebra and $M_i$ is a (unitary) $B_i$-module, $i=1,2$.

Notice that, if $A=R^a$ for some $V$-algebra $R$, $S$ (resp. 
$J$) is a $R$-algebra (resp. an $S$-module) and $X$ is any object 
of $\bExal_R(S,J)$, then by applying termwise the localisation 
functor we get an object $X^a$ of $\bExal_A(S^a,J^a)$.
With this notation we have the following lemma.
\begin{lemma}\label{lem_eq.Ext} 
i) Let $B$ be any $A$-algebra and $I$ a $B$-module. The natural 
functor
\set\begin{equation}\label{eq_funct.ext}
\bExal_{A_{!!}}(B_{!!},I_*)\to\bExal_A(B,I)\qquad X\mapsto X^a
\end{equation}
is an equivalence of categories.

ii) The equivalence \eqref{eq_funct.ext} induces a group isomorphism
$\Exal_{A_{!!}}(B_{!!},I_*)\stackrel{\sim}{\to}\Exal_A(B,I)$
functorial in all arguments.
\end{lemma}
\begin{proof} Of course (ii) is an immediate consequence of (i).
To show (i), let $X=(0\to I\to E\to B\to 0)$ be any object of
$\bExal_A(B,I)$. Using corollary \ref{cor_left.exact} one sees 
easily that the sequence $X_!=(0\to I_!\to E_{!!}\to B_{!!}\to 0)$ 
is right exact; $X_!$ won't be exact in general, unless $B$
(and therefore $E$) is an exact algebra. In any case, the 
kernel of $I_!\to E_{!!}$ is almost zero, so we get an extension
of $B_{!!}$ by a quotient of $I_!$ which maps to $I_*$. In
particular we get by pushout an extension $X_{!*}$ by $I_*$,
{\em i.e.\/} an object of $\bExal_{A_{!!}}(B_{!!},I_*)$ and 
in fact the assignment $X\mapsto X_{!*}$ is an essential 
inverse for the functor \eqref{eq_funct.ext}.
\end{proof}
\begin{remark}\label{rem_subseq} By inspecting the proof, 
we see that one can replace $I_*$ by $I_{!*}=\Img(I_!\to I_*)$ 
in (i) and (ii) above. When $B$ is exact, also $I_!$ will do.
\end{remark}

In \cite{Il} (II.1.2) it is shown how to associate to any ring
homomorphism $R\to S$ a natural simplicial complex of $S$-modules
denoted $\L_{S/R}$ and called the cotangent complex of $S$ over
$R$.

\begin{definition} Let $A\to B$ be a morphism of almost 
$V$-algebras. The {\em almost cotangent complex\/} of $B$ over
$A$ is the simplicial $B_{!!}$-module 
$$\L_{B/A}=B_{!!}\otimes_{(V^a\times B)_{!!}}
\L_{(V^a\times B)_{!!}/(V^a\times A)_{!!}}.$$
\end{definition}

Usually we will want to view $\L_{B/A}$ as an object of the derived 
category $\sD_\bullet(s.B_{!!})$ of simplicial $B_{!!}$-modules. 
Indeed, the hyperext functors computed in this category relate 
the cotangent complex to a number of important invariants. Recall
that, for any simplicial ring $R$ and any two $R$-modules $E,F$ 
the hyperext of $E$ and $F$ is the abelian group defined as
$$\hExt^p_R(E,F)=
\colim{n\ge -p}\Hom_{\sD_\bullet(R)}(\sigma^nE,\sigma^{n+p}F)$$
(where $\sigma$ is the suspension functor of \cite{Il} (I.3.2.1.4)).

Let us fix an almost algebra $A$. First we want to establish 
the relationship with differentials.

\begin{definition}\label{def_derivation}
Let $B$ be any $A$-algebra, $M$ any $B$-module. 

i) An {\em $A$-derivation\/} of $B$ with values in $M$ 
is an $A$-linear morphism $\partial:B\to M$ such that 
$\partial(b_1\cdot b_2)=b_1\cdot\partial(b_2)+
b_2\cdot\partial(b_1)$
for $b_1,b_2\in B_*$. The set of all $M$-valued
$A$-derivations of $B$ forms a $V$-module $\Der_A(B,M)$ and 
the almost $V$-module $\Der_A(B,M)^a$ has a natural structure 
of $B$-module.

ii) We reserve the notation $I_{B/A}$ for the ideal 
$\Ker(\mu_{B/A}:B\otimes_AB\to B)$. The {\em module 
of relative differentials\/} of $\phi$ is defined as the 
(left) $B$-module $\Omega_{B/A}=I_{B/A}/I^2_{B/A}$. 
It is endowed with a natural $A$-derivation 
$\delta:B\to\Omega_{B/A}$ defined by 
$b\mapsto\ubarold 1\otimes b-b\otimes\ubarold 1$ for all 
$b\in\ubar B$. The assignment $(A\to B)\mapsto\Omega_{B/A}$ 
defines a functor
$$\Omega:V^a\AlgMorph\to V^a\AlgMod$$
from the category of morphisms $A\to B$ of almost $V$-algebras
to the category $V^a\AlgMod$ consisting of all pairs $(B,M)$ 
where $B$ is an almost $V$-algebra and $M$ is a $B$-module.
The morphisms in $V^a\AlgMorph$ are the commutative squares; the 
morphisms $(B,M)\to(B',M')$ in $V^a\AlgMod$ are all pairs $(\phi,f)$
where $\phi:B\to B'$ is a morphism of almost $V$-algebras and
$f:B'\otimes_BM\to M'$ is a morphism of $B'$-modules.
\end{definition}

The module of relative differentials enjoys the familiar 
universal properties that one expects. In particular
$\Omega_{B/A}$ represents the functor $\Der_A(B,-)$, {\em i.e.}
for any (left) $B$-module $M$ the morphism
\set\begin{equation}\label{eq_repres}
\Hom_B(\Omega_{B/A},M)\to\Der_A(B,M)\qquad f\mapsto f\circ\delta
\end{equation}
is an isomorphism. As an exercise, the reader can supply the proof
for this claim and for the following standard proposition.
\begin{proposition}\label{prop_standard}
i) Let $B$ and $C$ be two $A$-algebras. Then there is a natural 
isomorphism: $$\Omega_{C\otimes_AB/C}\simeq C\otimes_A\Omega_{B/A}.$$

ii) Let $B$ be an $A$-algebra, $C$ a $B$-algebra.
There is a natural exact sequence of $C$-modules:
$$C\otimes_B\Omega_{B/A}\to\Omega_{C/A}\to\Omega_{C/B}\to 0.$$

iii) Let $I$ be an ideal of the $A$-algebra $B$ and let $C=B/I$
be the quotient $A$-algebra. Then there is a natural exact sequence:
$I/I^2\to C\otimes_B\Omega_{B/A}\to\Omega_{C/A}\to 0$.

iv) The functor $\Omega:V^a\AlgMorph\to V^a\AlgMod$ commutes with all
colimits.
\qed\end{proposition}

\begin{lemma}\label{lem_differ} For any $A$-algebra $B$ 
there is a natural isomorphism of $B_{!!}$-modules
$$(\Omega_{B/A})_!\simeq\Omega_{B_{!!}/A_{!!}}.$$
\end{lemma}
\begin{proof} Using the adjunction \eqref{eq_repres}
we are reduced to showing that the natural map
$$\phi_M:\Der_{A_{!!}}(B_{!!},M)\to\Der_A(B,M^a)$$ 
is a bijection for all $B_{!!}$-modules $M$. Given 
$\partial:B\to M^a$ we construct $\partial_!:B_!\to M^a_!\to M$. 
We extend $\partial_!$ to $V\oplus B_!$ by setting it equal 
to zero on $V$. Then it is easy to check that the 
resulting map descends to $B_{!!}$, hence giving an  
$A$-derivation $B_{!!}\to M$.
This procedure yields a right inverse $\psi_M$ to 
$\phi_M$. To show that $\phi_M$ is injective, suppose
that $\partial:B_{!!}\to M$ is an almost zero $A$-derivation.
Composing with the natural $A$-linear map $B_!\to B_{!!}$
we obtain an almost zero map $\partial':B_!\to M$. But 
$\fm\cdot B_!=B_!$, hence $\partial'=0$. This implies
that in fact $\partial=0$, and the assertion follows.
\end{proof}

\begin{proposition}\label{prop_H_0}
Let $M$ be a $B$-module. There exists a natural 
isomorphism of $B_{!!}$-modules
$$\hExt^0_{B_{!!}}(\L_{B/A},M_!)\simeq\Der_A(B,M).$$
\end{proposition}
\begin{proof} To ease notation, set $\tilde A=V^a\times A$
and $\tilde B=V^a\times B$. We have natural isomorphisms :
$$\begin{array}{l@{\:\simeq\:}ll}
\hExt^0_{B_{!!}}(\L_{B/A},M_!) &
\hExt^0_{\tilde B_{!!}}(\L_{\tilde B_{!!}/\tilde A_{!!}},M_!)
& \quad\text{by \cite{Il} (I.3.3.4.4)} \\
& \Der_{\tilde A_{!!}}(\tilde B_{!!},M_!)
& \quad\text{by \cite{Il} (II.1.2.4.2)} \\
& \Der_{\tilde A}(\tilde B,M) 
& \quad\text{by lemma \ref{lem_differ}.} 
\end{array}$$
But it is easy to see that the natural map
$\Der_A(B,M)\to\Der_{\tilde A}(\tilde B,M)$ is an 
isomorphism.
\end{proof}
\begin{theorem}\label{th_main}
There is a natural isomorphism
$$\Exal_A(B,M)\to\hExt^1_{B_{!!}}(\L_{B/A},M_!).$$
\end{theorem}
\begin{proof} With the notation of the proof of proposition
\ref{prop_H_0} we have natural isomorphisms 
$$\begin{array}{l@{\:\simeq\:}ll}
\hExt^1_{B_{!!}}(\L_{B/A},M_!) 
& \hExt^1_{\tilde B_{!!}}(\L_{\tilde B_{!!}/\tilde A_{!!}},M_!) 
& \quad\text{by \cite{Il} (I.3.3.4.4)} \\
& \Exal_{\tilde A_{!!}}(\tilde B_{!!},M_!) 
& \quad\text{by \cite{Il} (III.1.2.3)} \\
& \Exal_{\tilde A}(\tilde B,M) 
\end{array}$$
where the last isomorphism follows directly from lemma 
\ref{lem_eq.Ext}(ii) and the subsequent remark \ref{rem_subseq}.
Finally, \eqref{eq_equi.prod} shows that 
$\Exal_{\tilde A}(\tilde B,M)\simeq\Exal_A(B,M)$, as required.
\end{proof}
Moreover we have the following transitivity theorem as 
in \cite{Il} (II.2.1.2).
\begin{theorem}\label{th_transit}
Let $A\to B\to C$ be a sequence of morphisms of almost $V$-algebras. 
There exists a natural distinguished triangle of\/ 
$\sD_\bullet(s.C_{!!})$ 
$$\diagram
C_{!!}\otimes_{B_{!!}}\L_{B/A} \rto^-u & \L_{C/A} \rto^v & \L_{C/B}
\rto & C_{!!}\otimes_{B_{!!}}\sigma\L_{B/A}
\enddiagram$$
where the morphisms $u$ and $v$ are obtained by functoriality of\/ $\L$.
\end{theorem}
\begin{proof} It  follows directly from {\em loc. cit\/}. 
\end{proof}

\begin{proposition} Let $(A_\lambda\to B_\lambda)_{\lambda\in I}$ be
a system of almost $V$-algebra morphisms indexed by a small filtered
category $I$. Then there is a natural isomorphism in 
$\sD_\bullet(s.\colim{\lambda\in I}B_{\lambda!!})$
$$\colim{\lambda\in I}\L_{B_\lambda/A_\lambda}\simeq
\L_{\colim{\lambda\in I}B_\lambda/\colim{\lambda\in I}A_\lambda}.$$
\end{proposition}
\begin{proof} Remark \ref{rem_adjoints} gives an 
isomorphism : $\colim{\lambda\in I}A_{\lambda!!}
\stackrel{\sim}{\to}(\colim{\lambda\in I}A_\lambda)_{!!}$ 
(and likewise for $\colim{\lambda\in I}B_\lambda$). Then the 
claim follows from \cite{Il} (II.1.2.3.4).
\end{proof}

Next we want to prove the almost version of the flat base 
change theorem \cite{Il} (II.2.2.1). To this purpose we need 
some preparation.

\begin{proposition}\label{prop_Torindep} Let $B$ and $C$ be two 
$A$-algebras and set $T_i=\Tor_i^{A_{!!}}(B_{!!},C_{!!})$. 
If $A$, $B$, $C$ and $B\otimes_AC$ are all exact, 
then for every $i>0$ the natural morphism 
$\tilde\fm\otimes_VT_i\to T_i$ is an isomorphism.
\end{proposition}
\begin{proof} For any almost $V$-algebra $D$ we let $k_D$ denote
the complex of $D_{!!}$-modules 
$[\tilde\fm\otimes_VD_{!!}\to D_{!!}]$ placed in degrees $-1,0$;
we have a distiguished triangle
$$\cT(D)~:~\xymatrix{
\tilde\fm\otimes_VD_{!!} \ar[r] & D_{!!} \ar[r] 
& k_D \ar[r] & \tilde\fm\otimes_VD_{!!}[1].
}$$
By the assumption,  the natural map $k_A\to k_B$ is a 
quasi-isomorphism and $\tilde\fm\otimes_VB_{!!}\simeq B_!$.
On the other hand, for all $i\in\N$ we have 
$$\Tor ^{A_{!!}}_i(k_B,C_{!!})\simeq
\Tor ^{A_{!!}}_i(k_A,C_{!!})\simeq H^{-i}(k_A\otimes_{A_{!!}}C_{!!})=
H^{-i}(k_C).$$
In particular $\Tor ^{A_{!!}}_i(k_B,C_{!!})=0$ for all $i>1$.
As $\tilde\fm$ is flat over $V$, we have
$\tilde\fm\otimes_VT_i\simeq
\Tor^{A_{!!}}_i(\tilde\fm\otimes_VB_{!!},C_{!!})$.
Then by the long exact Tor sequence associated to 
$\cT(B)\derotimes_{A_{!!}}C_{!!}$ we get the assertion 
for all $i>1$. Next we consider the natural 
map of distinguished triangles 
$\cT(A)\derotimes_{A_{!!}}A_{!!}\to\cT(B)\derotimes_{A_{!!}}C_{!!}$;
writing down the associated morphism of long exact Tor
sequences, we obtain a diagram with exact rows :
$$\xymatrix{
0 \ar[r]  & 
\Tor ^{A_{!!}}_1(k_A,A_{!!}) \ar[r]^-\partial \ar[d] &
(\tilde\fm\otimes_VA_{!!})\otimes_{A_{!!}}A_{!!} \ar[r]^-i \ar[d] &
A_{!!}\otimes_{A_{!!}}A_{!!} \ar[d] \\ 
& \Tor ^{A_{!!}}_1(k_B,C_{!!}) \ar[r]^-{\partial'} &
(\tilde\fm\otimes_VB_{!!})\otimes_{A_{!!}}C_{!!} \ar[r]^-{i'} &
B_{!!}\otimes_{A_{!!}}C_{!!} 
}$$
By the above, the leftmost vertical map is an isomorphism;
moreover, the assumption gives 
$\Ker(i)\simeq\Ker(\tilde\fm\to V)\simeq\Ker(i')$. Then,
since $\partial$ is injective, also $\partial'$ must be
injective, which implies our assertion for the remaining
case $i=1$.
\end{proof}

\begin{corollary}\label{cor_Toindep} Keep the notation of proposition 
\ref{prop_Torindep} and suppose that $\Tor_i^A(B,C)\simeq 0$
 for some $i>0$. Then the corresponding $T_i$ vanishes.
\qed\end{corollary}

\begin{theorem}\label{th_flat.base.ch} Let $B$, $A'$ be 
two $A$-algebras. Suppose that the natural morphism 
$B\derotimes_AA'\to B'=B\otimes_AA'$ is an isomorphism in 
$\sD_\bullet(s.A)$. Then the natural morphisms
$$\begin{array}{l}
B'_{!!}\otimes_{B_{!!}}\L_{B/A}\longrightarrow\L_{B'/A'} \\
(B'_{!!}\otimes_{B_{!!}}\L_{B/A})\oplus
(B'_{!!}\otimes_{A'_{!!}}\L_{A'/A})\longrightarrow\L_{B'/A}
\end{array}$$
are quasi-isomorphisms.
\end{theorem}
\begin{proof} Let us remark that the functor 
$D\mapsto V^a\times D$ : $A\Alg\to(V^a\times A)\Alg$ 
commutes with tensor products; hence the same holds for the 
functor $D\mapsto(V^a\times D)_{!!}$ (see remark 
\ref{rem_adjoints}). Then, in view of corollary 
\ref{cor_Toindep}, the theorem is reduced immediately 
to \cite{Il} (II.2.2.1).
\end{proof}

As an application we obtain the vanishing of the almost cotangent 
complex for a certain class of morphisms.

\begin{theorem}\label{th_vanish.L} Let $R\to S$ be a morphism 
of almost algebras such that 
$$\Tor^R_i(S,S)\simeq 0\simeq\Tor^{S\otimes_RS}_i(S,S)
\qquad\text{for all $i>0$}$$
(for the natural $S\otimes_RS$-module structure induced by 
$\mu_{S/R}$). Then $\L_{S/R}\simeq 0$ in $\sD_\bullet(S_{!!})$.
\end{theorem}
\begin{proof} Since $\Tor^R_i(S,S)=0$ for all $i>0$, theorem
\ref{th_flat.base.ch} applies (with $A=R$ and $B=A'=S$), giving 
the natural isomorphisms
\set\begin{equation}\label{eq_two.equa}
\begin{array}{l}
(S\otimes_RS)_{!!}\otimes_{S_{!!}}\L_{S/R}\simeq\L_{S\otimes_RS/S} \\
((S\otimes_RS)_{!!}\otimes_{S_{!!}}\L_{S/R})\oplus
((S\otimes_RS)_{!!}\otimes_{S_{!!}}\L_{S/R})\simeq\L_{S\otimes_RS/R}
\end{array}
\end{equation}
Since $\Tor^{S\otimes_RS}_i(S,S)=0$, the same theorem also applies
with $A=S\otimes_RS$, $B=S$, $A'=S$, and we notice that in this
case $B'\simeq S$; hence we have 
\set\begin{equation}\label{eq_thanksSR}
\L_{S/S\otimes_RS}\simeq 
S_{!!}\otimes_{S_{!!}}\L_{S/S\otimes_RS}\simeq\L_{S/S}\simeq 0.
\end{equation}
Next we apply transitivity to the sequence $R\to S\otimes_RS\to S$,
to obtain (thanks to \eqref{eq_thanksSR}) 
\set\begin{equation}\label{eq_one.down}
S_{!!}\otimes_{S\otimes_RS_{!!}}\L_{S\otimes_RS/R}\simeq\L_{S/R}.
\end{equation}
Applying $S_{!!}\otimes_{S\otimes_RS_{!!}}-$ to the second isomorphism
\eqref{eq_two.equa} we obtain 
\set\begin{equation}\label{eq_two.to.go}
\L_{S/R}\oplus\L_{S/R}\simeq
S_{!!}\otimes_{S\otimes_RS_{!!}}\L_{S\otimes_RS/R}.
\end{equation}
Finally, composing \eqref{eq_one.down} and \eqref{eq_two.to.go}
we derive 
\set\begin{equation}\label{eq_over}
\L_{S/R}\oplus\L_{S/R}\stackrel{\sim}{\to}\L_{S/R}.
\end{equation}
However, by inspection, the isomorphism \eqref{eq_over} is the
sum map. Consequently $\L_{S/R}\simeq 0$, as claimed.
\end{proof}

Finally we have a fundamental spectral sequence as in \cite{Il} 
(III.3.3.2). 
\begin{theorem}\label{th_main.spec.seq} 
Let $\phi:A\to B$ be a morphism of almost algebras such that 
$B\otimes_AB\simeq B$ ({\em e.g.} such that $B$ is a quotient of $A$). 
Then there is a first quadrant homology spectral sequence of bigraded 
almost algebras
$$E^2_{pq}=H_{p+q}(\mbox{\rm Sym}^q_B(\L^a_{B/A}))\Rightarrow
\Tor^{A}_{p+q}(B,B).$$
\end{theorem}
\begin{proof} We replace $\phi$ by $\one_{V^a}\times\phi$ and
apply the functor $B\mapsto B_{!!}$ (which commutes with tensor 
products by remark \ref{rem_adjoints}) thereby reducing 
the assertion to the above mentioned \cite{Il} (III.3.3.2).
\end{proof}

\section{Almost ring theory} 
 
\subsection{Flat, unramified and \'etale morphisms}
Let $A\to B$ be a morphism of almost $V$-algebras. Using the natural
``multiplication'' morphism of $A$-algebras $\mu_{B/A}:B\otimes_AB\to B$ 
we can view $B$ as a $B\otimes_AB$-algebra.

\begin{definition}\label{def_morph} 
Let $\phi:A\to B$ be a morphism of almost $V$-algebras.

i) We say that $\phi$ is a {\em flat\/} (resp. {\em faithfully flat\/},
resp. {\em almost projective\/}) {\em morphism\/} if 
$B$ is a flat (resp. faithfully flat, resp. almost projective) 
$A$-module.

ii) We say that $\phi$ is {\em almost finite\/} (resp. {\em finite\/})
if $B$ is an almost finitely generated (resp. finitely generated) 
$A$-module.

iii) We say that $\phi$ is {\em weakly unramified\/} (resp. 
{\em unramified\/}) if $B$ is a flat (resp. almost projective) 
$B\otimes_AB$-module (via the morphism $\mu_{B/A}$ defined above).

iv) $\phi$ is {\em weakly \'etale\/} (resp. {\em \'etale\/}) 
if it is flat and weakly unramified (resp. unramified).
\end{definition}

\begin{lemma}\label{lem_itoiv}
Let $\phi:A\to B$ and $\psi:B\to C$ be morphisms of almost 
$V$-algebras. 

i) Any base change of a flat (resp. almost projective, resp. 
faithfully flat, resp. almost finite, resp. weakly unramified, 
resp. unramified, resp. weakly \'etale, resp. \'etale) morphism 
is flat (resp. almost projective, resp. faithfully flat, resp. 
almost finite, resp. weakly unramified, resp. unramified, resp.
weakly \'etale, resp. \'etale);

ii) if both $\phi$ and $\psi$ are flat (resp. almost projective,
resp. faithfully flat, resp. almost finite, resp. weakly unramified,
resp. unramified, resp. weakly \'etale, resp. \'etale), then so is 
$\psi\circ\phi$;

iii) if $\phi$ is flat and $\psi\circ\phi$ is faithfully flat, then 
$\phi$ is faithfully flat;

iv) if $\phi$ is weakly unramified and $\psi\circ\phi$ is flat 
(resp. weakly \'etale), then $\psi$ is flat (resp. weakly \'etale);

v) If $\phi$ is unramified and $\psi\circ\phi$ is \'etale, then
$\psi$ is \'etale;

vi) $\phi$ is faithfully flat if and only if it is a monomorphism
and $B/A$ is a flat $A$-module;

vii) if $\phi$ is almost finite and weakly unramified, then $\phi$
is unramified.
\end{lemma}
\begin{proof} For (vi) use the Tor sequences. In view of proposition 
\ref{prop_converse}(ii), to show (vii) it suffices to know that $B$ 
is an almost finitely presented $B\otimes_AB$-module; but 
this follows from the existence of an epimorphism of 
$B\otimes_AB$-modules $(B\otimes_AB)\otimes_AB\to\Ker(\mu_{B/A})$ 
defined by
$x\otimes b\mapsto x\cdot(\ubarold 1\otimes b-b\otimes\ubarold 1)$.
Of the remaining assertions, only (iv) and (v) are not obvious, but 
the proof is just the ``almost version'' of a well-known argument. 
Let us show (v); the same argument applies to (iv). We remark that 
$\mu_{B/A}$ is an \'etale morphism, since $\phi$ is unramified. 
Define $\Gamma_\psi=\one_C\otimes_B\mu_{B/A}$. By (i), 
$\Gamma_\psi$ is \'etale. Define also 
$p=(\psi\circ\phi)\otimes_A\one_B$. By (i), $p$ is flat (resp. 
\'etale). The claim follows by remarking that 
$\psi=\Gamma_\psi\circ p$ and applying (ii).
\end{proof}

\begin{remark}\label{rem_naive} 
i) Suppose we work in the classical limit case, that is, $V=\fm$
(cp. example \ref{ex_rings}(ii)). Then we caution the reader that
our notion of ``\'etale morphism'' is more general than the usual 
one, as defined in \cite{SGA1}. The relationship between the usual 
notion and ours is discussed in the digression at the end of
section \ref{sec_descent}.

ii) The naive hope that the functor $A\mapsto A_{!!}$ 
might preserve flatness is crushed by the following counterexample.
Let $(V,\fm)$ be as in example \ref{ex_rings}(i) and let $k$ 
be the residue field of $V$. Consider the flat map 
$V\times V\to V$ defined as $(x,y)\mapsto x$. We get a flat morphism 
$V^a\times V^a\to V^a$ in $V^a\Alg$; applying the left adjoint 
to localisation yields a map $V\times_kV\to V$ that is not flat.
On the other hand, faithful flatness {\em is\/} preserved.
Indeed, let $\phi:A\to B$ be a morphism of almost algebras.
Then $\phi\/$ is a monomorphism if and only if $\phi_{!!}\/$
is injective; moreover, $B_{!!}/\Img(A_{!!})\simeq B_!/A_!$,
which is flat over $A_{!!}$ if and only if $B/A$ is flat over
$A$, by proposition \ref{prop_comp.supp}.
\end{remark}

We will find useful to study certain ``almost idempotents'', 
as in the following proposition.

\begin{proposition}\label{prop_idemp}
A morphism $\phi:A\to B$ is unramified if and only if there exists 
an almost element $e_{B/A}\in\ubar{B\otimes_AB}$ such that

i) $e_{B/A}^2=e_{B/A}$;

ii) $\mu_{B/A}(e_{B/A})=\ubarold 1$;

iii) $x\cdot e_{B/A}=0$ for all $x\in I_{B/A*}$.
\end{proposition}
\begin{proof} Suppose that $\phi$ is unramified. We start by 
showing that for every $\eps\in\fm$ there exist almost elements 
$e_\eps$ of $B\otimes_AB$ such that 
\set\begin{equation}\label{eq_esta}
e_\eps^2=\eps\cdot e_\eps\qquad
\mu_{B/A}(e_\eps)=\eps\cdot\ubarold 1\qquad
I_{B/A*}\cdot e_\eps=0.
\end{equation}
Since $B$ is an almost projective $B\otimes_AB$-module, for 
every $\eps\in\fm$ there exists an ``approximate splitting'' for the
epimorphism $\mu_{B/A}:B\otimes_AB\to B$, {\em i.e.\/} a 
$B\otimes_AB$-linear morphism $u_\eps:B\to B\otimes_AB$ such that
$\mu_{B/A}\circ u_\eps=\eps\cdot\one_B$. Set 
$e_\eps=u_\eps\circ\ubarold 1:A\to B\otimes_AB$.
We see that $\mu_{B/A}(e_\eps)=\eps\cdot\ubarold 1$. To show that 
$e^2_\eps=\eps\cdot e_\eps$ we use the $B\otimes_AB$-linearity 
of $u_\eps$ to compute
$$e^2_\eps=e_\eps\cdot u_\eps(\ubarold 1)=
u_\eps(\mu_{B/A}(e_\eps)\cdot\ubarold 1)=
u_\eps(\mu_{B/A}(e_\eps))=\eps\cdot e_\eps.$$
Next take any almost element $x$ of $I_{B/A}$ and compute
$$x\cdot e_\eps=x\cdot u_\eps(\ubarold 1)=
u_\eps(\mu_{B/A}(x)\cdot\ubarold 1)=0.$$
This establishes \eqref{eq_esta}. Next let us take any other
$\delta\in\fm$ and a corresponding almost element $e_\delta$. 
Both $\eps\cdot\ubarold 1-e_\eps$ and 
$\delta\cdot\ubarold 1-e_\delta$ are elements of $I_{B/A*}$, 
hence we have $(\delta\cdot\ubarold 1-e_\delta)\cdot e_\eps=0
=(\eps\cdot\ubarold 1-e_\eps)\cdot e_\delta$
which implies 
\set\begin{equation}\label{eq_fund.e}
\delta\cdot e_\eps=\eps\cdot e_\delta\qquad
\text{for all $\eps,\delta\in\fm$.}
\end{equation} 
Let us define a map $e_{B/A}:\fm\otimes_V\fm\to B\otimes_AB_*$ 
by the rule
\set\begin{equation}\label{eq_define.e}
\eps\otimes\delta\mapsto\delta\cdot e_\eps\qquad
\text{for all $\eps,\delta\in\fm$}.
\end{equation}
To show that \eqref{eq_define.e} does indeed determine a well defined
morphism, we need to check that 
$\delta\cdot v\cdot e_\eps=\delta\cdot e_{v\cdot\eps}$ and
$\delta\cdot e_{\eps+\eps'}=\delta\cdot(e_\eps+e_{\eps'})$
for all $\eps,\eps',\delta\in\fm$ and all $v\in V$. However,
both identities follow easily by a repeated application of 
\eqref{eq_fund.e}. It is easy to see that $e_{B/A}$ defines 
an almost element with the required properties.

Conversely, suppose an almost element $e_{B/A}$ of $B\otimes_AB$ 
is given with the stated properties. We define $u:B\to B\otimes_AB$ 
by $b\mapsto e_{B/A}\cdot(1\otimes b)$ ($b\in\ubar B$) and $v=\mu_{B/A}$. 
Then (iii) says that $u$ is a $B\otimes_AB$-linear morphism and (ii) 
shows that $v\circ u=\one_B$. Hence, by lemma \ref{lem_thetwo}, 
$\phi$ is unramified.
\end{proof}

\begin{corollary}\label{cor_unram} Under the hypotheses and notation 
of the proposition, the ideal $I=I_{B/A}$ has a natural structure 
of $A$-algebra, with unit morphism given
by $\ubarold 1_{I/A}=\ubarold 1_{B\otimes_AB/A}-e_{B/A}$ and whose 
multiplication is the restriction of $\mu_{B\otimes_AB/A}$ to $I$. 
Moreover the natural morphism
$$B\otimes_AB\to I_{B/A}\oplus B \qquad x\mapsto 
(x\cdot\ubarold 1_{I/A}\oplus \mu_{B/A}(x))$$
is an isomorphism of $A$-algebras.
\end{corollary}
\begin{proof} Left to the reader as an exercise.
\end{proof}

\subsection{Almost traces}\label{sect_traces}
Let $A$ be an almost $V$-algebra. For any integer $n>0$, the 
standard direct sum decomposition of $A^n$ determines uniquely
$A$-linear morphisms
$A\stackrel{e_i^A}{\longrightarrow} 
A^n\stackrel{\pi_j^A}{\longrightarrow} A$ (for $i,j=1,...,n$)
such that $\pi^A_j\circ e^A_i=\delta_{ij}\cdot\one_A$ for all $i,j$
and $\sum_{i=1}^ne^A_i\circ\pi_i^A=\one_{A^n}$.
We can then define a natural {\em trace homomorphism\/}
\set\begin{equation}\label{eq_def.trace}
\Tr:\Alhom_A(A^n,A^n)\to A\qquad 
\phi\mapsto\sum_{i=1}^n\pi^A_i\circ\phi\circ e^A_i
\end{equation}
which is an $A$-linear morphism. For any 
$\phi,\psi\in\Alhom_A(A^n,A^n)_*$ we have
$\Tr(\phi\circ\psi)=\Tr(\psi\circ\phi)$.
It follows easily that $\Tr$ is independent of the given direct 
sum decomposition of $A^n$.
More generally, suppose that $M$ is an almost projective almost 
finitely generated $A$-module. Then for any 
$\eps\in\fm$ we can find $n=n(\eps)$ and morphisms
\set\begin{equation}\label{eq_usual}
{\diagram
M \rto^{u_\eps} & A^n \rto^{v_\eps} & M
\enddiagram}
\end{equation}
such that $v_\eps\circ u_\eps=\eps\cdot\one_M$. Let 
$E(M)=\Alhom_A(M,M)$; notice that $E(M)_*$ is naturally
isomorphic to $\Hom_A(M,M)$. We consider the $A$-linear morphism
\set\begin{equation}\label{eq_approx_trace}
t_\eps:E(M)\to A \qquad \phi\mapsto\Tr(u_\eps\circ\phi\circ v_\eps)
\qquad\text{($\phi\in E(M)_*$).}
\end{equation}
Now, pick any other $\delta\in\fm$. We compute
$$\begin{array}{r@{\:=\:}l}
\eps\cdot t_\delta(\phi)&\eps\cdot\Tr(u_\delta\circ\phi\circ v_\delta)
=\Tr(u_\delta\circ v_\eps\circ u_\eps\circ\phi\circ v_\delta) \\
&\Tr(u_\eps\circ\phi\circ v_\delta\circ u_\delta\circ v_\eps)
=\delta\cdot\Tr(u_\eps\circ\phi\circ v_\eps) 
=\delta\cdot t_\eps(\phi). 
\end{array}$$
This allows us to define a map
$$t_{M/A}:\fm\otimes_V\fm\otimes_VE(M)_*\to A_*$$
by setting 
$\eps\otimes\delta\otimes\phi\mapsto\eps\cdot t_\delta(\phi)$.
We leave to the reader the verification that $t_{M/A}$ is 
well defined and does not depend on the choice of $t_\delta$.
It induces a morphism $E(M)\to A$ that we denote again by
$t_{M/A}$ and we call the {\em almost trace morphism} for the 
almost $A$-module $M$. 

Let $f\in\ubar M^*$, $m\in\ubar M$ and define $\phi_{f,m}\in E(M)_*$ 
by the formula 
$\phi_{f,m}(m')=f(m')\cdot m$ for all $m'\in\ubar M$. We have 
the following :
\begin{lemma}\label{lem_Gabber}
With the above notation : $t_{M/A}(\phi_{f,m})=f(m)$.
\end{lemma}
\begin{proof} Let $f:M\to A$ and $m:A\to M$ be given. Obviously
we have $\phi_{f,m}=m\circ f$ and $f(m)=f\circ m$. Pick morphisms
$u_\eps$ and $v_\eps$ as in \eqref{eq_usual}. Using the foregoing
notation, we can write :
$$\begin{array}{r@{\:=\:}l}
t_\eps(\phi_{f,m})&
\sum_{i=1}^n(\pi_i^A\circ u_\eps\circ m)\circ(f\circ v_\eps\circ e_i^A)\\
&\sum_{i=1}^n(f\circ v_\eps\circ e_i^A)\circ(\pi_i^A\circ u_\eps\circ m)\\
&f\circ v_\eps\circ u_\eps\circ m=\eps\cdot f\circ m
\end{array}$$
from which the claim follows directly.
\end{proof}

\begin{lemma}\label{lem_swap.maps} Let $M$ and $N$ be 
almost finitely generated almost projective $A$-modules, 
and $\phi:M\to N$, $\psi:N\to M$ two $A$-linear morphisms. 
Then :

i) $t_{M/A}(\psi\circ\phi)=t_{N/A}(\phi\circ\psi)$.

ii) If $\psi\circ\phi=a\cdot\one_M$ and 
$\phi\circ\psi=a\cdot\one_N$ for some $a\in A_*$, and
if, furthermore, there exist $u\in\End_A(M)$, $v\in\End_A(N)$
such that $v\circ\phi=\phi\circ u$, then 
$a\cdot(t_{M/A}(u)-t_{N/A}(v))=0$.
\end{lemma}
\begin{proof} (i) is left to the reader as an exercise.
For (ii) we compute using (i) : 
$a\cdot t_{M/A}(u)=t_{M/A}(\psi\circ\phi\circ u)=
t_{M/A}(\psi\circ v\circ\phi)=t_{N/A}(v\circ\phi\circ\psi)=
a\cdot t_{N/A}(v)$.
\end{proof}

\begin{proposition} Let 
$\ubarold M=(0\to M_1\stackrel{i}{\to}M_2\stackrel{p}{\to}M_3\to 0)$ 
be an exact sequence of almost finitely generated almost
projective $A$-modules, and let 
$\ubarold u=(u_1,u_2,u_3):\ubarold M\to\ubarold M$ be 
an endomorphism of $\ubarold M$, given by endomorphisms 
$u_i:M_i\to M_i$ ($i=1,2,3$). Then we have 
$t_{M_2/A}(u_2)=t_{M_1/A}(u_1)+t_{M_3/A}(u_3)$.
\end{proposition}
\begin{proof} Suppose first that there exists a splitting
$s:M_3\to M_2$ for $p$, so that we can view $u_2$ as a matrix 
$\left(\begin{array}{ll} u_1 & v \\ 0 & u_3\end{array}\right)$,
where $v\in\Hom_A(M_3,M_1)$. By additivity of the trace, we are 
then reduced to show that $t_{M_2/A}(i\circ v\circ p)=0$.
By lemma \ref{lem_swap.maps}(i), this is the same as 
$t_{M_3/A}(p\circ i\circ v)$, which obviously vanishes.
In general, for any $a\in\fm$ we consider the morphism
$\mu_a=a\cdot\one_{M_3}$ and the pull back morphism 
$\ubarold{M}*\mu_a\to\ubarold M$ : 
$$\xymatrix{
0 \ar[r] & M_1 \ar[r]^i & M_2 \ar[r]^p & M_3 \ar[r] & 0 \\
0 \ar[r] & M_1 \ar[r] \udouble & P \ar[r] \ar[u]^{\phi} &
M_3 \ar[r] \ar[u]^{\mu_a} & 0}$$
Then $\ubarold{M}*\mu_a$ is a split exact sequence
with the endomorphism $\ubarold{u}*\mu_a=(u_1,v,u_3)$,
for a certain $v\in\End_A(P)$. The pair of morphisms
$(a\cdot\one_{M_2},p)$ induces a morphism $\psi:M_2\to P$,
and it is easy to check that 
$\phi\circ\psi=a\cdot\one_{M_2}$ and 
$\psi\circ\phi=a\cdot\one_P$. We can therefore apply 
lemma \ref{lem_swap.maps} to deduce that 
$a\cdot(t_{P/A}(v)-t_{M/A}(u))=0$. By the foregoing 
we know that $t_{P/A}(v)=t_{M_1/A}(u_1)+t_{M_3/A}(u_3)$,
so the claim follows.
\end{proof}

Suppose now that $B$ is an almost finitely generated almost 
projective $A$-algebra. For any $b\in\ubar B$, denote 
by $\mu_b:B\to B$ the $B$-linear morphism $b'\mapsto b\cdot b'$. 
The map $b\mapsto\mu_b$ defines a $B$-linear monomorphism 
$\mu:B\to E(B)$. The composition 
$$\Tr_{B/A}=t_{B/A}\circ\mu:B\to A$$
will also be called the almost trace morphism of the $A$-algebra $B$.

\begin{proposition}\label{prop_traces}
Let $A$ and $B$ be as in the above discussion.

i) If $\phi:A\to B$ is an isomorphism, then $\Tr_{B/A}=\phi^{-1}$.

ii) If $C$ any other $A$-algebra, then 
$\Tr_{C\otimes_AB/C}=\one_C\otimes_A\Tr_{B/A}$.

iii) If $C$ is an almost projective almost finite 
$B$-algebra, then $\Tr_{C/A}=\Tr_{B/A}\circ\Tr_{C/B}$.
\end{proposition}
\begin{proof} (i) and (ii) are left as exercises for the reader.
We verify (iii). For given $\eps,\delta\in\fm$ pick morphisms
$B\stackrel{u_\eps}{\longrightarrow} 
A^n\stackrel{v_\eps}{\longrightarrow} B$ and 
$C\stackrel{u'_\delta}{\longrightarrow}
B^m\stackrel{v'_\delta}{\longrightarrow} C$
such that $v_\eps\circ u_\eps=\eps\cdot\one_B$ and
$v'_\delta\circ u'_\delta=\delta\cdot\one_C$. If we set 
$u_\eps^{\oplus m}=u_\eps\otimes_A\one_{A^m}$,
$u''_{\delta\eps}=u_\eps^{\oplus m}\circ u'_\delta:C\to A^n\otimes_AA^m$,
$v_\eps^{\oplus m}=v_\eps\otimes_A\one_{A^m}$
and 
$v''_{\delta\eps}=v'_\delta\circ v^{\oplus m}_\eps:A^n\otimes_AA^m\to C$
then we have 
$v''_{\delta\eps}\circ u''_{\delta\eps}=\eps\cdot\delta\cdot\one_C$.
Define 
$$\begin{array}{r@{\: :\:}l}
t_{\eps,B/A}&~B\to A\qquad b\mapsto\Tr(u_\eps\circ\mu_b\circ v_\eps) \\
t_{\delta,C/B}&~C\to B\qquad c\mapsto\Tr(u'_\delta\circ\mu_c\circ v'_\delta) \\
t_{\delta\eps,C/A}&~C\to A\qquad c\mapsto
\Tr(u''_{\delta\eps}\circ\mu_c\circ v''_{\delta\eps}).
\end{array}$$
Using \eqref{eq_def.trace} we can write
$$\begin{array}{r@{\:=\:}l}
t_{\delta\eps,C/A}(c)&
\displaystyle{\mathop\sum_{i,j=1}^{n,m}}(\pi^A_i\otimes_A\pi^A_j)\circ 
                          u''_{\delta\eps}\circ\mu_c\circ v''_{\delta\eps}
                          \circ(e^A_i\otimes_Ae^A_j) \\
&\displaystyle{\mathop\sum_{i,j=1}^{n,m}}
  (\pi^A_i\otimes\pi^A_j)\circ u_\eps^{\oplus m}
  \circ u'_\delta\circ\mu_c\circ v'_\delta\circ
  v_\eps^{\oplus m}\circ(e^A_i\otimes e^A_j) \\
&\displaystyle{\mathop\sum_{i,j=1}^{n,m}}
  \pi^A_i\circ u_\eps\circ\pi^B_j\circ u'_\delta\circ\mu_c
  \circ v'_\delta\circ e^B_j\circ v_\eps\circ e^A_i \\
&\displaystyle{\mathop\sum_{i=1}^n}
  \pi^A_i\circ u_\eps\circ t_{\delta,C/B}(c)\circ v_\eps\circ e^A_i \\
&t_{\eps,B/A}\circ t_{\delta,C/B}(c)
\end{array}$$
which implies immediately the claim.
\end{proof}
\begin{corollary} Let $A\to B$ be a faithfully flat almost 
finitely presented  and \'etale morphism of almost 
$V$-algebras. Then $\Tr_{B/A}:B\to A$ is an epimorphism.
\end{corollary}
\begin{proof} Under the stated hypotheses, $B$ is an almost 
projective $A$-module (by proposition \ref{prop_converse}). 
Let $C=\Coker(\Tr_{B/A})$ and $\Tr_{B/B\otimes_AB}$ the trace 
morphism for the morphism of almost $V$-algebras $\mu_{B/A}$. 
By faithful flatness, the natural morphism 
$C\to C\otimes_AB=\Coker(\Tr_{B\otimes_AB/B})$ is a monomorphism, 
hence it suffices to show that $\Tr_{B\otimes_AB/B}$ is an 
epimorphism (here $B\otimes_AB$ is considered as a
$B$-algebra via the second factor). However, from proposition 
\ref{prop_traces}(i) and (iii) we see that $\Tr_{B/B\otimes_AB}$ 
is a right inverse for $\Tr_{B\otimes_AB/B}$. The claim follows.
\end{proof}
It is useful to introduce the $A$-linear morphism
$$\tr_{B/A}=\Tr_{B/A}\circ\mu_{B/A}:B\otimes_AB\to A.$$
We can view $\tr_{B/A}$ as a bilinear form; it
induces an $A$-linear morphism
$$\tau_{B/A}:B\to B^*=\Alhom_A(B,A)$$
characterized by the equality 
$\tr_{B/A}(b_1\otimes b_2)=\tau_{B/A}(b_1)(b_2)$ for all
$b_1,b_2\in\ubar B$. We say that $\tr_{B/A}$ is
{\em a perfect pairing\/} if $\tau_{B/A}$ is an isomorphism.

\begin{theorem}\label{th_proj.etale} An almost projective and almost finite 
morphism $\phi:A\to B$ of almost $V$-algebras is \'etale if and only
if the trace form $\tr_{B/A}$ is a perfect pairing.
\end{theorem}
\begin{proof} Suppose that $\phi$ is \'etale. Let $e_{B/A}$
be the idempotent almost element of $B\otimes_AB$ provided by proposition
\ref{prop_idemp}. We define a morphism 
$\sigma:B^*\to B$ by $f\mapsto (f\otimes_A\one_B)(e_{B/A})$.
To start with, we remark that both $\tau_{B/A}$ and $\sigma$ are
$B$-linear morphisms (for the natural $B$-module
structure of $B^*$ defined in remark \ref{rem_B.struct}).
Indeed, let us pick any $b,b',b''\in\ubar B$, $f\in\ubar B^*$ 
and compute directly
$$\begin{array}{r@{\:=\:}l}
(b\cdot\tau_{B/A}(b'))(b'')&\tau_{B/A}(b')(bb'')=
\Tr_{B/A}(bb'b'')=(\tau_{B/A}(bb'))(b''). \\
b\cdot\sigma(f)&b\cdot(f\otimes_A\one_B)(e_{B/A})=
(f\otimes_A\one_B)((\ubarold 1_{B/A}\otimes b)\cdot e_{B/A}) \\
&(f\otimes_A\one_B)((b\otimes\ubarold 1_{B/A})\cdot e_{B/A})=
((b\cdot f)\otimes_A\one_B)(e_{B/A}) \\
&\sigma(b\cdot f).
\end{array}$$
Next we show that $\sigma$ is a left inverse for $\tau_{B/A}$. 
In fact, let $b\in\ubar B$. We have
$$\begin{array}{r@{\:=\:}l}
\sigma\circ\tau_{B/A}(b)&(\tau_{B/A}(b)\otimes_A\one_B)(e_{B/A})=
(\Tr_{B/A}\otimes_A\one_B)((b\otimes\ubarold 1_{B/A})\cdot e_{B/A}) \\
&\Tr_{B\otimes_AB/B}((\ubarold 1_{B/A}\otimes b)\cdot e_{B/A})
=b\cdot\Tr_{B\otimes_AB/B}(e_{B/A}).
\end{array}$$
Therefore it suffices to show that 
$\Tr_{B\otimes_AB/B}(e_{B/A})=\ubarold 1$. However, by hypothesis 
$\phi$ is unramified, hence corollary \ref{cor_unram} gives a 
decomposition $B\otimes_AB\simeq B\oplus I_{B/A}$ such that $e_{B/A}$ 
acts as the identity on the first factor and as the zero morphism on 
the second factor. 

Now, let $X=\Ker(\sigma)$. From the above we derive a $B$-linear
isomorphism $B^*\simeq B\oplus X$. We dualize and apply 
lemma \ref{lem_alhom}(ii) to obtain another $B$-linear isomorphism
\set\begin{equation}\label{eq_split.epi}
B\simeq(B^*)^*\simeq(B\oplus X)^*\simeq 
B^*\oplus X^*\simeq B\oplus X\oplus X^*.
\end{equation}
Finally, composing the isomorphism \eqref{eq_split.epi} with
the projection on the first factor, we get a split $B$-linear
epimorphism $B\to B$, hence a {\em surjective} 
$\ubar B$-linear morphism $\ubar B\to\ubar B$. Such a morphism
is necessarily an isomorphism, and, tracing backward, the same 
must hold for $\tau_{B/A}$.

Conversely, suppose that the trace form is a perfect pairing.
By lemma \ref{lem_alhom}(i) the natural morphism
$\alpha:B^*\otimes_AB\to\Alhom_B(B\otimes_AB,B)$ is an isomorphism 
and one verifies easily that 
$\alpha\circ(\tau_{B/A}\otimes_A\one_B)=\tau_{B\otimes_AB/B}$.
In particular $\tau_{B\otimes_AB/B}$ is also an isomorphism.
The multiplication gives an almost element 
$\mu_{B/A}\in\Alhom_B(B\otimes_AB,B)_*$; let 
$e=\tau_{B\otimes_AB/B}^{-1}(\mu_{B/A})$. We derive
\set\begin{equation}\label{eq_e.trace}
\Tr_{B\otimes_AB/B}(e)=
\tau_{B\otimes_AB/B}(e)(\ubarold 1_{B\otimes_AB})=
\mu_{B/A}(\ubarold 1_{B\otimes_AB})=\ubarold 1_{B/A}.
\end{equation}
Furthermore, we have already remarked that $\tau_{B/A}$ is a $B$-linear
morphism, hence $\tau_{B\otimes_AB/B}$ is a $B\otimes_AB$-linear morphism. 
Consequently, for any almost element $x$
of $B\otimes_AB$ we have
$$\tau_{B\otimes_AB/B}(x\cdot e)=x\cdot\tau_{B\otimes_AB/B}(e)
=x\cdot\mu_{B/A}=\mu_{B/A}(x)\cdot\mu_{B/A}=
\mu_{B/A}(x)\cdot\tau_{B\otimes_AB/B}(e).$$
Since by hypothesis $\tau_{B/A}$ is an isomorphism, this
implies 
\set\begin{equation}\label{eq_e.mu}
x\cdot e=\mu_{B/A}(x)\cdot e.
\end{equation}
Consider the morphism $\mu_e:B\otimes_AB\to B\otimes_AB$ defined by 
$x\mapsto e\cdot x$; then $\mu_e$ is $B$-linear (for the 
$B$-module structure defined by the second factor).
Applying \eqref{eq_e.mu} and lemma \ref{lem_Gabber} we conclude 
that $t_{B\otimes_AB/B}(\mu_e)=\mu_{B/A}(e)$. On the other hand, 
\eqref{eq_e.trace} says that this trace is equal to $\ubarold 1_{B/A}$, 
hence 
\set\begin{equation}\label{eq_mu_e}
\mu_{B/A}(e)=\ubarold 1_{B/A}.
\end{equation}
Let $\beta:B\to B\otimes_AB$ be defined as $b\mapsto b\cdot e.$ From 
\eqref{eq_e.mu} we see that both $\beta$ and $\mu_{B/A}$ are 
$B\otimes_AB$-linear morphisms and from \eqref{eq_mu_e} we know 
moreover that $\mu_{B/A}\circ\beta=\one_B$.
By lemma \ref{lem_thetwo} we deduce that $B$ is an almost projective 
$B\otimes_AB$-module, {\em i.e.\/} $\phi$ is unramified, as claimed.
\end{proof}

\begin{definition} 
The {\em nilradical\/} of an almost algebra $A$ is the ideal 
{\em $\nil(A)=\nil(\ubar A)^a$} (where, for a ring $R$, we denote by
{\em $\nil(R)$} the ideal of nilpotent elements in $R$).
We say that $A$ is {\em reduced\/} if {\em $\nil(A)\simeq 0$}.
\end{definition}

Notice that, if $R$ is a $V$-algebra, then every nilpotent 
ideal in $R^a$ is of the form $I^a$, where $I$ is a nilpotent
ideal in $R$ (indeed, it is of the form $I^a$ where $I$ is
an ideal, and $\fm\cdot I$ is seen to be nilpotent). It
follows easily that $\nil(A)$ is the colimit of the nilpotent
ideals in $A$; moreover $\nil(R)^a=\nil(R^a)$. Using this
one sees that $A/\nil(A)$ is reduced.

\begin{proposition} Let $A\to B$ be an \'etale almost finitely
presented morphism of almost algebras. If $A$ is reduced then $B$ 
is reduced as well.
\end{proposition}
\begin{proof} Under the stated hypothesis, $B$ is an almost projective
$A$-module (by virtue of proposition \ref{prop_converse}(ii)). Hence, 
for given $\eps\in\fm$, pick a sequence of morphisms 
$B\stackrel{u_\eps}{\to}A^n\stackrel{v_\eps}{\to}B$ such that
$v_\eps\circ u_\eps=\eps\cdot\one_B$; with the notation
of \eqref{eq_approx_trace}, define $\nu_b:A^n\to A^n$
by $\nu_b=v_\eps\circ\mu_b\circ u_\eps$, so that 
$t_\eps(b)=\Tr(\nu_b)$. One verifies easily that 
$\nu^m_b=\eps^{m-1}\cdot\nu_{b^m}$ for all integers $m>0$.

Now, suppose that $b\in\nil(\ubar B)$. It follows that $b^m=0$ for $m$ 
sufficiently large, hence $\nu^m_b=0$ for $m$ sufficiently large.
Let $\fp$ be any prime ideal of $\ubar A$; let $\pi:\ubar A\to\ubar A/\fp$
be the natural projection and $F$ the fraction field of $\ubar A/\fp$. 
The $F$-linear morphism $\ubar{\nu_b}\otimes_{\ubar A}\one_F$ is 
nilpotent on the vector space $F^n$, hence 
$\pi\circ\Tr(\ubar{\nu_b})=\Tr(\ubar{\nu_b}\otimes_{\ubar A}\one_F)=0$.
This shows that $\Tr(\ubar{\nu_b})$ lies in the intersection of all
prime ideals of $\ubar A$, hence it is nilpotent. Since by hypothesis
$A$ is reduced, we get $\Tr(\ubar{\nu_b})=0$. Finally, this implies
$\Tr_{B/A}(b)=0$.  Now, for any $b'\in\ubar B$, the almost element $bb'$
will be nilpotent as well, so the same conclusion applies to it.
This shows that $\tau_{B/A}(b)=0$. But by hypothesis $B$ is \'etale
over $A$, hence theorem \ref{th_proj.etale} yields $b=0$, as required.
\end{proof}
\begin{remark} Let $M$ be an $A$-module. We say that an almost 
element $a$ of $A$ is $M$-regular if the multiplication morphism 
$m\mapsto am~:~M\to M$ is a monomorphism. Assume ({\bf A}) (cf.
section \ref{sec_ring.prel}) and suppose furthermore that $\fm$ 
is generated by a multiplicative system $\cS$ which is a cofiltered 
semigroup under the preorder structure $(\cS,\succ)$ induced by 
the divisibility relation in $V$. We say that $\cS$ is archimedean 
if, for all $s,t\in\cS$ there exists $n>0$ such that $s^n\succ t$. 
Suppose that $\cS$ is archimedean and that $A$ is a reduced 
almost algebra. Then $\cS$ consists of $A$-regular elements. 
Indeed, by hypothesis $\nil(\ubar A)^a=0$; since the annihilator 
of $\cS$ in $\ubar A$ is $0$ we get $\nil(\ubar A)=0$. Suppose 
that $s\cdot a=0$ for some $s\in\cS$ and $a\in\ubar A$. Let 
$t\in\cS$ be arbitrary and pick $n>0$ such that $t^n\succ s$. 
Then $(ta)^n=0$ hence $ta=0$ for all $t\in\cS$, hence $a=0$.
\end{remark}

\subsection{Lifting theorems}\label{sec_lift}
Throughout the following, the terminology ``epimorphism of 
$V^a$-alge\-bras'' will refer to a morphism of $V^a$-algebras
that induces an epimorphism on the underlying $V^a$-modules. 
\begin{lemma}\label{lem_epim}
Let $A\to B$ be an epimorphism of almost $V$-algebras with 
kernel $I$. Let $U$ be the $A$-extension $0\to I/I^2\to A/I^2\to B\to 0$.
Then the assignment $f\mapsto f*U$ defines a natural isomorphism
\set\begin{equation}\label{eq_Hom.Exal}
{\diagram
\Hom_B(I/I^2,M) \rto^{\sim} & \Exal_A(B,M).
\enddiagram}
\end{equation}
\end{lemma}
\begin{proof} Let $X=(0\to M\to E\stackrel{p}{\to}B\to 0)$ be any 
$A$-extension of $B$ by $M$. The composition $g:A\to E\stackrel{p}{\to}B$
of the structural morphism for $E$ followed by $p$ coincides with
the projection $A\to B$. Therefore $g(I)\subset M$ and $g(I^2)=0$.
Hence $g$ factors through $A/I^2$; the restriction of $g$ to $I/I^2$
defines a morphism $f\in\Hom_B(I/I^2,M)$ and a morphism of $A$-extensions
$f*U\to X$. In this way we obtain an inverse for \eqref{eq_Hom.Exal}.
\end{proof}
Now consider any morphism of $A$-extensions 
\set\begin{equation}\label{eq_A.ext}{
\diagram 
~\tilde B: \dto^{\tilde f} & \qquad 0 \rto & I \rto \dto^u & 
B \rto \dto^f & B_0 \rto \dto^{f_0} & 0 \\
~\tilde C: & \qquad 0 \rto & J \rto & C \rto & C_0 \rto & 0.
\enddiagram
}\end{equation}
The morphism $u$ induces by adjunction a morphism of $C_0$-modules
\set\begin{equation}\label{eq_adj} 
C_0\otimes_{B_0}I\to J
\end{equation}
whose image is the ideal $I\cdot C$, so that the square diagram of almost 
algebras defined by $\tilde f$ is cofibred ({\em i.e.\/} 
$C_0\simeq C\otimes_BB_0$) if and only if \eqref{eq_adj} is an epimorphism.
\begin{lemma}\label{lem_flat}
Let $\tilde f:\tilde B\to\tilde C$ be a morphism of $A$-extensions as above, 
such that the corresponding square diagram of almost algebras is cofibred. 
Then the morphism $f:B\to C$ is flat if and only if $f_0:B_0\to C_0$ is flat 
and \eqref{eq_adj} is an isomorphism.
\end{lemma}
\begin{proof} It follows directly from the (almost version of the) local
flatness criterion (see \cite{Mat} Th. 22.3).
\end{proof}

We are now ready to put together all the work done so far and
begin the study of deformations of almost algebras. 
\smallskip

The morphism $u:I\to J$ is an element in $\Hom_{B_0}(I,J)$; by lemma 
\ref{lem_epim} the latter group is naturally isomorphic to 
$\Exal_B(B_0,J)$. 
By applying transitivity (theorem \ref{th_transit}) to the sequence of 
morphisms $B\to B_0\stackrel{f_0}{\to} C_0$ we obtain an exact sequence 
of abelian groups
$$\Exal_{B_0}(C_0,J)\to\Exal_B(C_0,J)\to\Hom_{B_0}(I,J)
\stackrel{\partial}{\to}\hExt^2_{C_{0!!}}(\L_{C_0/B_0},J_!).$$
Hence we can form the element
$\omega(\tilde B,f_0,u)=\partial(u)\in
\hExt^2_{C_{0!!}}(\L_{C_0/B_0},J_!)$.
The proof of the next result goes exactly as in \cite{Il} (III.2.1.2.3).
\begin{proposition}\label{prop_obstr}
i) Let the $A$-extension $\tilde B$, the $B_0$-linear morphism $u:I\to J$ 
and the morphism of $A$-algebras $f_0:B_0\to C_0$ be given as above.
Then there exists an $A$-extension $\tilde C$ and a morphism 
$\tilde f:\tilde B\to\tilde C$ completing diagram \eqref{eq_A.ext} if and 
only if $\omega(\tilde B,f_0,u)=0$. (\/{\em i.e.\/} $\omega(\tilde B,f_0,u)$
is the obstruction to the lifting of $\tilde B$ over $f_0$.)

ii) Assume that the obstruction $\omega(\tilde B,f_0,u)$ vanishes. 
Then the set of isomorphism classes of $A$-extensions $\tilde C$ as 
in (i) forms a torsor under the group $\Exal_{B_0}(C_0,J)$ 
(\/$\simeq\hExt^1_{C_{0!!}}(\L_{C_0/B_0},J_!)$\/).

iii) The group of automorphisms of an $A$-extension $\tilde C$ as 
in (i) is naturally isomorphic to $\Der_{B_0}(C_0,J)$ 
(\/$\simeq\hExt^0_{C_{0!!}}(\L_{C_0/B_0},J_!)$\/).
\qed\end{proposition} 
The obstruction $\omega(\tilde B,f_0,u)$ depends functorially on $u$.
More exactly, if we denote by
$$\omega(\tilde B,f_0)\in
\hExt^2_{C_{0!!}}(\L_{C_0/B_0},(C_0\otimes_{B_0}I)_!)$$
the obstruction corresponding to the natural morphism 
$I\to C_0\otimes_{B_0}I$, then for any other morphism $u:I\to J$ we have 
$$\omega(\tilde B,f_0,u)=v_!\circ\omega(\tilde B,f_0)$$
where $v$ is the morphism \eqref{eq_adj}. Taking lemma \ref{lem_flat} 
into account we deduce 
\begin{corollary}\label{cor_flat.obstr}
Suppose that $B_0\to C_0$ is flat. Then

i) The class $\omega(\tilde B,f_0)$ is the obstruction to the existence
of a flat deformation of\/ $C_0$ over $B$, {\em i.e.\/} of a $B$-extension
$\tilde C$ as in \eqref{eq_A.ext} such that $C$ is flat over $B$ and 
$C\otimes_BB_0\to C_0$ is an isomorphism.

ii) If the obstruction $\omega(\tilde B,f_0)$ vanishes, then the set of 
isomorphism classes of flat deformations of\/ $C_0$ over $B$ forms a torsor 
under the group {\em $\Exal_{B_0}(C_0,C_0\otimes_{B_0}I)$}.

iii) The group of automorphisms of a given flat deformation of\/ $C_0$
over $B$ is naturally isomorphic to {\em $\Der_{B_0}(C_0,C_0\otimes_{B_0}I)$}.
\qed\end{corollary}

Now, suppose we are given two $A$-extensions $\tilde C^1,\tilde C^2$ 
with morphisms of $A$-extensions 
$$\diagram 
~\tilde B: \dto^{\tilde f^i} & \qquad 0 \rto & I \rto \dto^{u^i} & 
B \rto \dto^{f^i} & B_0 \rto \dto^{f^i_0} & 0 \\
~\tilde C^i: & \qquad 0 \rto & J^i \rto & C^i \rto & C^i_0 \rto & 0
\enddiagram$$
and morphisms $v:J^1\to J^2$, $g_0:C_0^1\to C_0^2$ such that 
\set\begin{equation}\label{eq_problem}
u^2=v\circ u^1 \qquad \text{and} \qquad f_0^2=g_0\circ f_0^1. 
\end{equation}
We consider the problem of finding a morphism of $A$-extensions 
\set\begin{equation}\label{eq_A.ext.C}{
\diagram
~\tilde C^1: \dto^{\tilde g} & \qquad 0 \rto & J^1 \rto \dto^v & 
C^1 \rto \dto^g & C^1_0 \rto \dto^{g_0}& 0 \\
~\tilde C^2: & \qquad 0 \rto & J^2 \rto & C^2 \rto & C^2_0 \rto & 0
\enddiagram
}\end{equation}
such that $\tilde f^2=\tilde g\circ\tilde f^1$. Let us denote 
by $e(\tilde C^i)\in\hExt^1_{C^i_{0!!}}(\L_{C^i_0/B},J^i_!)$ 
the classes defined by the $B$-extensions $\tilde C^1,\tilde C^2$ 
via the isomorphism of theorem \ref{th_main} and by
$$\begin{array}{r@{~:~}l}
v* & \hExt^1_{C^1_{0!!}}(\L_{C^1_0/B},J^1_!)\to
\hExt^1_{C^1_{0!!}}(\L_{C^1_0/B},J^2_!) \\
*g_0 & \hExt^1_{C^2_{0!!}}(\L_{C^2_0/B},J^2_!)\to
\hExt^1_{C^2_{0!!}}(C^2_{0!!}\otimes_{C^1_{0!!}}\L_{C^1_0/B},J^2_!)
\end{array}$$
the canonical morphisms defined by $v$ and $g_0$. Using the natural
isomorphism
$$\hExt^1_{C^1_{0!!}}(\L_{C^1_0/B},J^2_!)\simeq
\hExt^1_{C^2_{0!!}}(C^2_{0!!}\otimes_{C^1_{0!!}}\L_{C^1_0/B},J^2_!)$$
we can identify the target of both $v*$ and $*g$ with 
$\hExt^1_{C^1_{0!!}}(\L_{C^1_0/B},J^2_!)$.
It is clear that the problem
admits a solution if and only if the $A$-extensions $v*\tilde C^1$ 
and $\tilde C^2*g_0$ coincide, {\em i.e.\/} if and only if 
$v*e(\tilde C^1)-e(\tilde C^2)*g_0=0$. By applying transitivity 
to the sequence of morphisms $B\to B_0\to C^1_0$ we obtain an 
exact sequence
$$\hExt^1_{C^1_{0!!}}(\L_{C^1_0/B_0},J^2_!)\hookrightarrow
\hExt^1_{C^1_{0!!}}(\L_{C^1_0/B},J^2_!)
\to\Hom_{C^1_0}(C^1_0\otimes_{B_0}I,J^2)$$
It follows from \eqref{eq_problem} that the image of 
$v*e(\tilde C^1)-e(\tilde C^2)*g_0$ in the group 
$\Hom_{C^1_0}(C^1_0\otimes_{B_0}I,J^2)$ vanishes, therefore 
\set\begin{equation}\label{eq_obstr}
v*e(\tilde C^1)-e(\tilde C^2)*g_0\in
\hExt^1_{C^1_{0!!}}(\L_{C^1_0/B_0},J^2_!).
\end{equation}
In conclusion, we derive the following result as in \cite{Il} (III.2.2.2).
\begin{proposition}\label{prop_obstr2}
With the above notations, the class \eqref{eq_obstr} is 
the obstruction to the existence of a morphism of $A$-extensions 
$\tilde g:\tilde C^1\to\tilde C^2$ as in \eqref{eq_A.ext.C} such that
$\tilde f^2=\tilde g\circ\tilde f^1$. When the obstruction vanishes, the 
set of such morphisms forms a torsor under the group 
{\em $\Der_{B_0}(C^1_0,J^2)$} (the latter being identified with 
$\hExt^0_{C^2_{0!!}}(C^2_{0!!}\otimes_{C^1_{0!!}}\L_{C^1_0/B_0},J^2_!)$\/).
\qed\end{proposition}

For a given almost $V$-algebra $A$, we define the category $\wEt(A)$ 
as the full subcategory of $A\Alg$ consisting of all weakly \'etale
$A$-algebras. Notice that, by lemma \ref{lem_itoiv}(iv) all morphisms 
in $\wEt(A)$ are weakly \'etale.        

\begin{theorem}\label{th_liftetale}
i) Let $A\to B$ be a weakly \'etale morphism 
of almost algebras. Let $C$ be any $A$-algebra and 
$I\subset C$ a nilpotent ideal. Then the natural morphism
$$\Hom_{A\Alg}(B,C)\to\Hom_{A\Alg}(B,C/I)$$
is bijective.

ii) Let $A$ be a $V^a$-algebra, $I\subset A$ a nilpotent 
ideal and $A'=A/I$. Then the natural functor
$$\wEt(A)\to\wEt(A')\qquad
(\phi:A\to B)\mapsto(\one_{A'}\otimes_A\phi:A'\to A'\otimes_AB)$$
is an equivalence of categories.

iii) The equivalence of (ii) restricts to an equivalence 
$\Et(A)\to\Et(A')$.
\end{theorem}
\begin{proof}By induction we can assume $I^2=0$. Then (i) follows 
directly from proposition \ref{prop_obstr2} and theorem 
\ref{th_vanish.L}. We show (ii) : by corollary \ref{cor_flat.obstr} 
(and again theorem \ref{th_vanish.L}) a given weakly \'etale 
morphism $\phi':A'\to B'$ can be lifted to a {\em unique\/} 
flat morphism $\phi:A\to B$. We need to prove that $\phi$ is 
weakly \'etale, {\em i.e.\/} that $B$ is $B\otimes_AB$-flat. 
However, it is clear that $\mu_{B'/A'}:B'\otimes_{A'}B'\to B'$ 
is weakly \'etale, hence it has a flat lifting 
$\tilde\mu:B\otimes_AB\to C$. Then the composition 
$A\to B\otimes_AB\to C$ is flat and it is a lifting 
of $\phi'$. We deduce that there is an isomorphism of 
$A$-algebras $\alpha:B\to C$ lifting $\one_{B'}$ and moreover
the morphisms $b\mapsto\tilde\mu(b\otimes\ubarold 1)$ and
$b\mapsto\tilde\mu(\ubarold 1\otimes b)$ coincide with $\alpha$. 
Claim (ii) follows. To show (iii), suppose that $A'\to B'$ is \'etale
and let $I_{B'/A'}$ denote as usual the kernel of $\mu_{B'/A'}$.  
By corollary \ref{cor_unram} there is a natural morphism of 
almost algebras $B'\otimes_{A'}B'\to I_{B'/A'}$ which is clearly 
\'etale. Hence $I_{B'/A'}$ lifts to a weakly \'etale
$B\otimes_AB$-algebra $C$, and the isomorphism
$B'\otimes_{A'}B'\simeq I_{B'/A'}\oplus B'$ lifts to an isomorphism
$B\otimes_AB\simeq C\oplus B$ of $B\otimes_AB$-algebras.
It follows that $B$ is an almost projective 
$B\otimes_AB$-module, {\em i.e.\/} $A\to B$ is \'etale, as claimed.
\end{proof}

We conclude with some results on deformations of almost modules.
These can be established independently of the theory of the 
cotangent complex, along the lines of \cite{Il} (IV.3.1.12).

We begin by recalling some notation from {\em loc. cit.\/}
Let $R$ be a ring and $J\subset R$ an ideal with
$J^2=0$. Set $R'=R/J$; an extension of $R$-modules
$\ubarold M=(0\to K\to M\stackrel{p}{\to} M'\to 0)$ 
where $K$ and $M'$ are killed by $J$, defines a natural morphism
of $R'$-modules $u(\ubarold M):J\otimes_{R'}M'\to K$
such that $u(\ubarold M)(x\otimes m')=xm$ for $x\in J$, $m\in M$
and $p(m)=m'$. By the local flatness criterion (\cite{Mat} Th.22.3)
$M$ is flat over $R$ if and only if $M'$ is flat over $R'$ and
$u(\ubarold M)$ is an isomorphism. One can then show the following.
\begin{proposition}\label{prop_liftmodules}(cp. \cite{Il} (IV.3.1.5)) 

i) Given $R'$-modules $M'$ and $K$ and a morphism 
$u':J\otimes_{R'}M'\to K$ there exists an obstruction 
$\omega(R,u')\in\Ext_{R'}^2(M',K)$ whose vanishing is necessary 
and sufficient for the existence of an extension of $R$-modules 
$\ubarold M$ of $M'$ by $K$ such that $u(\ubarold M)=u'$. 

ii) When $\omega(R,u')=0$, the set of isomorphism classes of 
such extensions $\ubarold M$ forms a torsor under 
$\Ext^1_{R'}(M',K)$; the group of automorphisms of such
an extension is $\Hom_{R'}(M',K)$.
\qed\end{proposition}

\begin{lemma}\label{lem_frequent}
Let $A\to B$ be a finite morphism of almost algebras with nilpotent
kernel. Let $\phi:M\to N$ be an $A$-linear morphism and set 
$\phi_B=\phi\otimes_A\one_B:
M\otimes_AB\to N\otimes_AB$. Then there exists $m\geq 0$ such that

i) $\Ann_A(\Coker(\phi_B))^m\subset\Ann_A(\Coker(\phi))$.

ii) $(\Ann_V(\Ker(\phi_B))\cdot\Ann_V(\Tor^A_1(B,N))\cdot
\Ann_V(\Coker(\phi)))^m\subset\Ann_A(\Ker(\phi))$.

\noindent If $B=A/I$ for some nilpotent ideal $I$, and $I^n=0$, then
we can take $m=n$ in (i) and (ii).
\end{lemma}
\begin{proof} Under the assumptions, we can find a finitely 
generated $A_*$-module $Q$ such that $\fm\cdot B_*\subset Q\subset B_*$.
By \cite{Gru} (1.1.5), there exists a finite filtration 
$0=J_m\subset ...\subset J_1\subset J_0=A_*$ such that each $J_i/J_{i+1}$
is a quotient of a direct sum of copies of $Q$. This implies that, for
every $A$-module $M$, we have 
\set\begin{equation}\label{eq_Annihil}
\Ann_A(M\otimes_AB)^m\subset\Ann_A(M).
\end{equation}
(i) follows easily. Notice that if $B=A/I$ and $I^n=0$, then we can 
take $m=n$ in \eqref{eq_Annihil}. For (ii) let $C^\bullet=\Cone(\phi)$. 
We estimate $H=H^{-1}(C^\bullet\derotimes_AB)$ in two ways. 
By the first spectral sequence of hyperhomology we have an
exact sequence $\Tor_1^A(N,B)\to H\to\Ker(\phi_B)$.
By the second spectral sequence for hyperhomology 
we have an exact sequence 
$\Tor^A_2(\Coker(\phi),B)\to\Ker(\phi)\otimes_AB\to H$.
Hence $\Ker(\phi)\otimes_AB$ is annihilated by the product 
of the three annihilators in (ii) and the result follows by 
applying \eqref{eq_Annihil} with $M=\Ker(\phi)$.
\end{proof}

\begin{lemma}\label{lem_flat.proj} Keep the assumptions of lemma 
\ref{lem_frequent}, let $M$ be an $A$-module and set 
$M_B=B\otimes_AM$.

i) If $A\to B$ is an epimorphism, $M$ is flat and $M_B$ is almost 
projective over $B$, then $M$ is almost projective over $A$.

ii) If $M_B$ is an almost finitely generated $B$-module
then $M$ is an almost finitely generated $A$-module.

iii) If $\Tor_1^A(B,M)=0$ and $M_B$ is almost finitely presented
over $B$, then $M$ is almost finitely presented over $A$.
\end{lemma}
\begin{proof} (i) : we have to show that $\Ext^1_A(M,N)$ is almost
zero for every $A$-module $N$. Let $I=\Ker(A\to B)$; by 
assumption $I$ is nilpotent, so by the usual devissage we may 
assume that $I\cdot N=0$. If $\chi\in\Ext^1_A(M,N)$ 
is represented by an extension $0\to N\to Q\to M\to 0$ then 
after tensoring by $B$ and using the flatness of $M$ we get 
an exact sequence of $B$-modules 
$0\to N\to B\otimes_AQ\to M_B\to 0$. Thus $\chi$ comes from an 
element of $\Ext^1_{B}(M_B,N)$ which is almost zero by assumption. 

(ii) : let $\fm_0=(\eps_1,...,\eps_m)$ be a finitely generated 
subideal of $\fm$. By assumption there is a map 
$\phi':B^r\to M_B$ such that $\fm_0\cdot\Coker(\phi')=0$. 
For all $j\le m$ the morphism $\eps_j\cdot\phi'$ lifts to a 
morphism $\phi_j:A^r\to M$. Then 
$\phi=\phi_1\oplus...\oplus\phi_m:A^{rm}\to M$ satisfies
$\fm_0^2\cdot\Coker(\phi\otimes_A\one_B)=0$. By lemma 
\ref{lem_frequent}(i) it follows $\fm_0^{2n}\cdot\Coker(\phi)=0$ 
for some $n\geq 0$. As $\fm_0$ was arbitrary, the result follows.

(iii) Let $\fm_0$ be as above. By (ii), $M$ is almost 
finitely generated over $A$, so we can choose a morphism
$\phi:A^r\to M$ such that $\fm_0\cdot\Coker(\phi)=0$.
Consider $\phi_B=\phi\otimes_A\one_B:B^r\to M_B$.
By lemma \ref{lem_fin.pres}, there is a finitely generated
submodule $N$ of $\Ker(\phi_B)$ containing
$\fm_0^2\cdot\Ker(\phi_B)$. Notice that $\Ker(\phi)\otimes_AB$
maps onto $\Ker(B^r\to\Img(\phi)\otimes_AB)$ and
$\Ker(\Img(\phi)\otimes_AB\to M_B)\simeq
\Tor^A_1(B,\Coker(\phi))$ is annihilated by $\fm_0$.
Hence $\fm_0\cdot\Ker(\phi_B)$ is contained in the image
of $\Ker(\phi)$ and therefore we can lift a finite generating 
set $\{x_1',...,x_n'\}$ for $\fm_0^2\cdot N$ to almost 
elements $\{x_1,...,x_n\}$ of $\Ker(\phi)$. If we quotient 
$A^r$ by the span of these $x_i$, we get a finitely presented 
$A$-module $F$ with a morphism $\bar\phi:F\to M$
such that $\Ker(\bar\phi\otimes_AB)$ is annihilated by 
$\fm_0^4$ and $\Coker(\bar\phi)$ is annihilated by $\fm_0$.
By lemma \ref{lem_frequent}(ii) we derive 
$\fm_0^{5m}\cdot\Ker(\bar\phi)=0$ for some $m\geq 0$. 
Since $\fm_0$ is arbitrary, this proves the result.
\end{proof}

\begin{remark}\label{rem_descent} 
(i) Inspecting the proof, one sees that parts (ii)
and (iii) of lemma \ref{lem_flat.proj} hold whenever 
\eqref{eq_Annihil} holds. For instance, if $A\to B$ is any 
faithfully flat morphism, then \eqref{eq_Annihil} holds 
with $m=1$. 

ii) Consequently, if $A\to B$ is faithfully flat
and $M$ is an $A$-module such that $M_B$ is flat 
(resp. almost finitely generated, resp. almost finitely 
presented) over $B$, then $M$ is flat (resp. almost finitely
generated, resp. almost finitely presented) over $A$.

iii) On the other hand, we do not know whether a general
faithfully flat morphism $A\to B$ descends almost projectivity.
However, using (ii) and proposition \ref{prop_converse}
we see that if the $B$-module $M_B$ is almost finitely 
generated almost projective, then $M$ has the same property.

iv) However, if $B$ is faithfully flat and almost finitely
presented as an $A$-module, then $A\to B$ does descend 
almost projectivity, as can be easily deduced from lemma 
\ref{lem_alhom}(i) and proposition \ref{prop_converse}(ii).
\end{remark}

\begin{theorem}\label{th_liftmod} Let $I$ be a nilpotent ideal 
of the almost algebra $A$ and set $A'=A/I$. Suppose that $\tilde\fm$ 
is a (flat) $V$-module of homological dimension $\leq 1$. 
Let $P'$ be an almost projective $A'$-module.

i) There is an almost projective $A$-module $P$ 
with $A'\otimes_AP\simeq P'$. 

ii) If $P'$ is almost finitely presented, then $P$ is almost 
finitely presented.
\end{theorem}
\begin{proof} As usual we reduce to $I^2=0$. Then proposition 
\ref{prop_liftmodules}(i) applies with $R=A_*$, $J=I_*$, 
$R'=A_*/I_*$, $M'=P'_!$, $K=I_*\otimes_{R'}P'_!$ and $u'=\one_K$. 
We obtain a class 
$\omega(A_*,u')\in\Ext^2_{R'}(P'_!,I_*\otimes_{R'}P'_!)$ which 
gives the obstruction to the existence of a flat $A_*$-module 
$F$ lifting $P'_!$. Since $P'_!$ is almost projective, we know that
$\fm\cdot\Ext^2_{R'}(P'_!,I_*\otimes_{R'}P'_!)=0$, which says that
$0=\eps\cdot\omega(A_*,u')=\omega(A_*,\eps\cdot u')$ for all 
$\eps\in\fm$. In other words, for every $\eps\in\fm$ we can
find an extension of $A_*$-modules $\ubarold{P_\eps}$ of
$P'_!$ by $I_*\otimes_{R'}P'_!$ such that 
$u(\ubarold{P_\eps})=\eps\cdot\one_{I_*\otimes_{R'}P'_!}$. 
Let $\chi_\eps\in\Ext^1_{A_*}(P'_!,I_*\otimes_{R'}P'_!)$ be the 
class of $\ubarold{P_\eps}$. Notice that, for any $\delta\in\fm$,
$\delta\cdot\chi_\eps$ is the class of an extension $\ubarold X$
such that $u(\ubarold X)=\delta\cdot u(\ubarold{P_\eps})=
\delta\cdot\eps\cdot\one_{I_*\otimes_{R'}P'_!}$, hence, by 
proposition \ref{prop_liftmodules}(ii), 
$\gamma\cdot(\delta\cdot\chi_\eps-\chi_{\delta\cdot\eps})=0$
for all $\gamma\in\fm$. Hence we can define a morphism
$$\chi:\fm\otimes_V\fm\otimes_V\fm\to
\Ext^1_{A_*}(P'_!,I_*\otimes_{R'}P'_!)\qquad
\eps\otimes\delta\otimes\gamma\mapsto
\delta\cdot\gamma\cdot\chi_\eps.$$
However, one sees easily that 
$\fm\otimes_V\fm\otimes_V\fm\simeq\tilde\fm$
and $\tilde\fm\otimes_VP'_!\simeq P'_!$, hence we can view 
$\chi$ as an element of 
$\Hom_V(\tilde\fm,\Ext^1_{A_*}(P'_!,I_*\otimes_{R'}P'_!))$ 
and moreover we have a spectral sequence
$$E_2^{pq}=\Ext^p_V(\tilde\fm,\Ext^q_{A_*}(P'_!,I_*\otimes_{R'}P'_!))
\Rightarrow\Ext^{p+q}_{A_*}(P'_!,I_*\otimes_{R'}P'_!)$$
with $E_2^{pq}=0$ for all $p\geq 2$ (this spectral sequence
is constructed {\em e.g.\/} from the double complex 
$\Hom_V(F_p,\Hom_{A_*}(F'_q,I_*\otimes_{R'}P'_!))$ where
$F_\bullet$ (resp. $F'_\bullet$) is a projective resolution
of $\tilde\fm$ (resp. $P'_!$)). In particular, our $\chi$
is an element in $E_2^{01}$ which therefore survives in the
abutment as a class of $E_\infty^{01}$. The latter can be lifted
to an element $\tilde\chi$ via the surjection 
$\Ext^1_{A_*}(P'_!,I_*\otimes_{R'}P'_!)\to E_\infty^{01}$.
Let $0\to I_*\otimes_{R'}P'_!\to Q\to P'_!\to 0$ be an extension
representing $\tilde\chi$. Checking compatibilities, we see that 
$\delta\cdot\eps\cdot\tilde\chi=\delta\cdot\chi_\eps$ for every 
$\eps,\delta\in\fm$. Hence 
$u(\tilde\chi):I_*\otimes_{R'}P'_!\to I_*\otimes_{R'}P'_!$
coincides with the identity map on the submodule 
$\fm\cdot I_*\otimes_{R'}P'_!$. Since $\fm\cdot P_!=P_!$, we
see that $u(\tilde\chi)$ is actually the identity map.
By the local flatness criterion, it then follows
that $Q$ is flat over $R$, hence the $A$-module 
$P=Q^a$ is a flat lifting of $P'$, so it is almost projective, 
by lemma \ref{lem_flat.proj}(i). Now (ii) follows from
(i), lemma \ref{lem_flat.proj}(ii) and proposition 
\ref{prop_converse}(i).
\end{proof}

\begin{remark}\label{rem_counter} (i) According to proposition
\ref{prop_countably.pres}(ii), theorem \ref{th_liftmod} applies
especially when $\fm$ is countably generated as a $V$-module.

(ii) For $P$ and $P'$ as in theorem 
\ref{th_liftmod}(ii) let $\sigma_P:P\to P'$ be the projection. It is 
natural to ask whether the pair $(P,\sigma_P)$ is uniquely 
determined up to isomorphism, {\em i.e.\/} whether, for any other 
pair $(Q,\sigma_Q:Q\to P')$ for which theorem \ref{th_liftmod} 
holds, there exists an $A$-linear isomorphism $\phi:P\to Q$ such 
that $\sigma_Q\circ\phi=\sigma_P$. 
The answer is negative in general. Consider the case $P'=A'$. Take
$P=Q=A$ and let $\sigma_P$ be the natural projection, while 
$\sigma_Q=(u'\cdot\one_{A'})\circ\sigma_P$, where $u'$ is a unit 
in $\ubar A'$. Then the uniqueness question amounts to whether 
every unit in $\ubar A'$ lifts to a unit of $\ubar A$. The 
following counterexample is related to the fact that the completion 
of the algebraic closure $\bar\Q_p$ of $\Q_p$ is not maximally 
complete. Let $V=\bar\Z_p$, the integral closure of $\Z_p$ in 
$\bar\Q_p$. Then $V$ is a non-discrete valuation ring of rank 
one, and we take $\fm$ as in example \ref{ex_rings}(i), 
$A=(V/p^2V)^a$ and $A'=A/pA$. Choose a compatible system of 
roots of $p$. An almost element of $A'$ is just a $V$-linear 
morphism $\phi:\colim{n>0}p^{1/n!}V\to V/pV$.
Such a $\phi$ can be represented (in a non-unique way) by an 
infinite series of the form $\sum^\infty_{n=1}a_np^{1-1/n!}$
($a_n\in V$). The meaning of this expression is as follows. 
For every $m>0$, scalar multiplication by the element 
$\sum^{m}_{n=1}a_np^{1-1/n!}\in V$ defines a morphism 
$\phi_m:p^{1/m!}V\to V/pV$. For $m'>m$, let 
$j_{m,m'}:p^{1/m!}V\to p^{1/m'!}V$ be the imbedding. Then we
have $\phi_{m'}\circ j_{m,m'}=\phi_m$, so that we can define 
$\phi=\colim{m>0}\phi_m$. Similarly, every almost element of 
$A$ can be represented by an expression of the form 
$a_0+\sum^\infty_{n=1}a_np^{2-1/n!}$. Now, if $\sigma:A\to A'$
is the natural projection, the induced map 
$\sigma_*:A_*\to A'_*$ is given by:
$a_0+\sum^\infty_{n=1}a_np^{2-1/n!}\mapsto a_0$.
In particular, its image is the subring $V/p\subset(V/p)_*=A'_*$.
For instance, the unit $\sum^\infty_{n=1}p^{1-1/n!}$ of $A'_*$
does not lie in the image of this map.
\end{remark}
In the light of the above remark, the best one can achieve in general
is the following result.
\begin{proposition}\label{prop_best} Assume ({\bf A}) (see 
section \ref{sec_ring.prel}) and keep the notation of theorem 
\ref{th_liftmod}. Suppose moreover that $(Q,\sigma_Q:Q\to P')$ 
and $(P,\sigma_P:P\to P')$ are two pairs as in remark 
\ref{rem_counter}. Then for all $\eps\in\fm$ there exist 
$A$-linear morphisms $t_\eps:P\to Q$ and $s_\eps:Q\to P$ 
such that 
\smallskip

\noindent{\bf PQ}\rm{(}$\eps$\rm{)}$\qquad\qquad\qquad
\begin{array}{r@{\:=\:}lr@{\:=\:}l}
\sigma_Q\circ t_\eps & \eps\cdot\sigma_P & \qquad
\sigma_P\circ s_\eps & \eps\cdot\sigma_Q \\
s_\eps\circ t_\eps & \eps^2\cdot\one_P & \qquad
t_\eps\circ s_\eps & \eps^2\cdot\one_Q.
\end{array}$
\end{proposition}
\begin{proof}
Since both $Q$ and $P$ are almost projective and $\sigma_P,\sigma_Q$ 
are epimorphisms, there exist morphisms $\bar t_\eps:P\to Q$ and 
$\bar s_\eps:Q\to P$ such that 
$\sigma_Q\circ\bar t_\eps=\eps\cdot\sigma_P$ and 
$\sigma_P\circ\bar s_\eps=\eps\cdot\sigma_Q$. Then we have 
$\sigma_P\circ(\bar s_\eps\circ\bar t_\eps-\eps^2\cdot\one_P)=0$
and
$\sigma_Q\circ(\bar t_\eps\circ\bar s_\eps-\eps^2\circ\one_Q)=0$,
{\em i.e.\/} the morphism 
$u_\eps=\eps^2\cdot\one_P-\bar s_\eps\circ\bar t_\eps$ (resp.
$v_\eps=\eps^2\cdot\one_Q-\bar t_\eps\circ\bar s_\eps$) has image 
contained in the almost submodule $IP$ (resp. $IQ$). Since $I^m=0$ 
this implies $u^m_\eps=0$ and $v^m_\eps=0$. Hence
$$\eps^{2m}\cdot\one_P=(\eps^2\one_P)^m-u^m_\eps=
(\sum_{a=0}^{m-1}\eps^{2a}u^{m-1-a}_\eps)\circ
\bar s_\eps\circ\bar t_\eps.$$
Define $\bar s_{(2m-1)\eps}=
(\sum_{a=0}^{m-1}\eps^{2a}u^{m-1-a}_\eps)\circ\bar s_\eps$. 
Notice that  
$\bar s_{(2m-1)\eps}=
\bar s_\eps\circ(\sum_{a=0}^{m-1}\eps^{2a}v^{m-1-a}_\eps).$
This implies the equalities
$\bar s_{(2m-1)\eps}\circ\bar t_\eps=\eps^{2m}\cdot\one_P$ and
$\bar t_\eps\circ\bar s_{(2m-1)\eps}=\eps^{2m}\cdot\one_Q$.
Then the pair $(\bar s_{(2m-1)\eps},\eps^{2(m-1)}\cdot\bar t_\eps)$
satisfies {\bf PQ}($\eps^{2m-1}$). Under ({\bf A}), every element
of $\fm$ is a multiple of an element of the form $\eps^{2m-1}$, 
therefore the claim follows for arbitrary $\eps\in\fm$.
\end{proof}

\subsection{Descent}\label{sec_descent}
Faithfully flat descent in the almost setting presents no particular
surprises: since the functor $A\mapsto A_{!!}$ preserves faithful
flatness of morphisms (see remark \ref{rem_naive}) many well-known 
results for usual rings and modules extend {\em verbatim\/}
to almost algebras. So for instance, faithfully flat morphisms are
of universal effective descent for the fibred categories 
$F:V^a\AlgMod^o\to V^a\Alg^o$ and $G:V^a\AlgMorph^o\to V^a\Alg^o$
(see definition \ref{def_derivation}: for an almost $V$-algebra $B$, 
the fibre $F_B$ (resp. $G_B$) is the opposite of the category of 
$B$-modules (resp. $B$-algebras)).
Then, using remark \ref{rem_descent}, we deduce also universal 
effective descent for the fibred subcategories of flat (resp. almost 
finitely generated, resp. almost finitely presented, resp. almost 
projective almost finitely generated) modules.
Likewise, a faithfully flat morphism is of universal effective
descent for the fibred subcategories $\Et^o\to V^a\Alg^o$ of \'etale 
(resp. $\wEt^o\to V^a\Alg^o$ of weakly \'etale) algebras.

More generally, since the functor $A\mapsto A_{!!}$ preserves pure
morphisms in the sense of \cite{Oli}, and since, by a theorem of
Olivier ({\em loc. cit.\/}), pure morphisms are of universal effective 
descent for modules, the same holds for pure morphisms of almost algebras.

Non-flat descent is more delicate. Our results are not as complete
here as it could be wished, but nevertheless, they suffice for
current applications (namely, for the cases needed in \cite{Fa2}).

Our first statement is the almost version of a theorem of Gruson
and Raynaud (cp. \cite{Gr-Ra} (Part II, Th.1.2.4)).

\begin{proposition}\label{prop_descent}
A finite monomorphism of almost algebras descends flatness.
\end{proposition}
\begin{proof} Let $\phi:A\to B$ be such a morphism. Under the 
assumption, we can find a finite $A_*$-module $Q$ such that 
$\fm\cdot B_*\subset Q\subset B_*$. One sees easily that
$Q$ is a faithful $A_*$-module, so by \cite{Gr-Ra} (Part II, 
Th.1.2.4 and lemma 1.2.2), $Q$ satisfies the following condition :
%
\set\begin{equation}\label{eq_cond.Q}{
\parbox{14cm}{If\/ $(0\to N\to L\to P\to 0)$ is an exact sequence of 
$A_*$-modules with $L$ flat, such that $\Img(N\otimes_{A_*}Q)$ is 
a pure submodule of $L\otimes_{A_*}Q$, then $P$ is flat.}
}\end{equation}
%
Now let $M$ be an $A$-module such that $M\otimes_AB$ is flat.
Pick an epimorphism $p:F\to M$ with $F$ free over $A$.
Then ${\ubarold Y}=
(0\to\Ker(p\otimes_A\one_B)\to F\otimes_AB\to M\otimes_AB\to 0)$
is universally exact over $B$, hence over $A$. Consider the 
sequence 
${\ubarold X}=(0\to\Img(\Ker(p)_!\otimes_{A_*}Q)\to F_!\otimes_{A_*}Q
\to M_!\otimes_{A_*}Q\to 0)$.
Clearly ${\ubarold X}^a\simeq{\ubarold Y}$. However, it is easy to
check that a sequence $\ubarold E$ of $A$-modules is universally
exact if and only if the sequence ${\ubarold E}_!$ is universally exact
over $A_*$. We conclude that $\ubarold{X}=({\ubarold X}^a)_!$ is 
a universally exact sequence of $A_*$-modules, hence, by condition 
\eqref{eq_cond.Q}, $M_!$ is flat over $A_*$, {\em i.e.\/} $M$ is flat 
over $A$ as required. 
\end{proof}

\begin{corollary} Let $A\to B$ be a finite morphism of almost algebras,
with nilpotent kernel. If $C$ is a flat $A$-algebras such that 
$C\otimes_AB$ is weakly \'etale (resp. \'etale) over $B$, then $C$ is
weakly \'etale (resp. \'etale) over $A$.
\end{corollary}
\begin{proof} In the weakly \'etale case, we have to show that the
multiplication morphism $\mu:C\otimes_AC\to C$ is flat. As $N=\Ker(A\to B)$
is nilpotent, the local flatness criterion reduces the question
to the situation over $A/N$. So we may assume that $A\to B$ is a
monomorphism. Then $C\otimes_AC\to(C\otimes_AC)\otimes_AB$ is a
monomorphism, but $\mu\otimes_{C\otimes_AC}\one_{(C\otimes_AC)\otimes_AB}$
is the multiplication morphism of $C\otimes_AB$, which is flat
by assumption. Therefore, by proposition \ref{prop_descent}, $\mu$ 
is flat.

For the \'etale case, we have to show that $C$ is almost finitely 
presented as a $C\otimes_AC$-module. By hypothesis $C\otimes_AB$
is almost finitely presented as a $C\otimes_AC\otimes_AB$-module
and we know already that $C$ is flat as a $C\otimes_AC$-module, so 
by lemma \ref{lem_flat.proj}(iii) (applied to the finite morphism 
$C\otimes_AC\to C\otimes_AC\otimes_AB$) the claim follows. 
\end{proof}

Next we consider the following situation. We are given
a cartesian diagram of almost algebras
\set\begin{equation}\label{eq_cartesian}{
\diagram
A_0 \ar[r]^{f_2} \ar[d]_{f_1} & A_2 \ar[d]^{g_2} \\
A_1 \ar[r]^{g_1} & A_3
\enddiagram}\end{equation}
such that one of the morphisms $A_i\to A_3$ ($i=1,2$) is an 
epimorphism. We denote by $\cM_i$ (resp. $\cM_{i,\rm{fl}}$, resp. 
$\cM_{i,\rm{proj}}$) the category of all (resp. flat, resp. almost 
projective) $A_i$-modules, for $i=0,...,3$. Diagram \eqref{eq_cartesian} 
induces an essentially commutative diagram for the corresponding 
categories $\cM_i$, where the arrows are given by the ``extension 
of scalars'' functors. There follows a natural functor 
$$\pi:\cM_0\to\cM_1\times_{\cM_3}\cM_2$$
from $\cM_0$ to the 2-fibred products of $\cM_1$ and $\cM_2$ over
$\cM_3$. Recall (cp.~\cite{Bass} (Ch.VII \S 3)) that 
$\cM_1\times_{\cM_3}\cM_2$ is the category whose objects are the 
triples $(M_1,M_2,\xi)$, where $M_i$ is an $A_i$-module 
($i=1,2$) and 
$\xi:A_3\otimes_{A_1}M_1\stackrel{\sim}{\to}A_3\otimes_{A_2}M_2$
is an $A_3$-linear isomorphism. Given such an object $(M_1,M_2,\xi)$,
let us denote $M_3=A_3\otimes_{A_2}M_2$; we have a natural morphism
$M_2\to M_3$, and $\xi$ gives a morphism $M_1\to M_3$, so we can
form the fibre product $T(M_1,M_2,\xi)=M_1\times_{M_3}M_2$.
In this way we obtain a functor $T:\cM_1\times_{\cM_3}\cM_2\to\cM_0$,
and we leave to the reader the verification that $T$ is right
adjoint to $\pi$. Let us denote by $\eps:\one_{\cM_0}\to T\circ\pi$
and $\eta:\pi\circ T\to\one_{\cM_1\times_{\cM_3}\cM_2}$ the unit and
counit of the adjunction.

\begin{lemma}\label{lem_adjunct} The functor $\pi$ induces 
an equivalence of full subcategories :
$$\xymatrix{
\{X\in\rm{Ob}(\cM_0) | \eps_X\ \text{\rm{is an isomorphism}}\}
\ar[r]^-\pi &
\{Y\in\rm{Ob}(\cM_1\times_{\cM_3}\cM_2) | \eta_Y\ 
\text{\rm{is an isomorphism}}\}
}$$
having $T$ as essential inverse.
\end{lemma}
\begin{proof} General nonsense.
\end{proof}

\begin{lemma}\label{lem_eps.epi} Let $M$ be any $A_0$-module. 
Then $\eps_M$ is an epimorphism.  If $M$ is flat over $A_0$, $\eps_M$ 
is an isomorphism.
\end{lemma}
\begin{proof} Indeed, 
$\eps_M:M\to
(A_1\otimes_{A_0}M)\times_{A_3\otimes_{A_0}M}(A_2\otimes_{A_0}M)$
is the natural morphism. So, the assertions follow by applying 
$-\otimes_{A_0}M$ to the short exact sequence of $A_0$-modules
\set\begin{equation}\label{eq_short.ex.A_0}
0\to A_0\stackrel{f}{\to} A_1\oplus A_2\stackrel{g}{\to} A_3\to 0
\end{equation}
where $f(a)=(f_1(a),f_2(a))$ and $g(a,b)=g_1(a)-g_2(b)$.
\end{proof}
There is another case of interest, in which $\eps_M$ is an isomorphism.
Namely, suppose that one of the morphisms $A_i\to A_3$ ($i=1,2$),
say $A_1\to A_3$, has a section. Then also the morphism 
$A_0\to A_2$ gains a section $s:A_2\to A_0$ and we have the
following :
\begin{lemma}\label{lem_special.case} In the above situation, 
suppose that the $A_0$-module $M$ arises by extension 
of scalars from an $A_2$-module $M'$, via the section 
$s:A_2\to A_0$. Then $\eps_M$ is an isomorphism.
\end{lemma}
\begin{proof} Indeed, in this case, \eqref{eq_short.ex.A_0} 
is split exact as a sequence of $A_2$-modules, and it
remains such after tensoring by $M'$.
\end{proof}
 
\begin{lemma}\label{lem_key} $\eta_{(M_1,M_2,\xi)}$ is an 
isomorphism for all objects $(M_1,M_2,\xi)$.
\end{lemma}
\begin{proof} To fix ideas, suppose that $A_1\to A_3$ is an 
epimorphism. Consider any object $(M_1,M_2,\xi)$ 
of $\cM_1\times_{\cM_3}\cM_2$. Let $M=T(M_1,M_2,\xi)$; we deduce
a natural morphism 
$$\phi:(M\otimes_{A_0}A_1)\times_{M\otimes_{A_0}A_3}(M\otimes_{A_0}A_2)
\to M_1\times_{M_3}M_2$$ 
such that $\phi\circ\eps_M=\one_M$. It follows that $\eps_M$ is 
injective, hence it is an isomorphism, by lemma \ref{lem_eps.epi}.
We derive a commutative diagram with exact rows :
$$\xymatrix{
0 \ar[r] & M \ar[r] \ddouble & 
(M\otimes_{A_0}A_1)\oplus(M\otimes_{A_0}A_2) \ar[d]^{\phi_1\oplus\phi_2}
\ar[r] & M\otimes_{A_0}A_3 \ar[r] \ar[d]^{\phi_3} & 0 \\
0 \ar[r] & M \ar[r] & M_1\oplus M_2 \ar[r] & M_3 \ar[r] & 0.
}$$ 
From the snake lemma we deduce 
$$\begin{array}{ll}
(*)\qquad & \Ker(\phi_1)\oplus\Ker(\phi_2)\simeq\Ker(\phi_3) \\
(**)\qquad & \Coker(\phi_1)\oplus\Coker(\phi_2)\simeq\Coker(\phi_3).
\end{array}$$
Since $M_3\simeq M_1\otimes_{A_1}A_3$ we have 
$A_3\otimes_{A_1}\Coker(\phi_1)\simeq\Coker(\phi_3)$. But by 
assumption $A_1\to A_3$ is an epimorphism, so 
also $\Coker(\phi_1)\to\Coker(\phi_3)$ is an epimorphism. Then 
$(**)$ implies that $\Coker(\phi_2)=0$. But 
$\phi_3=\one_{A_3}\otimes_{A_2}\phi_2$, thus $\Coker(\phi_3)=0$
as well. We look at the exact sequence 
$0\to\Ker(\phi_1)\to M\otimes_{A_0}A_1\stackrel{\phi_1}{\to} M_1\to 0$ :
applying $A_3\otimes_{A_1}-$ we obtain an epimorphism 
$A_3\otimes_{A_1}\Ker(\phi_1)\to\Ker(\phi_3)$. From $(*)$ it follows
that $\Ker(\phi_2)=0$. In conclusion, $\phi_2$ is an isomorphism.
Hence the same is true for $\phi_3=\one_{A_3}\otimes_{A_2}\phi_2$,
and again $(*)$, $(**)$ show that $\phi_1$ is an isomorphism as well, 
which implies the claim.
\end{proof}

\begin{lemma}\label{lem_Ferrand} 
If $(A_1\times A_2)\otimes_{A_0}M$ is flat over $A_1\times A_2$, 
then $M$ is flat over $A_0$.
\end{lemma}
\begin{proof} Suppose that $A_1\to A_3$ is an epimorphism and let 
$I$ be its kernel. Let $\tilde A=A_{1!!}\times_{A_{3!!}}A_{2!!}$; 
it suffices to show that $M_!$ is a flat $\tilde A$-module. However, 
in view of proposition \ref{prop_comp.supp}, the assumption implies 
that $(A_{1!!}\times A_{2!!})\otimes_{\tilde A}M_!$ is a flat 
$A_{1!!}\times A_{2!!}$-module. $I_!$ is the kernel of the 
epimorphism $A_{1!!}\to A_{3!!}$. Moreover, $I_!$ identifies naturally 
with an ideal of $\tilde A$ and $\tilde A/I_!\simeq A_{2!!}$. Then the 
desired conclusion follows from \cite{Fer} (lemma in {\em loc. cit.\/}). 
\end{proof}

\begin{proposition}\label{prop_equiv.2-prod} The functor $\pi$ 
restricts to equivalences :
$$\begin{array}{c}
\cM_{0,\rm{fl}}\stackrel{\sim}{\to}
\cM_{1,\rm{fl}}\times_{\cM_{3,\rm{fl}}}\cM_{2,\rm{fl}} \\
\cM_{0,\rm{proj}}\stackrel{\sim}{\to}
\cM_{1,\rm{proj}}\times_{\cM_{3,\rm{proj}}}\cM_{2,\rm{proj}}.
\end{array}$$
\end{proposition}
\begin{proof} The assertion for flat almost modules follows 
directly from lemmata \ref{lem_adjunct}, \ref{lem_eps.epi}, 
\ref{lem_key} and \ref{lem_Ferrand}. Set $B=A_1\times A_2$. 
To establish the second equivalence, it suffices to show that, 
if $P$ is an $A_0$-module such that $B\otimes_{A_0}P$ 
is almost projective over $B$, then $P$ is almost projective 
over $A_0$, or which is the same, that $\AlExt^i_{A_0}(P,N)\simeq 0$ 
for all $i>0$ and any $A_0$-module $N$. We know already that 
$P$ is flat. Let $M$ be any $A_0$-module and $N$ any 
$B$-module. The standard isomorphism 
$R\Hom_B(B\derotimes_{A_0}M,N)\simeq R\Hom_{A_0}(M,N)$ yields 
a natural isomorphism 
$\AlExt^i_B(B\otimes_{A_0}M,N)\simeq\AlExt^i_{A_0}(M,N)$,
whenever $\Tor^{A_0}_j(B,M)=0$ for every $j>0$. In particular,
we have $\AlExt^i_{A_0}(P,N)\simeq 0$ whenever $N$ comes from 
either an $A_1$-module, or an $A_2$-module.
For a general $A_0$-module $N$ there is a 3-step filtration 
such that $\text{Fil}_0(N)=0$, 
$\gr_1(N)=\text{Fil}_1(N)=\Ker(\eps_N)$, 
$\gr_2(N)=\Ker(A_1\otimes_{A_0}N\to A_3\otimes_{A_0}N)$ 
and $\gr_3(N)=A_2\otimes_{A_0}N$. By an easy devissage,
we reduce to verify that $\AlExt^i_{A_0}(P,\gr_j(N))=0$
for every $i>0$ and $j=1,2,3$. However, $\gr_2(N)$ is an 
$A_1$-module and $\gr_3(N)$ is an $A_2$-module,
so the required vanishing follows for $j=2,3$. Moreover,
applying $-\otimes_{A_0}N$ to \eqref{eq_short.ex.A_0}, 
we derive a short exact sequence :
\set\begin{equation}\label{eq_gr1}
0\to\Tor^{A_0}_1(N,A_2)\to
\frac{\Tor^{A_0}_1(N,A_3)}{\Tor^{A_0}_1(N,A_1)}\to
\gr_1(N)\to 0.
\end{equation}
Here again, the leftmost term of \eqref{eq_gr1} is an
$A_2$-module, and the middle term is an $A_1$-module, 
so the same devissage yields the sought vanishing
also for $j=1$.
\end{proof}
\begin{corollary}\label{cor_equiv.2-prod} In the situation 
of \eqref{eq_cartesian}, denote by $\cA_{i,\rm{fl}}$ (resp. 
$\Et_i$, resp. $\wEt_i$) the category of flat (resp. \'etale, 
resp. weakly \'etale) $A_i$-algebras. The functor $\pi$ induces 
equivalences

$\cA_{0,\rm{fl}}\stackrel{\sim}{\to}
\cA_{1,\rm{fl}}\times_{\cA_{3,\rm{fl}}}\cA_{2,\rm{fl}}\qquad
\Et_0\stackrel{\sim}{\to}\Et_1\times_{\Et_3}\Et_2\qquad
\wEt_0\stackrel{\sim}{\to}\wEt_1\times_{\wEt_3}\wEt_2.
$\qed\end{corollary}

Next we want to reinterpret the equivalences of proposition
\ref{prop_equiv.2-prod} in terms of descent data. 
If $F:\cC\to V^a\Alg^o$ is a fibred category over the opposite
of the category of almost algebras, and if $X\to Y$ is a
given morphism of almost algebras, we shall denote by 
$\Desc(\cC,Y/X)$ the category of objects of the fibre category
$F_Y$, endowed with a descent datum relative to the morphism
$X\to Y$ (cp. \cite{Gi} (Ch.II \S1)). In the arguments hereafter, 
we consider morphisms of almost algebras and modules, and one 
has to reverse the direction of the arrows to pass to morphisms 
in the considered fibred category. 
Denote by $p_i:Y\to Y\otimes_XY$ ($i=1,2$), resp. 
$p_{ij}:Y\otimes_XY\to Y\otimes_XY\otimes_XY$ ($1\leq i<j\leq 3$) 
the usual morphisms. As an example, $\Desc(V^a\AlgMod^o,Y/X)$ 
consists of the pairs $(M,\beta)$ where $M$ is a $Y$-module 
and $\beta$ is a $Y\otimes_XY$-linear isomorphism 
$\beta:p_2^*(M)\stackrel{\sim}{\to}p_1^*(M)$ such that 
\set\begin{equation}\label{eq_cocycle}
p_{12}^*(\beta)\circ p_{23}^*(\beta)=p_{13}^*(\beta).
\end{equation} 
Let now $I\subset X$ be an ideal, and set $\bar X=X/I$, 
$\bar Y=Y/I\cdot Y$. For any $F:\cC\to V^a\Alg^o$ as above, 
one has an essentially commutative diagram:
$$\xymatrix{
\Desc(\cC,Y/X) \ar[r] \ar[d] & 
\Desc(\cC,\bar Y/\bar X) \ar[d] \\
F_Y \ar[r] & F_{\bar Y}.
}$$
This induces a functor :
\set\begin{equation}\label{eq_Desc}
\Desc(\cC,Y/X)\to
\Desc(\cC,\bar Y/\bar X)\times_{F_{\bar Y}}F_Y.
\end{equation}
\begin{lemma}\label{lem_Desc} With the above notation, suppose
moreover that  the natural morphism $I\to I\cdot Y$ is an isomorphism. 
Then the functor \eqref{eq_Desc} is an equivalence whenever
$\cC$ is one of the fibred categories $V^a\AlgMod^o$, 
$V^a\AlgMorph^o$, $\Et^o$, $\wEt^o$.
\end{lemma}
\begin{proof} For any $n>0$, denote by $Y^{\otimes n}$ (resp.
$\bar Y^{\otimes n}$) the $n$-fold tensor product of 
$Y$ (resp. $\bar Y$) with itself over $X$ (resp. $\bar X$), 
and by $\rho_n:Y^{\otimes n}\to\bar Y^{\otimes n}$ the natural
morphism. First of all we claim that, for every $n>0$, the 
natural diagram of almost algebras
\set\begin{equation}\label{eq_Y-cart}
{\diagram
Y^{\otimes n} \ar[r]^-{\rho_n} \ar[d]_{\mu_n} & 
\bar Y^{\otimes n} \ar[d]^{\bar\mu_n} \\
Y \ar[r]^-{\rho_1} & \bar Y
\enddiagram}\end{equation}
is cartesian (where $\mu_n$ and $\bar\mu_n$  are $n$-fold 
multiplication morphisms). For this, we need to verify that, 
for every $n>0$, the induced morphism $\Ker(\rho_n)\to\Ker(\rho_1)$
(defined by multiplication of the first two factors)
is an isomorphism. It then suffices to check that the
natural morphism $\Ker(\rho_n)\to\Ker(\rho_{n-1})$
is an isomorphism for all $n>1$. Indeed, consider the 
commutative diagram
$$\xymatrix{
I\otimes_XY^{\otimes n-1} \ddouble \ar[r]^-p & 
I\cdot Y^{\otimes n-1} \ar[r]^-i \ar[d]_-\psi & 
Y^{\otimes n-1} \ar[d]_{\phi\otimes\one_{Y^{\otimes n-1}}} 
\ar[rrd]^{\one_{Y^{\otimes n-1}}}  \\
I\otimes_XY^{\otimes n-1} \ar[r]_-{p'} & 
\Ker(\rho_n) \ar[r]_-{i'} & 
Y^{\otimes n} \ar[rr]_-{\mu_{Y/X}\otimes\one_{Y^{\otimes n-2}}} 
& & Y^{\otimes n-1}
}$$
From $I\cdot Y=\phi(Y)$, it follows that $p'$ is an epimorphism.
Hence also $\psi$ is an epimorphism. Since $i$ is a monomorphism,
it follows that $\psi$ is also a monomorphism, hence $\psi$
is an isomorphism and the claim follows easily.

We consider first the case $\cC=V^a\AlgMod^o$; we see that 
\eqref{eq_Y-cart} is a diagram of the kind considered in 
\eqref{eq_cartesian}, hence, for every $n>0$, we have the 
associated functor 
$\pi_n:Y^{\otimes n}\Mod\to
\bar Y^{\otimes n}\Mod\times_{\bar Y\Mod}Y\Mod$
and also its right adjoint $T_n$.
Denote by 
$\bar p_i:\bar Y\to\bar Y^{\otimes 2}$ ($i=1,2$) the usual 
morphisms, and similarly define  
$\bar p_{ij}:\bar Y^{\otimes 2}\to\bar Y^{\otimes 3}$.
Suppose there is given a descent datum $(\bar M,\bar\beta)$ 
for $\bar M$, relative to $\bar X\to\bar Y$. The 
cocycle condition \eqref{eq_cocycle} implies easily that 
$\bar\mu_2^*(\bar\beta)$ is the identity on 
$\bar\mu_2^*(\bar p_i^*(\bar M))=\bar M$. It follows
that the pair $(\bar\beta,\one_M)$ defines an isomorphism
$\pi_2(p_1^*M)\stackrel{\sim}{\to}\pi_2(p_2^*M)$ in 
the category
$\bar Y^{\otimes 2}\Mod\times_{\bar Y\Mod}Y\Mod$.
Hence 
$T_2(\bar\beta,\one_M):T_2\circ\pi_2(p_1^*M)\to
T_2\circ\pi_2(p_2^*M)$ is an isomorphism. However, 
we remark that either morphism $\bar p_i$ yields a section 
for $\mu_2$, hence we are in the situation 
contemplated in lemma \ref{lem_special.case}, and we derive
an isomorphism $\beta:p_2^*(M)\stackrel{\sim}{\to}p_1^*(M)$.
We claim that $(M,\beta)$ is an object of $\Desc(\cC,Y/X)$, 
{\em i.e.\/} that $\beta$ verifies the cocycle condition 
\eqref{eq_cocycle}. Indeed, we can compute: 
$\pi_3(p_{ij}^*\beta)=
(\rho_3^*(p_{ij}^*\beta),\mu_3^*(p_{ij}^*\beta))$ and by
construction we have 
$\rho_3^*(p_{ij}^*\beta)=\bar p_{ij}^*(\bar\beta)$ and
$\mu_3^*(p_{ij}^*\beta)=\mu_2^*(\beta)=\one_M$. Therefore,
the cocycle identity for $\bar\beta$ implies the equality
$\pi_3(p_{12}^*(\beta))\circ\pi_3(p_{23}^*(\beta))=
\pi_3(p_{13}^*(\beta))$. If we now apply the functor $T_3$
to this equality, and then invoke again lemma 
\ref{lem_special.case}, the required cocycle identity
for $\beta$ will ensue. Clearly $\beta$ is the only 
descent datum on $M$ lifting $\bar\beta$. This proves
that \eqref{eq_Desc} is essentially surjective. The same 
sort of argument also shows that the functor \eqref{eq_Desc} 
induces bijections on morphisms, so the lemma follows
in this case. Next, the case $\cC=V^a\AlgMorph^o$ can
be deduced formally from the previous case, by applying
repeatedly natural isomorphisms of the kind
$p_i^*(M\otimes_YN)\simeq 
p_i^*(M)\otimes_{Y\otimes_XY}p_i^*(N)$ ($i=1,2$).
Finally, the ``\'etaleness'' of an object of 
$\Desc(V^a\AlgMorph^o,Y/X)$ can be checked on its
projection onto $Y\Alg^o$, hence also the cases
$\cC=\wEt^o$ and $\cC=\Et^o$ follow directly.
\end{proof}

Now, let $B=A_1\times A_2$; to an objet $(M,\beta)$ in 
$\Desc(V^a\AlgMod^o,B/A)$ we assign an object $(M_1,M_2,\xi)$ 
of $\cM_1\times_{\cM_3}\cM_2$, as follows. Set $M_i=A_i\otimes_BM$ 
($i=1,2$) and $A_{ij}=A_i\otimes_{A_0}A_j$. We can write 
$B\otimes_{A_0}B=\prod_{i,j=1}^2A_{ij}$ and $\beta$
gives rise to the $A_{ij}$-linear isomorphisms
$\beta_{ij}:A_{ij}\otimes_{B\otimes_{A_0}B}p_2^*(M)
\stackrel{\sim}{\to}A_{ij}\otimes_{B\otimes_{A_0}B}p_1^*(M)$.
In other words, we obtain isomorphisms 
$\beta_{ij}:A_i\otimes_{A_0}M_j\to M_i\otimes_{A_0}A_j$.
However, we have a natural isomorphism $A_{12}\simeq A_3$
(indeed, suppose that $A_1\to A_3$ is an epimorphism with 
kernel $I$; then $I$ is also an ideal of $A_0$ and 
$A_0/I\simeq A_2$; now the claim follows by 
remarking that $I\cdot A_1=I$). Hence we can choose 
$\xi=\beta_{12}$. In this way we obtain a functor :
\set\begin{equation}\label{eq_funct.desc}
\Desc(V^a\AlgMod^o,B/A_0)\to(\cM_1\times_{\cM_3}\cM_2)^o.
\end{equation}
\begin{proposition}\label{prop_equiv.desc}
The functor \eqref{eq_funct.desc} is an equivalence of 
categories.
\end{proposition}
\begin{proof} Let us say that $A_1\to A_3$ is an epimorphism
with kernel $I$. Then $I$ is also an ideal of $B$ and we have
$B/I\simeq A_3\times A_2$ and $A_0/I\simeq A_2$. We intend to
apply lemma \ref{lem_Desc} to the morphism $A_0\to B$.
However, the induced morphism 
$\bar B=B/I\to\bar A_0=A_0/I$ in $V^a\Alg^o$ has a section, 
and hence it is of universal effective descent for every 
fibred category. Thus, we can replace in \eqref{eq_Desc} 
the category 
$\Desc(V^a\AlgMod^o,\bar B/\bar A_0)$ by $\bar A_0\Mod^o$,
and thereby, identify (up to equivalence) the target of 
\eqref{eq_Desc} with the 2-fibred product
$(\cM_1\times\cM_2)^o\times_{(\cM_3\times\cM_2)^o}\cM_2^o$.
The latter is equivalent to the category 
$\cM_1^o\times_{\cM_3^o}\cM_2^o$ and the resulting
functor 
$\Desc(V^a\AlgMod^o,B/A_0)\to\cM_1^o\times_{\cM_3^o}\cM_2^o$
is canonically isomorphic to \eqref{eq_funct.desc}, which 
gives the claim.
\end{proof}

Putting together propositions \ref{prop_equiv.2-prod} and
\ref{prop_equiv.desc} we obtain the following :
\begin{corollary}\label{cor_eff-desc} In the situation of 
\eqref{eq_cartesian}, the morphism $A_0\to A_1\times A_2$ is 
of effective descent for the fibred categories of flat almost 
modules and of almost projective almost modules.
\qed\end{corollary}

Next we would like to give sufficient conditions to ensure
that a morphism of almost algebras is of effective descent for
the fibred category $\wEt^o\to V^a\Alg^o$ of weakly \'etale
algebras (resp. for \'etale algebras). To this aim we are led 
to the following :
\begin{definition} A morphism $\phi:A\to B$ of almost algebras is
said to be {\em strictly finite\/} if $\Ker(\phi)$ is nilpotent
and $B\simeq R^a$, where $R$ is a finite $A_*$-algebra.
\end{definition}

\begin{theorem}\label{th_strictly} Let $\phi:A\to B$ be a strictly 
finite morphism of almost algebras. Then :

i) For every $A$-algebra $C$, the induced morphism 
$C\to C\otimes_AB$ is again strictly finite.

ii) If $M$ is a flat $A$-module and $B\otimes_AM$ is
almost projective over $B$, then $M$ is almost projective 
over $A$.

iii) $A\to B$ is of universal effective descent for the fibred 
categories of weakly \'etale (resp. \'etale) almost algebras.
\end{theorem}
\begin{proof} (i): suppose that $B=R^a$ for a finite $A_*$-algebra
$R$; then $S=C_*\otimes_{A_*}R$ is a finite $C_*$-algebra and
$S^a\simeq C$. It remains to show that $\Ker(C\to C\otimes_AB)$
is nilpotent. Suppose that $R$ is generated by $n$ elements as 
an $A_*$-module and let $F_{A_*}(R)$ (resp. $F_{C_*}(S)$) be 
the Fitting ideal of $R$ (resp.of $S$); we have 
$\Ann_{C_*}(S)^n\subset F_{C_*}(S)\subset\Ann_{C_*}(S)$ (see
\cite{Lan} (Chap.XIX Prop.2.5)); on the other hand
$F_{C_*}(S)=F_{A_*}(R)\cdot C_*$, so the claim is clear.

(iii): we shall consider the fibred category  
$F:\wEt^o\to V^a\Alg^o$; the same argument applies also to 
\'etale almost algebras. We begin by establishing a very 
special case :
\begin{claim}\label{cl_special} Assertion (iii) holds when 
$B=(A/I_1)\times(A/I_2)$, where $I_1$ and $I_2$ are ideals 
in $A$ such that $I_1\cap I_2$ is nilpotent.
\end{claim}
\begin{pfclaim} First of all we remark that the situation 
considered in the claim is stable under arbitrary base change, 
therefore it suffices to show that $\phi$ is of $F$-2-descent 
in this case. Then we factor $\phi$ as a composition
$A\to A/\Ker(\phi)\to B$ and we remark that $A\to A/\Ker(\phi)$
is of $F$-2-descent by theorem \ref{th_liftetale}; since a composition
of morphisms of $F$-2-descent is again of $F$-2-descent, we are reduced
to show that $A/\Ker(\phi)\to B$ is of $F$-2-descent, {\em i.e.\/} 
we can assume that $\Ker(\phi)\simeq 0$. However, in this case the
claim follows easily from corollary \ref{cor_eff-desc}.
\end{pfclaim}

\begin{claim}\label{cl_more-ideals} More generally, assertion (iii) 
holds when $B=\prod_{i=1}^nA/I_i$, where $I_1,...,I_n$ 
are ideals of $A$, such that $\bigcap^n_{i=1}I_i$ is nilpotent.
\end{claim}
\begin{pfclaim} We prove this by induction on $n$, the case $n=2$ 
being covered by claim \ref{cl_special}. Therefore, suppose that 
$n>2$, and set $B'=A/(\bigcap_{i=1}^{n-1}I_j)$. By induction,
the morphism $B'\to\prod_{i=1}^{n-1}A/I_i$ is of universal 
$F$-2-descent. However, according to \cite{Gi} (Chap.II Prop.1.1.3), 
the sieves of universal $F$-2-descent form a topology on $V^a\Alg^o$;
for this topology, $\{A,B\}$ is a covering family of $A\times B$ and
$(A\to B'\times(A/I_n))^o$ is a covering morphism, hence $\{B',A/I_n\}$
is a covering family of $A$, and then, by composition of covering 
families, $\{\prod_{i=1}^{n-1}A/I_i,A/I_n\}$ is a covering family
of $A$, which is equivalent to the claim.
\end{pfclaim}

Now, let $A\to B$ be a general strictly finite morphism, so
that $B=R^a$ for some finite $A_*$-algebra $R$. Pick generators 
$f_1,...,f_m$ of the $A_*$-module $R$, and monic
polynomials $p_1(X),...,p_m(X)$ such that $p_i(f_i)=0$ 
for $i=1,...,m$. 
\begin{claim}\label{cl_splitter} There exists a finite and faithfully 
flat extension $C$ of $A_*$ such that the images  in $C[X]$ of 
$p_1(X)$,...,$p_m(X)$ split as products of monic linear factors. 
\end{claim}
\begin{pfclaim} This extension $C$ can be obtained as follows. 
It suffices to find, for each $i=1,...,m$, an extension $C_i$ 
that splits $p_i(X)$, because then 
$C=C_1\otimes_{A_*}...\otimes_{A_*}C_m$ will split 
them all, so we can assume that $m=1$ and $p_1(X)=p(X)$; moreover, 
by induction on the degree of $p(X)$, it suffices to find an extension 
$C'$ such that $p(X)$ factors in $C'[X]$ as a product of the form
$p(X)=(X-\alpha)\cdot q(X)$, where $q(X)$ is a monic polynomial
of degree $\deg(p)-1$. Clearly we can take $C'=A_*[T]/(p(T))$.
\end{pfclaim}

Given a $C$ as in claim \ref{cl_splitter}, we remark that the morphism 
$A\to C^a$ is of universal $F$-2-descent. Considering again the 
topology of universal $F$-2-descent, it follows that $A\to B$ 
is of universal $F$-2-descent if and only if the same holds for the 
induced morphism $C^a\to C^a\otimes_AB$. Therefore, in proving 
assertion (iii) we can replace $\phi$ by $\one_C\otimes_A\phi$ 
and assume from start that the polynomials $p_i(X)$ factor in $A_*[X]$ 
as product of linear factors. Now, let $\deg(p_i)=d_i$
and $p_i(X)=\prod^{d_i}_j(X-\alpha_{ij})$ (for $i=1,...,m$). 
We get a surjective homomorphism of $A_*$-algebras 
$D=A_*[X_1,...,X_m]/(p_1(X_1),...,p_m(X_m))\to R$ by
the rule $X_i\mapsto f_i$ ($i=1,...,m$). Moreover, any 
sequence $\ubarold\alpha=
(\alpha_{1,j_1},\alpha_{2,j_2},...,\alpha_{m,j_m})$
yields a homomorphism $\psi_{\ubarold\alpha}:D\to A_*$, 
determined by the assignment $X_i\mapsto\alpha_{i,j_i}$.
A simple combinatorial argument shows that 
$\prod_{\ubarold\alpha}\Ker(\psi_{\ubarold\alpha})=0$, where 
$\ubarold\alpha$ runs over all the sequences as above.
Hence the product map
$\prod_{\ubarold\alpha}\psi_{\ubarold\alpha}:
D\to\prod_{\ubarold\alpha}A_*$
has nilpotent kernel. We notice that the $A_*$-algebra 
$(\prod_{\ubarold\alpha}A_*)\otimes_DR$ is a quotient of 
$\prod_{\ubarold\alpha}A_*$, hence it can be written as a 
product of rings of the form $A_*/I_{\ubarold\alpha}$, 
for various ideals $I_{\ubarold\alpha}$. By (i), the
kernel of the induced homomorphism 
$R\to\prod_{\ubarold\alpha}A_*/I_{\ubarold\alpha}$ is nilpotent,
hence the same holds for the kernel of the 
composition $A\to\prod_{\ubarold\alpha}A/I^a_{\ubarold\alpha}$, 
which is therefore of the kind considered in claim 
\ref{cl_more-ideals}. Hence  
$A\to\prod_{\ubarold\alpha}A/I^a_{\ubarold\alpha}$
is of universal $F$-2-descent. Since such 
morphisms form a topology, it follows that 
also $A\to B$ is of universal $F$-2-descent, which 
concludes the proof of (iii).

Finally, let $M$ be as in (ii) and pick again $C$ as in 
the proof of claim \ref{cl_splitter}. By remark 
\ref{rem_descent}(iv), $M$ is almost 
projective over $A$ if and only if $C^a\otimes_AM$ 
is almost projective over $C^a$; hence we can replace 
$\phi$ by $\one_{C^a}\otimes_A\phi$, and by arguing 
as in the proof of (iii), we can assume from start 
that $B=\prod_{j=1}^n(A/I_j)$ for ideals
$I_j\subset A$, $j=1,...,n$ such that 
$I=\bigcap^n_{j=1}I_j$ is nilpotent. By an easy induction, 
we can furthermore reduce to the case $n=2$. We factor
$\phi$ as $A\to A/I\to B$; by proposition 
\ref{prop_equiv.2-prod} it follows that $(A/I)\otimes_AM$
is almost projective over $A/I$, and then lemma 
\ref{lem_flat.proj}(i) says that $M$ itself is almost 
projective.
\end{proof}

\begin{remark} It is natural to ask whether theorem 
\ref{th_strictly} holds if we replace everywhere ``strictly 
finite'' by ``finite with nilpotent kernel'' (or even by 
``almost finite with nilpotent kernel''). 
We do not know the answer to this question.
\end{remark}

We conclude with a digression to explain the relationship
between our results and related facts that can be extracted
from the literature. So, we now place ourselves in the 
``classical limit'' $V=\fm$ (cp. example \ref{ex_rings}(ii)). 
In this case, weakly \'etale morphisms had already been 
considered in some earlier work, and they were called 
``absolutely flat'' morphisms. A ring homomorphism $A\to B$
is \'etale in the usual sense of \cite{SGA1} if and only if 
it is absolutely flat and of finite presentation. Let us denote
by $\uEt^o$ the fibred category over $V\Alg^o$, whose fibre
over a $V$-algebra $A$ is the opposite of the category of 
\'etale $A$-algebras in the usual sense. We claim that, if 
a morphism $A\to B$ of $V$-algebras is of universal effective 
descent for the fibred category $\wEt^o$ (resp. $\Et^o$), 
then it is a morphism of universal effective descent for 
$\uEt^o$. Indeed, let $C$ be an \'etale $A$-algebra (in the 
sense of definition \ref{def_morph}) and such that $C\otimes_AB$ 
is \'etale over $B$ in the usual sense. We have to show that $C$
is \'etale in the usual sense, {\em i.e.\/} that it is of
finite presentation over $A$. This amounts to showing that,
for every filtered inductive system 
$(A_\lambda)_{\lambda\in\Lambda}$ of $A$-algebras, we have
$\colim{\lambda\in\Lambda}\Hom_{A\Alg}(C,A_\lambda)\simeq
\Hom_{A\Alg}(C,\colim{\lambda\in\Lambda}A_\lambda)$. Since, 
by assumption, this is known after extending scalars to $B$
and to $B\otimes_AB$, it suffices to show that, for any 
$A$-algebra $D$, the natural sequence
$$\xymatrix{\Hom_{A\Alg}(C,D)\ar[r] & \Hom_{B\Alg}(C_B,D_B)
\ar@<.5ex>[r] \ar@<-.5ex>[r] &
\Hom_{B\otimes_AB\Alg}(C_{B\otimes_AB},D_{B\otimes_AB})
}$$
is exact. For this, note that 
$\Hom_{A\Alg}(C,D)=\Hom_{D\Alg}(C_D,D)$ (and similarly for 
the other terms) and by hypothesis $(D\to D\otimes_AB)^o$ is
a morphism of 1-descent for the fibred category $\wEt^o$ 
(resp. $\Et^o$). 

As a consequence of these observations and 
of theorem \ref{th_strictly}, we see that any finite ring 
homomorphism $\phi:A\to B$ with nilpotent kernel is of universal 
effective descent for the fibred category of \'etale algebras. 
This fact was known as follows. By \cite{SGA1} (Exp.IX, 4.7), 
$\Spec(\phi)$ is of universal effective descent for the fibred 
category of separated \'etale morphisms of finite type. One 
has to show that if $X$ is such a scheme over $A$, such that 
$X\otimes_AB$ is affine, then $X$ is affine. This follows by 
reduction to the noetherian case and \cite{EGA} (Chap.II, 6.7.1).

\subsection{Behaviour of \'etale morphisms under Frobenius}
We consider the following category $\cB$ of base rings. The 
objects of $\cB$ are the pairs $(V,\fm)$, where $V$ is a ring
and $\fm$ is an ideal of $V$ with $\fm=\fm^2$ and $\tilde\fm$ 
is flat. The morphisms $(V,\fm_V)\to(W,\fm_W)$ between two 
objects of $\cB$ are the ring homomorphisms $f:V\to W$ such 
that $\fm_W=f(\fm_V)\cdot W$.

We have a fibred and cofibred category $\cB\Mod\to\cB$ (see 
\cite{SGA1} (Exp.VI \S5,6,10) for generalities on fibred
categories). An object of $\cB\Mod$ (which we may call a 
``$\cB$-module'') consists of a pair $((V,\fm),M)$, where $(V,\fm)$
is an object of $\cB$ and $M$ is a $V$-module. Given two objects
$X=((V,\fm_V),M)$ and $Y=((W,\fm_W),N)$, the morphisms $X\to Y$
are the pairs $(f,g)$, where $f:(V,\fm_V)\to(W,\fm_W)$ is a 
morphism in $\cB$ and $g:M\to N$ is an $f$-linear map.

Similarly one has a fibred and cofibred category $\cB\Alg\to\cB$
of $\cB$-algebras. We will also need to consider the fibred and
cofibred category $\cB\Mon\to\cB$ of non-unitary commutative
$\cB$-monoids: an object of $\cB\Mon$ is a pair $((V,\fm),A)$ where 
$A$ is a $V$-module endowed with a morphism $A\otimes_VA\to A$
subject to associativity and commutativity conditions, as 
discussed in section \ref{sec_alm.cat}. The fibre over
an object $(V,\fm)$ of $\cB$, is the category of $V$-monoids
denoted $(V,\fm)\Mon$ or simply $V\Mon$. 

The almost isomorphisms in the fibres of $\cB\Mod\to\cB$ give a 
multiplicative system $\Sigma$ in $\cB\Mod$, admitting a calculus
of both left and right fractions. The ``locally small''
conditions are satisfied (see \cite{We} p.381), so that one 
can form the localised category $\cB^a\Mod=\Sigma^{-1}(\cB\Mod)$. 
The fibres of the localised category over the objects of 
$\cB$ are the previously considered categories of almost 
modules. Similar considerations hold for $\cB\Alg$ and $\cB\Mon$, 
and we get the fibred and cofibred categories $\cB^a\Mod\to\cB$, 
$\cB^a\Alg\to\cB$ and $\cB^a\Mon\to\cB$. In particular, for 
every object $(V,\fm)$ of $\cB$, we have an obvious notion of
almost $V$-monoid and the category consisting of these is 
denoted $V^a\Mon$. The localisation functors 
\smallskip

\centerline{$\cB\Mod\to\cB^a\Mod\ :\ M\mapsto M^a \qquad 
   \cB\Alg\to\cB^a\Alg\ :\ B\mapsto B^a$}
\smallskip

have left and right adjoints. These adjoints can be chosen
as functors of categories over $\cB$ such that the adjunction
units and counits are morphisms over identity arrows in $\cB$.
On the fibres these induce the previously considered left and 
right adjoints $M\mapsto M_!$, $M\mapsto M_*$, $B\mapsto B_{!!}$, 
$B\mapsto B_*$.
We will use the same notation for the corresponding functors 
on the larger categories. Then it is easy to check that the 
functor $M\mapsto M_!$ is cartesian and cocartesian ({\em i.e.\/}
it sends cartesian arrows to cartesian arrows and cocartesian
arrows to cocartesian arrows), $M\mapsto M_*$ and $B\mapsto B_*$
are cartesian, and $B\mapsto B_{!!}$ is cocartesian.

Let $\cB/\F_p$ be the full subcategory of $\cB$ consisting of
all objects $(V,\fm)$ where $V$ is an $\F_p$-algebra. Define 
similarly $\cB\Alg/\F_p$, $\cB\Mon/\F_p$ and $\cB^a\Alg/\F_p$, 
$\cB^a\Mon/\F_p$, so that we have again fibred and cofibred categories 
$\cB^a\Alg/\F_p\to\cB/\F_p$ and $\cB^a\Alg/\F_p\to\cB/\F_p$ (resp.
the same for non-unitary monoids). We remark that the categories
$\cB^a\Alg/\F_p$ and $\cB^a\Mon/\F_p$ have small limits and colimits, 
and these are preserved by the projection to $\cB/\F_p$. Especially,
if $A\to B$ and $A\to C$ are two morphisms in $\cB^a\Alg/\F_p$
or $\cB^a\Mon/\F_p$, we can define $B\otimes_AC$ as such a colimit.

If $A$ is a (unitary or non-unitary) $\cB$-monoid over $\F_p$, 
we denote by $\phi_A:A\to A$ the Frobenius endomorphism 
$x\mapsto x^p$. If $(V,\fm)$ is an object of $\cB/\F_p$, it 
follows from proposition \ref{prop_less.obv}(ii) that 
$\phi_V:(V,\fm)\to(V,\fm)$ is a morphism in $\cB$. 
If $B$ is an object of $\cB\Alg/\F_p$ (resp. $\cB\Mon/\F_p$) 
over $V$, then the Frobenius map induces a morphism 
$\phi_B:B\to B$ in $\cB\Alg/\F_p$ (resp. $\cB\Mon/\F_p$) over 
$\phi_V$. In this way we get a natural transformation from the 
identity functor of $\cB\Alg/\F_p$ (resp. $\cB\Mon/\F_p$) to 
itself that induces a natural transformation on the identity 
functor of $\cB^a\Alg/\F_p$ (resp. $\cB^a\Mon/F_p$).

Using the pull-back functors, any object $B$ of $\cB\Alg$ over 
$V$ defines new objects $B_{(m)}$ of $\cB\Alg$ ($m\in\N$) over 
$V$, where $B_{(m)}=(\phi^m_V)^*(B)$, which is just $B$
considered as a $V$-algebra via the homomorphism 
$V\stackrel{\phi^m}{\longrightarrow} V\to B$. These 
operations also induce functors $B\mapsto B_{(m)}$ on 
almost $\cB$-algebras.

\begin{definition}\label{def_invert.phi} 
i) Let $(V,\fm)$ be an object of $\cB/\F_p$; 
we say that a morphism $f:A\to B$ of almost $V$-algebras 
(resp. almost $V$-monoids) is {\em invertible up to\/} $\phi^m$ 
if there exists a morphism $f':B\to A$ in $\cB^a\Alg$ (resp. 
$\cB^a\Mon$) over $\phi^m_V$, such that $f'\circ f=\phi^m_A$ 
and $f\circ f'=\phi^m_B$.

ii) We say that an almost $V$-monoid $I$ ({\em e.g.\/} an
ideal in a $V^a$-algebra) is {\em Frobenius nilpotent\/}
if $\phi_I$ is nilpotent.
\end{definition}
Notice that a morphism $f$ of $V^a\Alg$ (or $V^a\Mon$) is 
invertible up to $\phi^m$ if and only if $f_*:A_*\to B_*$ 
is so as a morphism of $\F_p$-algebras.

\begin{lemma}\label{lem_invert} Let $(V,\fm)$ be an 
object of $\cB/\F_p$ and let $f:A\to B$, $g:B\to C$ 
be  morphisms of almost $V$-algebras or almost $V$-monoids. 

i) If $f$ is invertible up to $\phi^n$ and $g$ is
invertible up to $\phi^m$, then $g\circ f$ is invertible 
up to $\phi^{m+n}$.

ii) If $f$ is invertible up to $\phi^n$ and $g\circ f$ 
is invertible up to $\phi^m$, then $g$ is invertible up 
to $\phi^{m+n}$.

iii) If $g$ is invertible up to $\phi^n$ and $g\circ f$ 
is invertible up to $\phi^m$, then $f$ is invertible up
to $\phi^{m+n}$.

iv) The Frobenius morphisms induce $\phi_V$-linear morphisms
({\em i.e.\/} morphisms in $\cB^a\Mod$ over $\phi_V$)
$\phi':\Ker(f)\to\Ker(f)$ and 
$\phi'':\Coker(f)\to\Coker(f)$, and $f$ is
invertible up to some power of $\phi$ if and only
if both $\phi'$ and $\phi''$ are nilpotent.

v) Consider a map of short exact sequences of almost 
$V$-monoids :
$$\xymatrix{0 \ar[r] & A' \ar[r] \ar[d]_{f'} & 
                       A \ar[r] \ar[d]_f &
                       A'' \ar[r] \ar[d]_{f''} & 0 \\
            0 \ar[r] & B' \ar[r] & B \ar[r] & B'' \ar[r] & 0
}$$
and suppose that two of the morphisms $f',f,f''$ are
invertible up to a power of $\phi$. Then also the third
morphism has this property.
\end{lemma}
\begin{proof} (i): if $f'$ is an inverse of $f$ up to
$\phi^n$ and $g'$ is an inverse of $g$ up to $\phi^m$,
then $f'\circ g'$ is an inverse of $g\circ f$ up to
$\phi^{m+n}$. (ii): given an inverse $f'$ of $f$ up 
to $\phi^n$ and an inverse $h'$ of $h=g\circ f$ up to 
$\phi^m$, let $g'=\phi^n_B\circ f\circ h'$. We compute :
$$\begin{array}{l@{\: =\:}l}
g\circ g' & g\circ\phi^n_B\circ f\circ h'=
\phi^n_C\circ g\circ f\circ h=\phi^n_C\circ\phi^m_C \\
g'\circ g & \phi^n_B\circ f\circ h'\circ g=
f\circ h'\circ g\circ\phi^n_B=
f\circ h'\circ g\circ f\circ f'\\
& f\circ\phi^m_A\circ f'=\phi^m_B\circ f\circ f'=
\phi^m_B\circ\phi^n_B.
\end{array}$$
(iii) is similar and (iv) is an easy diagram chasing
left to the reader. (v) follows from (iv) and the snake 
lemma.
\end{proof}

\begin{lemma}\label{lem_pull.and.push}
Let $(V,\fm)$ be an object of $\cB/\F_p$.

(i) If $f:A\to B$ is a morphism of almost $V$-algebras, 
invertible up to $\phi^n$, then so is $A'\to A'\otimes_AB$ 
for every morphism $A\to A'$ of almost algebras.

(ii) If $f:(V,\fm_V)\to(W,\fm_W)$ is a morphism in 
$\cB/\F_p$, the functors 
$f_*:(V,\fm_V)^a\Alg\to(W,\fm_W)^a\Alg$ 
and $f^*:(W,\fm_W)^a\Alg\to(V,\fm_V)^a\Alg$ preserve 
the class of morphisms invertible up to $\phi^n$.
\end{lemma}
\begin{proof} (i): given $f':B\to A_{(m)}$, construct a 
morphism $A'\otimes_AB\to A'_{(m)}$ using the morphism 
$A'\to A'_{(m)}$ coming from $\phi^m_{A'}$ and $f'$.
(ii): the assertion for $f^*$ is clear, and the assertion
for $f_*$ follows from (i).
\end{proof}

\begin{remark}\label{rem_pull.push} Statements like those 
of lemma \ref{lem_pull.and.push} hold for the classes of 
flat, (weakly) unramified, (weakly) \'etale morphisms.
\end{remark}

\begin{theorem}\label{th_Frobenius} 
Let $(V,\fm)$ be an object of $\cB/\F_p$ and 
$f:A\to B$ a weakly \'etale morphism of almost 
$V$-algebras.

(i) If $f$ is invertible up to $\phi^n$ ($n\geq 0$), 
then it is an isomorphism.

(ii) For every integer $m\geq 0$ the natural square 
diagram
\set\begin{equation}\label{eq_Frobenius}{
\diagram A \ar[r]^f \ar[d]_{\phi^m_A} & 
         B \ar[d]^{\phi^m_B} \\
         A_{(m)} \ar[r]^{f_{(m)}} & B_{(m)}
\enddiagram
}\end{equation}
is cocartesian.
\end{theorem}
\begin{proof} (i): we first show that $f$ is faithfully 
flat. Since $f$ is flat, it remains to show that if $M$ 
is an $A$-module such that $M\otimes_AB=0$, then 
$M=0$. It suffice to do this for $M=A/I$, for an arbitrary 
ideal $I$ of $A$. After base change by $A\to A/I$, we 
reduce to show that $B=0$ implies $A=0$. However, 
$A_*\to B_*$ is invertible up to $\phi^n$, so 
$\phi^n_{A_*}=0$ which means $A_*=0$. In particular, 
$f$ is a monomorphism, hence the proof is complete in 
case that $f$ is an epimorphism. In general, consider 
the composition
$B\stackrel{\one_B\otimes f}{\longrightarrow}
B\otimes_AB\stackrel{\mu_{B/A}}{\longrightarrow}B$.
From lemma \ref{lem_pull.and.push}(i) it follows that 
$\one_B\otimes f$ is invertible up to $\phi^n$; then 
lemma \ref{lem_invert}(ii) says that $\mu_{B/A}$ is 
invertible up to $\phi^n$. The latter is also weakly
\'etale; by the foregoing we derive that it is an
isomorphism. Consequently $\one_B\otimes f$ is an
isomorphism, and finally, by faithful flatness, $f$
itself is an isomorphism.

(ii): the morphisms $\phi^m_A$ and $\phi^m_B$ are 
invertible up to $\phi^m$. By lemma 
\ref{lem_pull.and.push}(i) it follows that 
$\one_B\otimes\phi^m_A:B\to B\otimes_AA_{(m)}$ is
invertible up to $\phi^m$; hence, by lemma 
\ref{lem_invert}(ii), the morphism 
$h:B\otimes_AA_{(m)}\to B_{(m)}$ induced by $\phi^m_B$
and $f_{(m)}$ is invertible up to $\phi^{2m}$ (in fact 
one verifies that it is invertible up to $\phi^m$). 
But $h$ is a morphism of weakly \'etale $A_{(m)}$-algebras, 
so it is weakly \'etale, so it is an isomorphism by (i).
\end{proof}

\begin{remark} Theorem \ref{th_Frobenius}(ii) extends
a statement of Faltings (\cite{Fa2} p.10) for his notion
of almost \'etale extensions. 
\end{remark} 

We recall (cp. \cite{Gi} (Chap.0, 3.5)) that a morphism 
$f:X\to Y$ of objects in a site is called {\em bicovering\/} 
if the induced map of associated sheaves of sets is an 
isomorphism; if $f$ is squarable (``quarrable'' in French), 
this is equivalent to the condition that both $f$ and the 
diagonal morphism $X\to X\times_YX$ are covering morphisms.

Let $F\to E$ be a fibered category and $f:P\to Q$ a
squarable morphism of $E$. Consider the following
condition:
\set\begin{equation}\label{eq_bicover}{
\parbox{14cm}
{for every base change $P\times_QQ'\to Q'$ of $f$, the 
inverse image functor $F_{Q'}\to F_{P\times_QQ'}$ is 
an equivalence of categories.}
}\end{equation}
Inspecting the arguments in \cite{Gi} (Chap.II,\S 1.1) 
one can show:
\begin{lemma}\label{lem_bicover} With the above notation, 
let $\tau$ be the topology of  universal effective descent 
relative to $F\to E$. Then we have :

i) if \eqref{eq_bicover} holds, then $f$ is a covering morphism 
for the topology $\tau$;

ii) $f$ is bicovering for $\tau$ if and only if \eqref{eq_bicover} 
holds both for $f$ and for the diagonal morphism $P\to P\times_QP$. 
\end{lemma}
\begin{remark} In \cite{Gi} (Chap.II, 1.1.3(iv)) it
is stated that ``la r\'eciproque est vraie si $i=2$'', 
meaning that \eqref{eq_bicover} is equivalent to the condition 
that $f$ is bicovering for $\tau$. (Actually the cited statement
is given in terms of presheaves, but one can show that
\eqref{eq_bicover} is equivalent to the corresponding condition 
for the fibered category $F^+\to\hat E_U$ considered in 
{\em op.cit.}) However, this fails in general : as a counterexample 
we can give the following. Let $E$ be the category of
schemes of finite type over a field $k$; set $P=\A^1_k$,
$Q=\Spec(k)$. Finally let $F\to E$ be the discretely
fibered category defined by the presheaf 
$X\mapsto H^0(X,\Z)$. Then it is easy to show that $f$
satisfies \eqref{eq_bicover} but the diagonal map does not, 
so $f$ is not bicovering. The mistake in the proof is in 
\cite{Gi} (Chap.II, 1.1.3.5), where one knows that 
$F^+(d)$ is an equivalence of categories (notation 
of {\em loc.cit.}) but one needs it also after base 
changes of $d$. 
\end{remark}

\begin{lemma}\label{lem_invert.Et} (i)
If $f:A\to B$ is a morphism of $V^a$-algebras which is 
invertible up to $\phi^m$, then the induced functors 
$\Et(A)\to\Et(B)$ and $\wEt(A)\to\wEt(B)$ are equivalences 
of categories.

ii) If $A\to B$ is weakly \'etale and $C\to D$ is a morphism 
of $A$-algebras invertible up to $\phi^m$, then the induced
map $\Hom_{A\Alg}(B,C)\to\Hom_{A\Alg}(B,D)$ is bijective.
\end{lemma} 
\begin{proof} We first consider (i) for the special case 
where $f=\phi^m_A:A\to A_{(m)}$. The functor 
$(\phi^m_V)^*:V^a\Alg\to V^a\Alg$ induces a functor
$(-)_{(m)}:A\Alg\to A_{(m)}\Alg$, and by restriction 
(see remark \ref{lem_pull.and.push}) we obtain a functor 
$(-)_{(m)}:\Et(A)\to\Et(A_{(m)})$; by theorem 
\ref{th_Frobenius}(ii), the latter is isomorphic to
the functor $(\phi^m)_*:\Et(A)\to\Et(A_{(m)})$ of the 
lemma. Furthermore, from remark \ref{rem_almost.zero}(ii) and
\eqref{eq_alm.morph} we derive a natural ring isomorphism 
$\omega:A_{(m)*}\simeq A_*$, hence an essentially
commutative diagram
$$\xymatrix{
\Et(A) \ar[r] \ar[d]_{(\phi^m)_*} & 
A\Alg \ar[r]^-\alpha \ar[d]_{(-)_{(m)}} & 
(A_*,\fm\cdot A_*)^a\Alg \ar[d]_{\omega^*} \\
\Et(A_{(m)}) \ar[r] & A_{(m)}\Alg \ar[r]^-\beta &
(A_{(m)*},\fm\cdot A_{(m)*})^a\Alg 
}$$
where $\alpha$ and $\beta$ are the equivalences 
of remark \ref{rem_base.comp}. Clearly $\alpha$
and $\beta$ restrict to equivalences on the 
corresponding categories of \'etale algebras, 
hence the lemma follows in this case. 

For the general case of (i), let $f':B\to A_{(m)}$ be
a morphism as in definition \ref{def_invert.phi}.
Diagram \eqref{eq_Frobenius} induces an essentially
commutative diagram of the corresponding categories
of algebras, so by the previous case, the functor
$(f')_*:\Et(B)\to\Et(A_{(m)})$ has both a left
essential inverse and a right essential inverse;
these essential inverses must be isomorphic, so
$f_*$ has an essential inverse as desired.
Finally, we remark that the map in (ii) is the same
as the map 
$\Hom_{C\Alg}(B\otimes_AC,C)\to\Hom_{D\Alg}(B\otimes_AD,D)$,
and the latter is a bijection in view of (i).
\end{proof}
\begin{remark} Notice that lemma \ref{lem_invert.Et}(ii)
generalises the lifting theorem \ref{th_liftetale}(i) (in 
case $V$ is an $\F_p$-algebra). Similarly, it follows 
from lemmata \ref{lem_invert.Et}(i) and \ref{lem_invert}(iv) 
that, in case $V$ is an $\F_p$-algebra, one can replace 
``nilpotent'' in theorem \ref{th_liftetale} parts (ii) and (iii) 
by ``Frobenius nilpotent''. 
\end{remark}
In the following $\tau$ will denote indifferently the topology 
of universal effective descent defined by either of the fibered 
categories $\wEt^o\to V^a\Alg^o$ or $\Et^o\to V^a\Alg^o$.
\begin{proposition}\label{prop_invert.bicov}
If $f:A\to B$ is a morphism of almost 
$V$-algebras which is invertible up to $\phi^m$, then $f^o$\/
is bicovering for the topology $\tau$.
\end{proposition}
\begin{proof} In light of lemmata \ref{lem_bicover}(ii) and 
\ref{lem_invert.Et}(i), it suffices to show that $\mu_{B/A}$ 
is invertible up to a power of $\phi$. For this, factor the 
identity morphism of $B$ as 
$B\stackrel{\one_B\otimes f}{\longrightarrow}
B\otimes_AB\stackrel{\mu_{B/A}}{\longrightarrow}B$ and argue 
as in the proof of theorem \ref{th_Frobenius}.
\end{proof}
\begin{proposition} Let $A\to B$ be a morphism of almost 
$V$-algebras and $I\subset A$ an ideal. Set $\bar A=A/I$ 
and $\bar B=B/I\cdot B$. Suppose that either

a)\quad $I\to I\cdot B$ is an epimorphism with nilpotent kernel, or

b)\quad  $V$ is an $\F_p$-algebra and $I\to I\cdot B$ is invertible
up to a power of $\phi$.

Then we have :

i) conditions (a) and (b) are stable under any base change 
$A\to C$.

ii) $(A\to B)^o$ is covering (resp. bicovering) for $\tau$
if and only if $(\bar A\to\bar B)^o$ is.
\end{proposition}
\begin{proof} Suppose first that $I\to I\cdot B$ is an isomorphism;
in this case we claim that $I\cdot C\to I\cdot(C\otimes_AB)$ is
an epimorphism and  $\Ker(I\cdot C\to I\cdot(C\otimes_AB))^2=0$ for 
any $A$-algebra $C$. Indeed, since by assumption $I\simeq I\cdot B$, 
$C\otimes_AB$ acts on $C\otimes_AI$, hence $\Ker(C\to C\otimes_AB)$
annihilates $C\otimes_AI$, hence annihilates its image $I\cdot C$, 
whence the claim. If, moreover, $V$ is an $\F_p$-algebra, lemma 
\ref{lem_invert}(iv) implies that $I\cdot C\to I\cdot(C\otimes_AB)$
is invertible up to a power of $\phi$.

In the general case, consider the intermediate almost $V$-algebra 
$A_1=\bar A\times_{\bar B}B$ equipped with the ideal 
$I_1=0\times_{\bar B}(I\cdot B)$. In case (a), $I_1=I\cdot A_1$ 
and $A\to A_1$ is an epimorphism with nilpotent kernel, hence 
it remains such after any base change $A\to C$. To prove (i)
in case (a), it suffices then to consider the morphism $A_1\to B$, 
hence we can assume from start that $I\to I\cdot B$ is an isomorphism,
which is the case already dealt with. To prove (i) in case (b),
it suffices to consider the cases of $(A,I)\to(A_1,I_1)$
and $(A_1,I_1)\to(B,I\cdot B)$. The second case is treated above.
In the first case, we do not necessarily have $I_1=I\cdot A_1$
and the assertion to be checked is that, for every $A$-algebra 
$C$, the morphism $I\cdot C\to I_1\cdot(A_1\otimes_AC)$ is
invertible up to a power of $\phi$. We apply lemma 
\ref{lem_invert}(v) to the commutative diagram with exact rows:
$$\xymatrix{
0 \ar[r] & I \ar[r] \ar[d] & A \ar[r] \ar[d] & 
A/I \ar[r] \ddouble & 0 \\
0 \ar[r] & I\cdot B \ar[r] & A_1 \ar[r] & A/I \ar[r] & 0 \\
}$$
to deduce that $A\to A_1$ is invertible up to some power of 
$\phi$, hence so is $C\to A_1\otimes_AC$, which implies the
assertion.

As for (ii), we remark that  the ``only if'' part is trivial;
and we assume therefore that $(\bar A\to\bar B)^o$ is 
$\tau$-covering (resp. $\tau$-bicovering). Consider first the 
assertion for ``covering''. We need to show that $(A\to B)^o$
is of universal effective descent for $F$, where $F$ is either
one of our two fibered categories. In light of (i), this is 
reduced to the assertion that $(A\to B)^o$ is of effective 
descent for $F$. We notice that $(A\to A_1)^o$ is bicovering 
for $\tau$ (in case (a) by theorem \ref{th_liftetale} and lemma 
\ref{lem_bicover}(ii), in case (b) by proposition 
\ref{prop_invert.bicov}). As $(\bar A\to A_1/I_1)^o$ is 
an isomorphism, the assertion is reduced to the case where 
$I\to I\cdot B$ is an isomorphism. In this case, by lemma 
\ref{lem_Desc}, there is a natural equivalence: 
$\Desc(F,B/A)\stackrel{\sim}{\to}
\Desc(F,\bar B/\bar A)\times_{F_{\bar B}}F_B$.
Then the assertion follows easily from corollary 
\ref{cor_equiv.2-prod}. Finally suppose that 
$(\bar A\to\bar B)^o$ is bicovering. The foregoing 
already says that $(A\to B)^o$ is covering, so
it remains to show that $(B\otimes_AB\to B)^o$ is
also covering. The above argument again reduces to 
the case where $I\to I\cdot B$ is an isomorphism.
Then, as in the proof of lemma \ref{lem_Desc}, the induced 
morphism $I\cdot(B\otimes_AB)\to I\cdot B$ is an isomorphism
as well. Thus the assertion for ``bicovering'' is reduced to 
the assertion for ``covering''.
\end{proof}

\section{Appendix}
\subsection{}
In this appendix we have gathered a few miscellaneous results
that were found in the course of our investigation, and which
may be useful for other applications.

We need some preliminaries on simplicial objects : first of all,
a simplicial almost algebra is just an object in the category 
$s.(V^a\Alg)$. Then for a given simplicial almost algebra $A$ 
we have the category $A\Mod$ of $A$-modules : it 
consists of all simplicial almost $V$-modules $M$ such that 
$M[n]$ is an $A[n]$-module and such that the face and 
degeneracy morphisms $d_i:M[n]\to M[n-1]$ and 
$s_i:M[n]\to M[n+1]$ $(i=0,1,...,n)$ are $A[n]$-linear.  We will 
need also the derived category of $A$-modules; it is 
defined as follows.

A bit more generally, let $\cC$ be any abelian category. For an object 
$X$ of $s.\cC$ let $N(X)$ be the normalized chain complex (defined 
as in \cite{Il} (I.1.3)). By the theorem of Dold-Kan (\cite{We} th.8.4.1)  
$X\mapsto N(X)$ induces an equivalence $N:s.\cC\to\sC_\bullet(\cC)$.
Now we say that a morphism $X\to Y$ in $s.\cC$ is a 
{\em quasi-isomorphism\/} if the induced morphism $N(X)\to N(Y)$ 
is a quasi-isomorphism of chain complexes.
\medskip
 
In the following we fix a simplicial almost algebra $A$.
\begin{definition} We say that $A$ is {\em exact\/} if the 
almost algebras $A[n]$ are exact for all $n\in\N$. A morphism 
$\phi:M\to N$ of $A$-modules (or $A$-algebras) 
is a {\em quasi-isomorphism\/} if the morphism $\phi$ of underlying 
simplicial almost $V$-modules is a quasi-isomorphism. We define the 
category $\sD_\bullet(A)$ (resp. the category $\sD_\bullet(A\Alg)$) 
as the localization of the category $A\Mod$ (resp. $A\Alg$) with 
respect to the class of quasi-isomorphisms. 
\end{definition}

As usual, the morphisms in $\sD_\bullet(A)$ can be computed via
a calculus of fraction on the category $\Hot_\bullet(A)$
of simplicial complexes up to homotopy. Moreover, if $A_1$
and $A_2$ are two simplicial almost algebras, then the 
``extension of scalars'' functors define equivalences of
categories
$$\begin{array}{l}
\sD_\bullet(A_1\times A_2)\stackrel{\sim}{\longrightarrow}
\sD_\bullet(A_1)\times\sD_\bullet(A_2) \\
\sD_\bullet(A_1\times A_2\Alg)\stackrel{\sim}{\longrightarrow}
\sD_\bullet(A_1\Alg)\times\sD_\bullet(A_2\Alg).
\end{array}$$
\begin{proposition}\label{prop_preserve} (i) The functor on 
$A$-algebras given by $B\mapsto(s.V^a\times B)_{!!}$ preserves 
quasi-isomorphisms and therefore induces a functor 
$\sD_\bullet(A\Alg)\to\sD_\bullet((s.V^a\times A)_{!!}\Alg)$.

(ii) The localisation functor $R\mapsto R^a$ followed by ``extension
of scalars'' via $s.V^a\times A\to A$ induces a functor
$\sD_\bullet((s.V^a\times A)_{!!}\Alg)\to\sD_\bullet(A\Alg)$
and the composition of this and the above functor is naturally
isomorphic to the identity functor on $\sD_\bullet(A\Alg)$.
\end{proposition}
\begin{proof} (i) : let $B\to C$ be a quasi-isomorphism 
of $A$-algebras. Clearly the induced morphism 
$s.V^a\times B\to s.V^a\times C$ is still a quasi-isomorphism 
of $V$-algebras. But by remark \ref{rem_exact.alg}, 
$s.V^a\times B$ and $s.V^a\times C$ are exact simplicial 
almost $V$-algebras; moreover, it follows from corollary 
\ref{cor_left.exact} that $(s.V^a\times B)_!\to(s.V^a\times C)_!$ 
is a quasi-isomorphism of $V$-modules. Then the claim 
follows easily from the exactness of the sequence 
\eqref{eq_left.adj.algebras}. Now (ii) is clear.
\end{proof}

\begin{remark} In case $\fm$ is flat, then all
$A$-algebras are exact, and the same argument shows that
the functor $B\mapsto B_{!!}$ induces a functor 
$\sD_\bullet(A\Alg)\to\sD_\bullet(A_{!!}\Alg)$.
In this case, composition with localisation is naturally
isomorphic to the identity functor on $\sD_\bullet(A\Alg)$.
\end{remark}

\begin{proposition}\label{prop_Gabber} 
Let $f:R\to S$ be a map of $V$-algebras such that $f^a:R^a\to S^a$ 
is an isomorphism. Then $\L_{S/R}^a\simeq 0$ in $\sD_\bullet(s.S^a)$.
\end{proposition}
\begin{proof} We show by induction on $q$ that 

\medskip

\noindent{\bf VAN}($q;S/R$){\hskip 2.5cm} $H_q(\L_{S/R}^a)=0.$
\medskip

For $q=0$ the claim follows immediately from \cite{Il} (II.1.2.4.2).
Therefore suppose that $q>0$ and that {\bf VAN}($j;D/C$) is known for
all almost isomorphisms of $V$-algebras $C\to D$ and all $j<q$.
Let $\bar R=f(R)$. Then by transitivity (\cite{Il} (II.2.1.2)) 
we have a distinguished triangle in $\sD_\bullet(s.S^a)$
$$\diagram
(S\otimes_{\bar R}\L_{\bar R/R})^a \rto^-u & 
\L^a_{S/R} \rto^v & \L^a_{S/\bar R} \rto & 
\sigma(S\otimes_{\bar R}\L_{\bar R/R})^a.
\enddiagram$$
We deduce that {\bf VAN}($q;\bar R/R$) and {\bf VAN}($q;S/\bar R$) imply
{\bf VAN}($q;S/R$), thus we can assume that $f$ is either injective or
surjective. Let $S_\bullet\to S$ be the simplicial $V$-algebra
augmented over $S$ defined by $S_\bullet=P_V(S)$. It is a simplicial
resolution of $S$ by free $V$-algebras, in particular the augmentation
is a quasi-isomorphism of simplicial $V$-algebras. Set 
$R_\bullet=S_\bullet\times_SR$. This is a simplicial $V$-algebra
augmented over $R$ via a quasi-isomorphism. Moreover, the induced
morphisms $R[n]^a\to S[n]^a$ are isomorphisms. By \cite{Il} (II.1.2.6.2)
there is a quasi-isomorphism 
$\L_{S/R}\simeq\L^\Delta_{S_\bullet/R_\bullet}$.
On the other hand we have a spectral sequence
$$E^1_{ij}=H_j(\L_{S[i]/R[i]})\Rightarrow 
H_{i+j}(L^\Delta_{S_\bullet/R_\bullet}).$$
It follows easily that 
{\bf VAN}$(j;S[i]/R[i])$ for all $i\ge 0,j\le q$
implies {\bf VAN}$(q;S/R)$.
Therefore we are reduced to the case where
$S$ is a free $V$-algebra and $f$ is either injective or surjective.
We examine separately these two cases. If $f:R\to V[T]$ is surjective,
then we can find a right inverse $s:V[T]\to R$ for $f$. By applying
transitivity to the sequence $V[T]\to R\to V[T]$ we get a distinguished
triangle
$$(V[T]\otimes_R\L_{R/V[T]})^a \stackrel{u}{\to}\L^a_{V[T]/V[T]}
\stackrel{v}{\to}\L^a_{V[T]/R}\to\sigma(V[T]\otimes_R\L_{R/V[T]})^a.$$
Since $\L^a_{V[T]/V[T]}\simeq 0$ there follows an isomorphism :
$H_q(\L_{V[T]/R})^a\simeq H_{q-1}(V[T]\otimes_R\L_{R/V[T]})^a$.
Furthermore, since $f^a$ is an isomorphism, $s^a$ is an isomorphism as well,
hence by induction (and by a spectral sequence of the type \cite{Il} 
(I.3.3.3.2)) $H_{q-1}(V[T]\otimes_R\L_{R/V[T]})^a\simeq 0$.
The claim follows in this case.

Finally suppose that $f:R\to V[T]$ is injective. Write $V[T]=\Sym(F)$,
for a free $V$-module $F$ and set $\tilde F=\tilde\fm\otimes_VF$; since 
$f^a$ is an isomorphism, 
$\Img(\Sym(\tilde F)\to\Sym(F))\subset R$.
We apply transitivity to the sequence 
$\Sym(\tilde F)\to R\to\Sym(F)$. By arguing as above 
we are reduced to showing that 
$\L^a_{\Sym(F)/\Sym(\tilde F)}\simeq 0.$ We know that 
$H_0(\L^a_{\Sym(F)/\Sym(\tilde F)})\simeq 0$ and we will 
show that  
$H_q(\L^a_{\Sym(F)/\Sym(\tilde F)})\simeq 0$ for $q>0$. 
To this purpose we apply transitivity to the sequence 
$V\to\Sym(\tilde F)\to\Sym(F)$. As $F$ and 
$\tilde F$ are flat $V$-modules, \cite{Il} (II.1.2.4.4) 
yields $H_q(\L_{\Sym(F)/V})\simeq
H_q(\L_{\Sym(\tilde F)/V})\simeq 0$ for $q>0$ and 
$H_0(\L_{\Sym(\tilde F)/V})$ is a flat 
$\Sym(\tilde F)$-module.
In particular 
$H_j(\Sym(F)\otimes_{\Sym(\tilde F)}\L_{\Sym(\tilde F)/V})
\simeq 0$ for all $j>0$. 
Consequently $H_{j+1}(\L_{\Sym(F)/\Sym(\tilde F)})\simeq 0$ 
for all $j>0$ and 
$H_1(\L_{\Sym(F)/\Sym(\tilde F)})\simeq
\Ker(\Sym(F)\otimes_{\Sym(\tilde F)}
\Omega_{\Sym(\tilde F)/V}\to\Omega_{\Sym(F)/V})$.
The latter module is easily seen to be almost zero.
\end{proof}

\begin{theorem}\label{th_Gabber} Let $\phi:R\to S$ be a map 
of simplicial $V$-algebras inducing an isomorphism 
$R^a\stackrel{\sim}{\to}S^a$ in $\sD_\bullet(R^a)$. 
Then $(\L^\Delta_{S/R})^a\simeq 0$ in $\sD_\bullet(S^a)$.
\end{theorem}
\begin{proof} Apply the base change theorem (\cite{Il} II.2.2.1)
to the (flat) projections of $s.V\times R$ onto $R$ and respectively
$s.V$ to deduce that the natural map 
$\L^\Delta_{s.V\times S/s.V\times R}\to
\L^\Delta_{S/R}\oplus\L^\Delta_{s.V/s.V}\to\L^\Delta_{S/R}$ is a 
quasi-isomorphism in $\sD_\bullet(s.V\times S)$. By proposition 
\ref{prop_preserve} the induced morphism 
$(s.V\times R)^a_{!!}\to(s.V\times S)^a_{!!}$ is still a 
quasi-isomorphism. There are spectral sequences
$$\begin{array}{c}
E^1_{ij}=H_j(\L_{(V\times R[i])/(V\times R[i])_{!!}^a})\Rightarrow 
H_{i+j}(\L^\Delta_{(s.V\times R)/(s.V\times R)_{!!}^a}) \\
F^1_{ij}=H_j(\L_{(V\times S[i])/(V\times S[i])_{!!}^a})\Rightarrow 
H_{i+j}(\L^\Delta_{(s.V\times S)/(s.V\times S)_{!!}^a}).
\end{array}$$
On the other hand, by proposition \ref{prop_Gabber} we
have
$\L^a_{(V\times R[i])/(V\times R[i])_{!!}^a}\simeq 0
\simeq\L^a_{(V\times S[i])/(V\times S[i])_{!!}^a}$
for all $i\in\N$. Then the theorem follows directly from 
\cite{Il} (II.1.2.6.2(b)) and transitivity. 
\end{proof}

\begin{proposition}\label{prop_quasi.exact}
Let $A\to B$ be a morphism of exact almost $V$-algebras. 
Then the natural map 
$\tilde\fm\otimes_V\L_{B_{!!}/A_{!!}}\to\L_{B_{!!}/A_{!!}}$ 
is a quasi-isomorphism.
\end{proposition}
\begin{proof} By transitivity we may assume $A=V^a$. Let 
$P_\bullet=P_V(B_{!!})$ be the standard resolution of $B_{!!}$
(see \cite{Il} II.1.2.1). Each $P[n]^a$ contains $V$ as a 
direct summand, hence it is exact, so that we have an exact 
sequence of simplicial $V$-modules 
$0\to s.\tilde\fm\to s.V\oplus(P^a_\bullet)_!\to
(P^a_\bullet)_{!!}\to 0$.
The augmentation $(P^a_\bullet)_!\to(B^a_{!!})_!\simeq B_!$
is a quasi-isomorphism and we deduce that 
$(P^a_\bullet)_{!!}\to B_{!!}$ is a quasi-isomorphism; hence 
$(P^a_\bullet)_{!!}\to P_\bullet$ is a quasi-isomorphism as 
well. We have $P[n]\simeq\Sym(F_n)$ for a free $V$-module 
$F_n$ and the map $(P[n]^a)_{!!}\to P[n]$ is identified with 
$\Sym(\tilde\fm\otimes_VF_n)\to\Sym(F_n)$, whence 
$\Omega_{P[n]^a_{!!}/V}\otimes_{P[n]^a_{!!}}P[n]\to
\Omega_{P[n]/V}$ is identified with 
$\tilde\fm\otimes_V\Omega_{P[n]/V}\to\Omega_{P[n]/V}$.
By \cite{Il} (II.1.2.6.2) the map 
$\L^\Delta_{(P^a_\bullet)_{!!}/V}\to\L^\Delta_{P_\bullet/V}$
is a quasi-isomorphism. In view of \cite{Il} (II.1.2.4.4) we 
derive that $\Omega_{(P_\bullet^a)_{!!}/V}\to\Omega_{P_\bullet/V}$
is a quasi-isomorphism, {\em i.e.\/} 
$\tilde\fm\otimes_V\Omega_{P_\bullet/V}\to\Omega_{P_\bullet/V}$
is a quasi-isomorphism. Since $\tilde\fm$ is flat and 
$\Omega_{P_\bullet/V}\to
\Omega_{P_\bullet/V}\otimes_{P_\bullet}B_{!!}=\L_{B_{!!}/V}$
is a quasi-isomorphism, we get the desired conclusion.
\end{proof}

In view of proposition \ref{prop_Gabber} we have  
$\L^a_{(V^a\times A)_{!!}/V\times A_{!!}}\simeq 0$
in $\sD_\bullet(V^a\times A)$. By this, transitivity 
and localisation (\cite{Il} II.2.3.1.1) we derive that 
$\L^a_{B/A}\to\L^a_{B_{!!}/A_{!!}}$ is a quasi-isomorphism 
for all $A$-algebras $B$. If $A$ and $B$ are exact 
({\em e.g.\/} if $\fm$ is flat), we conclude from proposition 
\ref{prop_quasi.exact} that the natural map 
$\L_{B/A}\to\L_{B_{!!}/A_{!!}}$ is a quasi-isomorphism.

Finally we want to discuss left derived functors of (the almost version
of) some notable non-additive functors that play a role in deformation
theory. Let $R$ be a simplicial $V$-algebra. Then we have an obvious 
functor $G:\sD_\bullet(R)\to\sD_\bullet(R^a)$ obtained by applying 
dimension-wise the localisation functor. Let $\Sigma$ be the 
multiplicative set of morphisms of $\sD_\bullet(R)$ that induce almost 
isomorphisms on the cohomology modules. An argument as in section 
\ref{sec_homol} shows that $G$ induces an equivalence of categories
$\Sigma^{-1}\sD_\bullet(R)\to\sD_\bullet(R^a)$.

Now let $R$ be a $V$-algebra and $\cF_p$ one of the functors 
$\otimes^p$, $\Lambda^p$, $\Sym^p$, $\Gamma^p$ defined in 
\cite{Il} (I.4.2.2.6).
\begin{lemma}\label{lem_cF} Let $\phi:M\to N$ be an almost 
isomorphism of $R$-modules. Then $\cF_p(\phi):\cF_p(M)\to\cF_p(N)$ 
is an almost isomorphism.
\end{lemma}
\begin{proof} Let $\psi:\tilde\fm\otimes_VN\to M$ be the 
map corresponding to $(\phi^a)^{-1}$ under the bijection
\eqref{eq_alm.morph}. By inspection, the compositions 
$\phi\circ\psi:\tilde\fm\otimes_VN\to N$ and 
$\psi\circ(\one_{\tilde\fm}\otimes\phi):\tilde\fm\otimes_VM\to M$ 
are induced by scalar multiplication. Pick any $s\in\fm$ and lift 
it to an element $\tilde s\in\tilde\fm$; define $\psi_s:N\to M$
by $n\mapsto\psi(\tilde s\otimes n)$ for all $n\in N$.
Then $\phi\circ\psi_s=s\cdot\one_N$ and 
$\psi_s\circ\phi=s\cdot\one_M$. This easily implies that 
$s^p$ annihilates $\Ker\cF_p(\phi)$ and $\Coker\cF_p(\phi)$.
In light of proposition \ref{prop_less.obv}(ii), the claim 
follows. 
\end{proof}
Let $B$ be an almost $V$-algebra. We define a functor 
$\cF_p^a$ on $B\Mod$ by $M\mapsto(\cF_p(M_!))^a$, where 
$M_!$ is viewed as a $B_{!!}$-module or a $B_*$-module 
(to show that these choices define the same functor it 
suffices to observe that $B_*\otimes_{B_{!!}}N\simeq N$ 
for all $B_*$-modules $N$ such that $N=\fm\cdot N$).
For all $p>0$ we have diagrams :
\set\begin{equation}\label{eq_derotimes}{
\diagram
R\Mod \ar[r]^{\cF_p} \ar@<.5ex>[d] & R\Mod \ar@<.5ex>[d] \\
R^a\Mod \ar[r]^{\cF_p^a} \ar@<.5ex>[u] & 
R^a\Mod \ar@<.5ex>[u]
\enddiagram}\end{equation}
where the downward arrows are localisation and the upward
arrows are the functors $M\mapsto M_!$.
Lemma \ref{lem_cF} implies that the downward arrows in the
diagram commute (up to a natural isomorphism) with the 
horizontal ones. It will follow from the following proposition 
\ref{prop_upward} that the diagram commutes also going upward. 

For any $V$-module $N$ we have an exact sequence 
$\Gamma^2N\to\otimes^2N\to\Lambda^2N\to 0$. As observed in
the proof of proposition \ref{prop_less.obv}, the symmetric
group $S_2$ acts trivially on $\otimes^2\tilde\fm$ and
$\Gamma^2\tilde\fm\simeq\otimes^2\tilde\fm$, so 
$\Lambda^2\tilde\fm=0$. Also we have natural isomorphisms
$\Gamma^p\tilde\fm\simeq\tilde\fm$ for all $p>0$.
\begin{proposition}\label{prop_upward}
Let $R$ be a commutative ring and $L$ a flat $R$-module 
with $\Lambda^2L=0$. Then for $p>0$ and for all $R$-modules
$N$ we have natural isomorphisms
$$\Gamma^p(L)\otimes_R\cF_p(N)
\stackrel{\sim}{\to}{\cF_p}(L\otimes_RN).$$
\end{proposition}
\begin{proof} Fix an element $x\in\cF_p(N)$. 
For each $R$-algebra $R'$ and each element $l\in R'\otimes_RL$
we get a map $\phi_l:R'\otimes_RN\to R'\otimes_RL\otimes_RN$
by $y\mapsto l\otimes y$, hence a map 
$\cF_p(\phi_l):R'\otimes_R\cF_p(N)\simeq\cF_p(R'\otimes_RN)\to
\cF_p(R'\otimes_RL\otimes_RN)\simeq 
R'\otimes_R\cF_p(L\otimes_RN)$. For varying $l$ we obtain
a map of sets $\psi_{R',x}:R'\otimes_RL\to 
R'\otimes_R\cF_p(L\otimes_RN)$ : 
$l\mapsto\cF_p(\phi_l)(1\otimes x)$. According to the 
terminology of \cite{Ro}, the system of maps $\psi_{R',x}$ 
for $R'$ ranging over all $R$-algebras forms a homogeneous 
polynomial law of degree $p$ from $L$ to $\cF_p(L\otimes_RN)$, 
so it factors through the universal homogeneous degree $p$ 
polynomial law $\gamma_p:L\to\Gamma^p(L)$ . The resulting 
$R$-linear map $\bar\psi_x:\Gamma^p(L)\to\cF_p(L\otimes_RN)$ 
depends $R$-linearly on $x$, hence we derive an $R$-linear 
map $\psi:\Gamma^p(L)\otimes_R\cF_p(N)\to{\cF_p}(L\otimes_RN)$.
Next notice that by hypothesis $S_2$ acts trivially on $\otimes^2L$ 
so $S_p$ acts trivially on 
$\otimes^pL$ and we get an isomorphism
$\beta:\Gamma^p(L)\stackrel{\sim}{\longrightarrow}\otimes^pL$.
We deduce a natural map 
$(\otimes^pL)\otimes_R\cF_p(N)\to\cF_p(L\otimes_RN)$.
Now, in order to prove the proposition for the case 
$\cF_p=\otimes^p$, it suffices to show that this last 
map is just the natural isomorphism that ``reorders
the factors''. Indeed, let $x_1,...,x_n\in L$ and 
$q=(q_1,...,q_n)\in\N^n$ such that $|q|=\sum_iq_i=p$;
then $\beta$ sends the generator 
$x_1^{[q_1]}\cdot...\cdot x_n^{[q_n]}$ to 
${p\choose q_1,...,q_n}\cdot 
x_1^{\otimes q_1}\otimes...\otimes x_n^{\otimes q_n}$.
On the other hand, pick any $y\in\otimes^pN$ and let 
$R[T]=R[T_1,...,T_r]$ be the polynomial $R$-algebra in $n$ 
variables; write 
$(T_1\otimes x_1+...+T_n\otimes x_n)^{\otimes p}\otimes y=
\psi_{R[T],y}(T_1\otimes x_1+...+T_n\otimes x_n)=
\sum_{r\in\N^n}T^r\otimes w_r$ with $w_r\in\otimes^p(L\otimes_RN)$.
Then $\psi((x_1^{[q_1]}\cdot...\cdot x_n^{[q_n]})\otimes y)=w_q$
(see \cite{Ro} pp.266-267)
and the claim follows easily. Next notice 
that $\Gamma^p(L)$ is flat, so that tensoring with $\Gamma^p(L)$ 
commutes with taking coinvariants (resp. invariants) under the 
action of the symmetric group; this implies the assertion for 
$\cF_p=\Sym^p$ (resp. $\cF_p=\text{TS}^p$). 
To deal with $\cF_p=\Lambda^p$ recall that for any $V$-module 
$M$ and $p>0$ we have the antisymmetrization  operator 
$a_M=\sum_{\sigma\in S_p}\text{sgn}(\sigma)\cdot\sigma:
\otimes^pM\to\otimes^pM$ and a surjection $\Lambda^p(M)\to\Img(a_M)$ 
which is an isomorphism for $M$ free, hence for $M$ flat.
The result for $\cF_p=\otimes^p$ (and again the flatness of 
$\Gamma^p(L)$) then gives 
$\Gamma^p(L)\otimes\Img(a_N)\simeq\Img(a_{L\otimes_RN})$,
hence the assertion for $\cF_p=\Lambda^p$ and $N$ flat.
For general $N$ let $F_1\stackrel{\partial}{\to}
F_0\stackrel{\eps}{\to}N\to 0$ be a presentation with $F_i$
free. Define $j_0,j_1:F_0\oplus F_1\to F_0$ by 
$j_0(x,y)=x+\partial(y)$ and $j_1(x,y)=x$. By functoriality
we derive an exact sequence 
$$\xymatrix{
\Lambda^p(F_0\oplus F_1) \ar@<.5ex>[r] \ar@<-.5ex>[r] & 
\Lambda^p(F_0) \ar[r] & \Lambda^p(N) \ar[r] & 0
}$$
which reduces the assertion to the flat case. For $\cF_p=\Gamma^p$
the same reduction argument works as well (cf. \cite{Ro} p.284)
and for flat modules the assertion for $\Gamma^p$ follows
from the corresponding assertion for $\text{TS}^p$.
\end{proof}
\begin{lemma} Let $A$ be a simplicial almost algebra, $L,E$ 
and $F$ three $A$-modules, $f:E\to F$ a quasi-isomorphism. 
If $L$ is flat or $E,F$ are flat, then 
$L\otimes_Af:L\otimes_AE\to L\otimes_AF$ is a quasi-isomorphism.
\end{lemma}
\begin{proof} It is deduced directly from \cite{Il} (I.3.3.2.1)
by applying $M\mapsto M_!$.
\end{proof}
As usual, this allows one to show that 
$\otimes:\Hot_\bullet(A)\times\Hot_\bullet(A)\to
\Hot_\bullet(A)$ admits a left derived functor 
$\derotimes:\sD_\bullet(A)\times\sD_\bullet(A)\to\sD_\bullet(A)$.
If $R$ is a simplicial $V$-algebra then we have essentially
commutative diagrams
$$\xymatrix{
\sD_\bullet(R)\times\sD_\bullet(R) \ar[r]^-{\derotimes} \ar@<.5ex>[d] &
\sD_\bullet(R) \ar@<.5ex>[d] \\
\sD_\bullet(R^a)\times\sD_\bullet(R^a)\ar[r]^-{\derotimes}\ar@<.5ex>[u]&
\sD_\bullet(R^a) \ar@<.5ex>[u]
}$$
where again the downward (resp. upward) functors are induced by 
localisation (resp. by $M\mapsto M_!$).

We mention the derived functors of the non-additive functor $\cF_p$ 
defined above in the simplest case of modules over a constant
simplicial ring. Let $A$ be a (commutative) almost algebra.
\begin{lemma} If $u:X\to Y$ is a quasi-isomorphism of flat 
$s.A$-modules then $\cF_p^a(u):\cF_p^a(X)\to\cF_p^a(Y)$ is 
a quasi-isomorphism.
\end{lemma}
\begin{proof} This is deduced from \cite{Il} (I.4.2.2.1) applied
to $N(X_!)\to N(Y_!)$ which is a quasi-isomorphism of chain complexes
of flat $A_{!!}$-modules. We note that {\em loc. cit.\/} deals 
with a more general mixed simplicial construction of $\cF_p$ which
applies to bounded above complexes, but one can check that it
reduces to the simplicial definition for complexes in 
$\cC_\bullet(A_{!!})$.
\end{proof}
Using the lemma one can construct 
$L\cF_p^a:\sD_\bullet(s.A)\to\sD_\bullet(s.A)$. If $R$ is a
$V$-algebra we have the derived category version of the 
essentially commutative squares \eqref{eq_derotimes}, 
relating $L\cF_p:\sD_\bullet(s.R)\to\sD_\bullet(s.R)$
and $L\cF_p^a:\sD_\bullet(s.R^a)\to\sD_\bullet(s.R^a)$.

\begin{acknowledgement}
The second author is very much indebted to Gerd Faltings 
for many patient explanations on the method of almost 
\'etale extensions.
Next he would like to acknowledge several interesting 
discussions with Ioannis Emmanouil. He is also much obliged 
to Pierre Deligne, for a useful list of critical remarks.
Finally, he owes a special thank to Roberto Ferretti, who has read 
the first tentative versions of this article, has corrected many 
slips and has made several valuable suggestions.
\end{acknowledgement}

\vskip 1cm

\noindent
\begin{tabular}{lll}
 Ofer Gabber              &           & Lorenzo Ramero \\
 I.H.E.S.                 &           & Universit\'e de Bordeaux I \\
 Le Bois-Marie            &           & Laboratoire de Math\'ematiques Pures \\
 35, route de Chartres    &           & 351, cours de la Liberation \\
 F-91440 Bures-sur-Yvette &           & F-33405 Talence Cedex \\
 {\footnotesize\em E-mail address:} \texttt{gabber@ihes.fr} 
 & \qquad &
 {\footnotesize\em E-mail address:} \texttt{ramero@math.u-bordeaux.fr} \\
 & \qquad &
 {\footnotesize\em web:} \texttt{http://www.math.u-bordeaux.fr/$\sim$ramero}
\end{tabular}

\end{document}